\documentclass[a4paper,10pt]{article}
\usepackage{graphicx,multirow}
\usepackage{amssymb,amsmath,amsthm}
\usepackage{lscape,lipsum,microtype}
\usepackage[margin=2.54cm]{geometry}  
\usepackage{hyperref}

\usepackage{ifluatex}
\ifluatex
 \usepackage{unicode-math}
\fi

\usepackage{latexsym}
\usepackage{amsmath}
\usepackage{amsfonts}
\usepackage{amssymb}
\usepackage{amscd}

\renewcommand{\d}{\delta }

\renewcommand{\l }{\lambda }

\newcommand{\n }{\nabla }

\newcommand{\intbar}{\mathop{\int\makebox(-13.5,0){\rule[4pt]{.7em}{0.3pt}}%
\kern-6pt}\nolimits}

\newcommand{\be}{\begin{equation}}
\newcommand{\ee}{\end{equation}}
\newcommand{\bes}{\begin{equation*}}
\newcommand{\ees}{\end{equation*}}
\newcommand{\ba}{\begin{eqnarray}}
\newcommand{\ea}{\end{eqnarray}}
\newcommand{\bas}{\begin{eqnarray*}}
\newcommand{\eas}{\end{eqnarray*}}
\newenvironment{pf}{\noindent{\sc Proof}.\enspace}{\rule{2mm}{2mm}\medskip}
\newenvironment{pfn}{\noindent{\sc Proof}}{\rule{2mm}{2mm}\medskip}

\newcommand{\R}{\mathbb{R}}

\newcommand{\Z}{\mathbb{Z}}

\newcommand{\N}{\mathbb{N}}

\author{Mohammed ALDAWOOD}

\date{}

\title{\bf A Bahri-Brezis type problem on asymptotically hyperbolic manifolds. }
\begin{document}

\newtheorem{lem}{Lemma}[section]
\newtheorem{pro}[lem]{Proposition}
\newtheorem{thm}[lem]{Theorem}
\newtheorem{rem}[lem]{Remark}
\newtheorem{cor}[lem]{Corollary}
\newtheorem{df}[lem]{Definition}

\maketitle

\begin{center}
{\small

\noindent
Department of Mathematics Howard University \\  204 Academic Service Building B\\
Washington, D.C. 20059, USA.
}

\end{center}

\footnotetext[1]{E-mail addresses: mohammed.aldawood@bison.howard.edu}

\
\

\begin{center}
{\bf Abstract}
\end{center}
In this paper, we study a Bahri-Brezis type problem on a compact $3$-dimensional  asymptotically hyperbolic manifold. Using the celebrated Algebraic Topological argument of Bahri-Coron\cite{bc}, we show the existence of at least one solution under the assumption that the
corresponding degenerate boundary value problem has a positive first eigenvalue and a positive Green’s function.
 \begin{center}


\bigskip

\centerline{\bf AMS subject classification: 53C21, 35C60, 58J60, 55N10.}

\end{center}

\section{Introduction and statement of the results}
Given \;$X=X^{n+1}$\; a smooth manifold with boundary \;$M=M^{n}$,\; closure \;$\overline{X}:=X\cup M$,\; and \;$n\geq 2$\; we say that \;$\rho$\; is a defining function of the boundary \;$M$\; in \;$\overline{X}$,\; if
$$\rho>0\;\;\;\text{in}\;\;\;X,\;\;\;\rho=0\;\;\;\text{on}\;\;\;M\;\;\;\text{and}\;\;\;d\rho\neq0\;\;\;\text{on}\;\;\;M.$$
A Riemannian metric \;$g^{+}$\; on \;$X$\; is said to be conformally compact, if for some defining function \;$\rho$\; the
Riemannian metric
\begin{equation}\label{metric.g}
    g:=\rho^{2}g^{+}
\end{equation}
extends to \;$\overline{X}$\; so that \;$(\overline{X},g)$\; is a compact Riemannian manifold with boundary \;$M$\; and interior \;$X$.\; Clearly this induces a conformal class of Riemannian metrics
$$[h]=[g|_{TM}]$$
on \;$M$,\; where \;$TM$\; denotes the tangent bundle of \;$M$,\; when the defining functions \;$\rho$\; vary and the resulting conformal manifold \;$(M,[h])$\; is called conformal infinity of \;$(X,g^{+})$.\; Moreover a Riemannian metric \;$g^{+}$\; in \;$X$\; is said to be asymptotically hyperbolic, if it is conformally compact and its sectional
curvature tends to \;$-1$\; as one approaches the conformal infinity of \;$(X,g^{+})$,\; which is equivalent to
$$|d\rho|_{g}=1$$
on \;$M$\;, see \cite{mr}, and in such a case \;$(X,g^{+})$\; is called an asymptotically hyperbolic manifold. 
\vspace{4pt}

\noindent
 For every asymptotically hyperbolic manifold \;$(X,g^{+})$\; and every choice of the representative \;$h$\; of its conformal infinity \;$(M, [h])$,\; there exists a geodesic defining function \;$y$\; of \;$M$\; in \;$\overline{X}$\; such that
in a tubular neighborhood of \;$M$\; in \;$\overline{X}$,\; the Riemannian metric \;$g^{+}$\; takes the following normal form
\begin{equation}\label{metric.g+}
    g^{+}=\frac{dy^{2}+h_{y}}{y^{2}},
\end{equation}
where \;$h_y$\; is a family of Riemannian metrics on \;$M$\; satisfying \;$h_{0}=h$\; and \;$y$\; is the unique such a one in a tubular neighborhood of \;$M$.\;We recall \;$g=y^{2}g^{+}$\;. 
\vspace{4pt}

\noindent
\noindent
In this paper, we study the Boundary Value Problem (BVP) with critical Sobolev nonlinearity on the boundary defined by \begin{equation}\label{BVP1}
\left\{
\begin{split}
-div_{g}(y^{1-2\gamma}\nabla_{g}u)+qu&=0&\;\;\;\;\text{in}\;\;\;\;&X,\\
-d^{*}_{\gamma}\lim_{y\to 0}y^{1-2\gamma}\frac{\partial u}{\partial y}&=u^{\frac{n+2\gamma}{n-2\gamma}}&\;\;\;\;\text{on}\;\;\;\;&M,\\
u&>0\;\;\;\;&\text{on}\;\;\;\;&\overline{X},
\end{split}
\right.
\end{equation}
where \;$\left(X,g^{+}\right)$\; is an asymptotically hyperbolic manifold with conformal infinity \;$(M,[h])$\; of dimensional \;$n\geq2$,\;\;$y$\; is the unique geodesic defining function associated to \;$h$\; verifying \eqref{metric.g+}, \;$q$\; is a bounded smooth potential defined on \;$X$,\;\;$\frac{\partial}{\partial y}$\; is the derivative with respect to \;$y$,\;\;$d^{*}_{\gamma}$\; is as in \eqref{d*}, \;$div_{g}$\; is the divergent with respect to \;$g$\;,\;$\nabla_{g}$\; is  covariant derivative with respect to \;$g$\;,  and \;$\gamma\in(0,1).$\;\\
We specialize to the case \;$dim\;M=n=2$.\; Hence, the BVP of interest becomes
\begin{equation}\label{BVP2}
\left\{
\begin{split}
-div_{g}(y^{1-2\gamma}\nabla_{g}u)+qu&=0&\;\;\;\;\text{in}\;\;\;\;&X,\\
-d^{*}_{\gamma}\lim_{y\to 0}y^{1-2\gamma}\frac{\partial u}{\partial y}&=u^{\frac{1+\gamma}{1-\gamma}}&\;\;\;\;\text{on}\;\;\;\;&M,\\
u&>0&\;\;\;\;\text{on}\;\;\;\;&\overline{X}.
\end{split}
\right.
\end{equation}
\noindent
We will also assume that the Green's function \;$G$\; defined by \eqref{eqgreen} is strictly positive and hence \eqref{BVP2} satisfies the strong maximum principle in the following sense: if \;$u$\; satisfies
\begin{equation}\label{pq2}
\left\{
\begin{split}
-div_{g}(y^{1-2\gamma}\nabla_{g}u)+qu&=0\;\;\;\;&\text{in}\;\;\;\;X,\\
-d^{*}_{\gamma}\lim_{y\to 0}y^{1-2\gamma}\frac{\partial u}{\partial y}&\geq0\;\;\;\;&\text{on}\;\;\;\;M,\\
u&\geq0\;\;\;\;&\text{on}\;\;\;\;M,\\
u&\neq0\;\;\;\;&\text{on}\;\;\;\;M,
\end{split}
\right.
\end{equation}
then 
\begin{equation}\label{pq3}
    u>0\;\;\;\;\;\;\text{on}\;\;\;\;\;\overline{X}.
\end{equation}
The BVP \eqref{BVP2} is equivalent to a nonlocal Sobolev critical problem. Indeed, letting \;$P^{\gamma}_q$\; be the Dirichlet-to-Neumann map type operator defined by \(\ P^{\gamma}_q:{C^{\infty}(M)} \to {C^{\infty}(M)}\), \;$u\to -d^{*}_{\gamma}\lim_{y\to 0}y^{1-2\gamma}\frac{\partial U_{q}}{\partial y}$\; where \;$U_q$\; is the unique solution of
\begin{equation}\label{BVP3}
\left\{
\begin{split}
-div_{g}(y^{1-2\gamma}\nabla_{g}U_q)+qU_q&=0&\;\;\;\;\text{in}\;\;\;\;&X,\\
U_q&=u&\;\;\;\;\text{on}\;\;\;\;& M,
\end{split}
\right.
\end{equation}
then, clearly \;$u$\; is solution of BVP \eqref{BVP2} implies 
\begin{equation}\label{PQ}
\left\{
\begin{split}
 P^{\gamma}_qu&=u^{\frac{1+\gamma}{1-\gamma}}&\;\;\;\;\text{on}\;\;\;\;& M,\\
u&>0&\;\;\;\;\text{on}\;\;\;\;& M,
\end{split}
\right.
\end{equation}
and \;$u$\; is solution of \eqref{PQ} implies that \;$U_q$\; defined by \eqref{BVP3} is solution of BVP \eqref{BVP2} thanks to strong maximum principle.
\vspace{4pt}

\noindent
It is easy to see that \;$\lambda_{1}(P^{\gamma}_q)$\; the first eigenvalue of \;$P^{\gamma}_q$\; is the same as the first eigenvalue of the eigenvalue problem
\begin{equation}\label{BVP4}
    \left\{
\begin{split}
-div_{g}(y^{1-2\gamma}\nabla_{g}u)+qu&=0&\;\;\;\;\text{in}\;\;\;\;&X,\\
-d^{*}_{\gamma}\lim_{y\to 0}y^{1-2\gamma}u&=\lambda u&\;\;\;\;\text{on}\;\;\;\;& M,
\end{split}
\right.
\end{equation}
 and \;$\lambda \in \R$\;.
\vspace{4pt}

\noindent
Moreover, it is easy to see that problem \eqref{BVP2} satisfies the strong maximum principle in the sense \eqref{pq2} is equivalent to \;$P^{\gamma}_{q}$\; satisfies the strong maximum principle in the  sense that: if \;$u$\; satisfies
\begin{equation}\label{pq3}
\left\{
\begin{split}
P^{\gamma}_{q} u&\geq0\;\;\;\;&\text{on}\;\;\;\;M,\\
u&\geq0\;\;\;\;&\text{on}\;\;\;\;M,\\
u&\neq0\;\;\;\;&\text{on}\;\;\;\;M,
\end{split}
\right.
\end{equation}
then 
\begin{equation}\label{pq4}
    U_{q}>0\;\;\;\;\;\;\text{on}\;\;\;\;\;\overline{X},\;\;\;\text{where}\;\;\;\;\;U_{q}\;\;\;\;\;\text{as in }\;\;\;\;\;\eqref{BVP3}.
\end{equation}
\vspace{4pt}

\noindent
As in \cite{aldawood} and \cite{aldawood1}, it is easy that a necessary condition of the existence of solution to \eqref{BVP2} is that the first eigenvalue \;$\lambda_{1}(P^{\gamma}_q)>0$\;.
\\
\vspace{4pt}

\noindent
On the other hand, the BVP \eqref{BVP2} has a variational structure. Indeed, thanks to elliptic regularity theory, and the strong maximum principle a smooth solution of \eqref{BVP2} can be found by looking at critical points of the functional \;$J_q$\; defined by 
\begin{equation}\label{eq:dfjq}
J_q(u):=\frac{\left<u, u\right>_q}{\left(\oint_{M}u^{\frac{2}{1-\gamma}}\;dS_g\right)^{1-\gamma}}, \;\;\;\;u\in W^{1,2}_{y^{1-2\gamma},+}(X)=\{u\in W^{1,2}_{y^{1-2\gamma}}(X),\;u\ge 0,\;u\neq0\},
\end{equation}
and  
\begin{equation}\label{scalq}
\left<u, u\right>_q=d^{*}_{\gamma}\int_{X}y^{1-2\gamma}\left(|\nabla_g u|^{2}+qu^{2}\right)\;dV_g.
\end{equation}
In the definition \eqref{eq:dfjq}, 
\begin{equation}
  W^{1,2}_{y^{1-2\gamma}}(X)=\overline{C^{\infty}(X)}^{\|\cdot\|_{W^{1,2}_{y^{1-2\gamma}}(X)}},
\end{equation}
where
\begin{equation}
   \|u\|^{2}_{W^{1,2}_{y^{1-2\gamma}}(X)}=\int_{X}y^{1-2\gamma}\left(|\nabla_g u|^{2}+u^{2}\right)\;dV_g.
\end{equation}
\vspace{4pt}

\noindent
Moreover, \;$dV_{g}$\; denotes the volume form with respect to \;$g,$\; and \;$dS_{g}$\; denotes the volume form with respect to the Riemannian metric induced by \;$g$\; on \;$M.$\;\\
\vspace{4pt}

\noindent
Denoting by \;$H_{g}$\; the mean curvature of \;$M$\; in \;$(\overline{X},g)$.\;We prove the following theorem.
\begin{thm}\label{thm1}
Let \;$\left(X,g^{+}\right)$\; be an asymptotically hyperbolic manifold with conformal infinity \;$(M,[h])$\; of dimensional \;$n=2$,\; interior \;$X$,\; \;$H_{g}=0$,\; and \;$y$\; be the unique  geodesic defining function associated to \;$h$\; given by \eqref{metric.g+}. Let also \;$\gamma \in (0,1)$,\; and \;$q:X \rightarrow \R$\; be a smooth function verifying \;$q=O\left(y^{1-2\gamma}\right)$\; in a tubular neighborhood of \;$M$\;. Assuming that \;$\lambda_{1}(P^{\gamma}_q)>0$\;, and the Green's function \;$G$\; defined by \eqref{eqgreen} is strictly positive, then the BVP \eqref{BVP2} has a least one solution.
\end{thm}
\vspace{4pt}

\noindent
As in \cite{aldawood} and \cite{aldawood1}, to prove Theorem \ref{thm1}, we will use the Barycenter Technique of Bahri-Coron\cite{bc} which is possible since \;$dim\;M=2$\; imply the problem under study is a Global one (for the definition of "Gobal" for Yamabe and Cherrier-Escobar type problems, see \cite{nss}). Indeed, as in \cite{aldawood}, \cite{aldawood1}, and \cite{nss}, we will follow the scheme of the Algebraic topological argument of Bahri-Coron\cite{bc} as performed in the work \cite{mmc}. As in \cite{aldawood} and \cite{aldawood1}, one of the main difficulty with respect to the works  \cite{mmc} and \cite{nss} is the presence of the linear term \;$"qu"$\; and the lack of conformal invariance. As in \cite{aldawood} and \cite{aldawood1}, to deal with such a difficulty, we use the fact that \;$dim\;M=2$\; and the Brendle\cite{bre1}-Schoen\cite{sc}' s bubble construction to run the scheme of the Algebraic topological argument of Bahri-Coron\cite{bc} for  the existence. However, to deal with the fractional aspect of the problem, we use the assumption \;$q=O\left(y^{1-2\gamma}\right)$\; different from the case \cite{nss} where the geometry is used.  
\vspace{4pt}

\noindent
The paper is divided into the following sections. In Section \ref{NP}, we discuss some preliminaries and fix some notations. In Section \ref{PSS}, we recall the profile decomposition of Palais-Smale (PS)-sequences for \;$J_q.$\; We also introduce the neighborhoods of potential critical points at infinity of \;$J_q$\; and  their associated selection maps. In addition, we state a Deformation Lemma for \;$J_q$\; taking into account the possible bubbling phenomena involved in the study of \eqref{BVP2}. In Section \ref{SAE}, we derive some sharp self-action estimates needed for the application of the Barycenter Technique of Bahri-Coron\cite{bc} for existence. In Section \ref{IE}, we derive some sharp inter-action estimates needed for the application of the Barycenter Technique of Bahri-Coron\cite{bc} for existence. In Section \ref{ATA},  we present the algebraic topological argument for existence. In Section \ref{APP}, we collect some technical estimates.
\section*{Acknowledgement}
We would like to thank Dr. Cheikh Ndiaye for proposing the problem and for his daily support.
%
%
%
%
%
%
\section{Notations and preliminaries}\label{NP}
In this Section, we introduce some preliminary notions and fix some notations.
\vspace{4pt}

\noindent
For \;$\zeta=\left(\zeta_1,\zeta_2,\zeta_3\right)\in \bar{\R}^3_{+}=\R^2\times\bar{\R}_{+},$\; with \;$\bar{\R}_{+}=[0,+\infty),$\; we set \;$\bar{\zeta}=(\zeta_1,\zeta_2)\in \R^{2},$\; so that \;$\zeta=\left(\bar{\zeta},\zeta_3\right),$\; and we set also \;$r=|\zeta|=\sqrt{\zeta^{2}_1+\zeta^{2}_2+\zeta^{2}_3}.$\; We define \;$\R^3_{+}=\R^2\times\R_{+}$\; with \;$\R_{+}=(0,\infty).$\; We identify the boundary of \;$\R^3_{+}$\; denoted \;$\partial \R^3_{+}$\; with \;$\R^{2}.$\;
\vspace{6pt}

\noindent
For \;$r>0,$\; and \;$a\in M$\; and for a Riemannian metric \;$\bar{h}$\; defined on \;$M,$\; we denote by   \;$B^{\bar{h}}_{r}(a)$\; the geodesic ball with respect to \;$\bar{h}$\; of radius \;$r$\; centered at \;$a$.\; We also denote by \;$d_{\bar{h}}(x,y)$\; the geodesic distance with respect to \;$\bar{h}$\; between two points \;$x$\; and \;$y$\; of \;$M$,\; and \;$inj_{\bar{h}}(M)$\; stands for the
injectivity radius of \;$(M,\bar{h})$.\; Moreover, \;$dV_{\bar{h}}$\; denotes the Riemannian measure associated to the metric \;$\bar{h}$\; on \;$M$,\; and for \;$a\in M$\; we use the notation \;$exp^{a}_{\bar{h}}$\; to denote the exponential map with respect to \;$\bar{h}$\; on \;$M$\; centered at \;$a$.\;
\vspace{6pt}

\noindent
Similarly, for \;$r>0,$\; and \;$a\in M$\; and for a Riemannian metric \;$\bar{g}$\; defined \;$\overline{X},$\; we denote by \;$B^{\bar{g},+}_{a}(r)$\; the geodesic half ball in \;$\overline{X}$\; with respect to \;$\bar{g}$\; of radius \;$r$\; centered at \;$a$.\; We also denote  by \;$d_{\bar{g}}(x,y)$\; the geodesic distance with respect to \;$\bar{g}$\; between two points \;$x \in M$\; and \;$y \in \overline{X}$,\; and \;$inj_{\bar{g}}(\overline{X})$\; stands for the
injectivity radius of \;$(\overline{X},\bar{g})$.\; Moreover, \;$dV_{\bar{g}}$\; denotes the Riemannian measure associated to the metric \;$\bar{g}$\; on \;$\overline{X}$,\; and for \;$a\in M$\; we use the notation \;$exp^{\bar{g,+}}_{a}$\; to denote the exponential map with respect to \;$\bar{g}$\; on \;$\overline{X}$\; centered at \;$a$.\;
\vspace{6pt}

\noindent
$\mathbb{N}$\; denotes the set of nonnegative integers, \;$\mathbb{N}^{*}$\; the set of positive integers and for \;$k\in \mathbb{N}^{*},\;\;\R^{k}$\; stands for
the standard $k$-dimensional Euclidean space, \;$\R^{k}_{+}$\; the open positive half-space of \;$\R^{k}$,\; and \;$\bar{\R}^{k}_{+}$\; its closure in \;$\R^{k}$.\; For simplicity we use the notation \;$\R_{+}:= \R^{1}_{+}$,\; and \;$\bar{\R}_{+}:=\bar{\R}^{1}_{+}$.\; For \;$r>0$\; and \;$k=3$\; we denote respectively
$$B^{\R^{3}}_{r}(0)\;\;\text{and}\;\;B^{\R^{3}_{+}}_{r}(0)=B^{\R^{3}}_{r}(0)\cap \R^{3}_{+}\simeq]0,r[\times B^{\R^{2}}_{r}(0)$$
the open and open upper half ball of \;$\R^{3}$\; of center \;$0$\; and radius \;$r$,\; and set \;$B_{r}=B^{\R^{2}}_{r}(0)$,\; \;$B^{+}_{r}=B^{\R^{3}_{+}}_{r}(0)$,\; and \;$\partial B^{+}_{r}=B^{\R^{3}_{+}}_{r}(0) \cap \R^{2}$.\;
\vspace{6pt}

\noindent
The letter \;$C$\; is typically used to represent large positive constants, and the value of \;$C$\; can verify from formula to formula as well as within the same line.\\
We set 
\begin{equation}\label{uqsqrt}
    \|u\|_{q}=\sqrt{\left<u, u\right>_q}\;,\;\;u\in W^{1,2}_{y^{1-2\gamma}}(X),
\end{equation}
and \;$\left<u, u\right>_q$\; is as in \eqref{scalq}.\\ 
\vspace{4pt}

\noindent
For \;$a\in\R^{2}$\; and \;$\lambda>0,$\; we define the standard bubbles on \;$\R^{2}$\; as  
\begin{equation}\label{bubble-M3}
     \delta_{a,\lambda}(\bar{\zeta})=c_0\left[\frac{\lambda}{1+\lambda^2\left|\bar{\zeta}-a\right|^2}\right]^{1-\gamma},\;\;\;\;\bar{\zeta}\in \R^{2}
\end{equation}
with \;$c_0>0$\; is such that \;$\delta_{a,\lambda}$\; satisfies 
\begin{equation}\label{delta1}
(-\Delta_{\R^{2}})^{\gamma} \delta_{a,\lambda}=\delta^{\frac{1+\gamma}{1-\gamma}}_{a,\lambda}\qquad\qquad\;\;\;\;\text{on}\;\;\;\;\R^{2}.   
\end{equation}
We define \;$\bar{\delta}_{a,\lambda}(\bar{\zeta},y)$\; with \;$(\bar{\zeta},y)\in \bar{\R}^{3}_{+}$\; as the unique solution to
\begin{equation*}
    \left\{
\begin{split}
 -div(y^{1-2\gamma}\nabla\bar{\delta}_{a,\lambda})&=0&\;\;\;\;\text{in}\;\;\;\;&\R^3_{+},\\
\bar{\delta}_{a,\lambda}&=\delta_{a,\lambda}\;\;\;\;&\text{on}\;\;\;\;&\R^2.
\end{split}
\right.
\end{equation*}
\vspace{2pt}

\noindent
We have \;$\bar{\delta}_{a,\lambda}$\; verifies 
\begin{equation}\label{deltabar}
\left\{
\begin{split}
 -div(y^{1-2\gamma}\nabla\bar{\delta}_{a,\lambda})&=0&\;\;\;\;\text{in}\;\;\;\;&\R^3_{+},\\
-d^{*}_{\gamma}\lim_{y\to 0}y^{1-2\gamma}\frac{\partial \bar{\delta}_{a,\lambda}}{\partial y}&=\bar{\delta}^{\frac{1+\gamma}{1-\gamma}}_{a,\lambda}\;\;\;\;&\text{on}\;\;\;\;&\R^2,
\end{split}
\right.
\end{equation}
where
\begin{equation}\label{d*}
  d^{*}_{\gamma}=-\frac{d_{\gamma}}{2\gamma},  
\end{equation}
\;$d_{\gamma}=2^{2\gamma}\frac{\Gamma(\gamma)}{\Gamma(-\gamma)}$,\; and \;$\Gamma$\; is the standard  Gamma function. Also,
\begin{equation*}
    div(X)=\sum_{i=1}^3 \frac{\partial X^{i}}{\partial x_{i}},\;\;\; X=(x_{1},x_{2},x_{3})
\end{equation*}
is the Divergent on \;$\R^{3}.$\;\\
\vspace{4pt}

\noindent
We define 
\begin{equation}\label{c3}
c_1=\int_{\R^2}\left(\frac{1}{1+|y|^2}\right)^{1+\gamma}dy.
\end{equation}
\vspace{4pt}

\noindent
Furthermore, it is a well-know fact that  
\begin{equation}\label{S3}
d^{*}_{\gamma}\int_{\R^3_{+}}y^{1-2\gamma}|\nabla \bar{\delta}_{0,\lambda}|^{2}=d^{*}_{\gamma}\int_{\R^3_{+}}y^{1-2\gamma}|\nabla \bar{\delta}_{0,1}|^{2}=\int_{\R^2}\delta_{0,\lambda}^{\frac{2}{1-\gamma}}=\int_{\R^2}\delta_{0,1}^{\frac{2}{1-\gamma}}.
\end{equation}
\vspace{4pt}

\noindent
We set
\begin{equation}\label{S}
\mathcal{S}=\frac{d^{*}_{\gamma}\int_{\R^3_{+}}y^{1-2\gamma}|\nabla \bar{\delta}_{0,\lambda}|^{2}}{\left(\int_{\R^2}\delta_{0,\lambda}^{\frac{2}{1-\gamma}}\right)^{1-\gamma}}. 
\end{equation}

\vspace{4pt}

\noindent
For \;$a\in M$\; and \;$x \in \overline{X}$,\; we let \;$G(a,x)$\; be the unique solution of
\begin{equation}\label{eqgreen}
\left\{
\begin{split}
-div_{g}(y^{1-2\gamma}\nabla_{g}G(a,x))+qG(a,x)&=0,&\;\;\;&x \in X,\\
-d^{*}_{\gamma}\lim_{y\to 0} y^{1-2\gamma}\;\frac{\partial G(a,(z,y))}{\partial y}&=\delta_a(z),&\;\;\;& z\in M,\;\; x=(z,y).
\end{split}
\right.
\end{equation}
\vspace{4pt}

\noindent
Since \;$q$\; is a bounded smooth function defined on \;$X$\; and \;$dim\;M=2,$\; then Green's function \;$G(a,x)$\; satisfies the following estimates (see \cite{mmn}) 
  \begin{equation}\label{estg}
 \left|G(a,x)-\frac{1}{d_g(a,x)^{2-2\gamma}}\right|\leq C ,\;\;\;\;\text{for}\;\;\;\;x\neq a\in\overline{X},
 \end{equation}
 \begin{equation}\label{estgg}
 \left|\nabla \left(G(a,x)-\frac{1}{d_g(a,x)^{2-2\gamma}}\right)\right|\le \frac{C}{d_g(a,x)^{2-2\gamma}},\;\;\;\; \text{for}\;\;\;\;x\neq a\in \overline{X},
 \end{equation}
 and 
 \begin{equation}\label{estggg}
 \left|\nabla^{2} \left(G(a,x)-\frac{1}{d_g(a,x)^{2-2\gamma}}\right)\right|\le \frac{C}{d_g(a,x)^{3-2\gamma}},\;\;\;\; \text{for}\;\;\;\;x\neq a\in \overline{X},
 \end{equation}
 with \;$C$\; a positive constant.\\
\vspace{2pt}

\noindent
Moreover, we have \;$\overline{X}$\; compact implies that there exists 
 \begin{equation}\label{delta0}
     \delta_0>0
 \end{equation}
  such that \;$\forall$\;\;$a\in \partial M$\; and \;$\forall$\;\;$0<2\delta\le \delta_0,$\; the Fermi coordinates centered at \;$a$\; defines a smooth map
 \begin{equation}\label{psi}
     \Psi_a:B^{+}_{2\delta}\to \overline{X},
 \end{equation}
 identifying a neighborhood \;$O(a)$\; of \;$a$\; in \;$\overline{X}.$\; We will identify a point \;$x=\Psi_a(\zeta)\in O(a)$\; with \;$\zeta\in B^{+}_{2\delta}.$\; With this agreement and recalling that \;$H_g=0$\;, we have that an expansion of the Riemannian metric\;$g$\; and \;$\sqrt{|g|}$\; (where \;$|g|$\; denotes the modulus of the determinant of \;$g$\;) on \;$B^{g}_{2\delta}$\; is given by the following formulas 
 \begin{equation}\label{E_Metric}
 \begin{split}
     g^{ij}(\zeta)&=\delta_{ij}+2(L_g)_{ij}\zeta_n(0)+\frac{1}{3}R_{ikjl}[\hat{g}]\zeta_k\zeta_l(0)+g^{ij}_{nk}\zeta_n\zeta_k(0)+\left\{3(L_g)_{ij}(L_g)_{kj}+R_{injn}[g]\right\}\zeta^{2}_n(0)\\
     &+o(|\zeta|^3),\;\;\;\;\;\zeta\in B^{+}_{2\delta}\\
     \sqrt{|g(\zeta)|}&=1-\frac{1}{6}R_{ic}[\hat{g}]_{ij}\zeta_i \zeta_j(0)-\left[\frac{1}{2}\Arrowvert L_g\Arrowvert^2+R_{ij}[g]_{nn}\right]\zeta^2_{n}(0)+o(|\zeta|^3),\;\;\;\;\;\zeta\in B^{+}_{2\delta}.
     \end{split}
 \end{equation}
 In the formulas in \eqref{E_Metric}, \;$n=2,$\;\;$(L_g)_{ij}$\; denotes the component of the second fundamental form of \;$M$\; with respect to \;$g,$\; \;$(R_{ic}[g])_{ab},$ $a,b=1,..,n$\; denotes the component of the Ricci tensor of \;$\overline{X}$\; with respect to \;$g,$\;\;$\hat{g}:=g|_{M}$\; is the Riemannian metric induced by \;$g$\; on \;$M,$\;\;$(R_{ic}[\hat{g}])_{ij},$\;\;$i,j=1,..,n-1$\; are the components of the Ricci tensor of \;$M$\; with respect to \;$\hat{g},$\;\;$R_{abcd}[g],$\;\;$a,b,c,d=1,..,n$\; denotes the components of Riemann tensor of \;$\overline{X}$\; with respect to \;$g,$\; and \;$R_{ijkl}[\hat{g}],$\;\;$i,k=1,..,n-1$\; denotes the components of Riemann curvature tensor of \;$M$\; with respect to \;$\hat{g}.$\; All the tensors in the right of \eqref{E_Metric} are evaluated at \;$0,$\;and we also use Einstein summation convention for repeated indexes.
 \vspace{4pt}
 
 \noindent
 Let \(\ \chi:\mathbb{R} \to\R\)\; be a smooth cut-off function satisfying
\begin{equation}\label{chi}
\chi(t)= \left\{ \begin{array}{ll}
         1,&\;\;\mbox{if\;\;\;$t\leq  1$},\\
        0,&\;\;\mbox{if\;\;\;$t\geq 2$}.\end{array} \right.
        \end{equation}
\vspace{6pt}
 
 \noindent
For \;$0<2\delta<\delta_0,$\; we define
\begin{equation}\label{chid}  
        \chi_{\delta}(\zeta)=\chi\left(\frac{|\zeta|}{\delta}\right),\;\;\;\;\zeta\in \bar{\R}^3_{+}.
        \end{equation}
\vspace{6pt}
 
 \noindent
For \;$0<2\delta<\delta_0,$\; \;$a\in M,$\; and \;$\lambda>0,$\; we define the Brendle\cite{bre1}-Schoen\cite{sc}'s bubble 
\begin{equation}\label{uald}
       u_{a, \lambda}(x)=u_{a,\lambda,\delta}(x):=\chi^{a}_\delta(x)\hat{\delta}_{a,\lambda}(x) +(1-\chi^{a}_\delta(x)) \frac{c_0}{\lambda^{1-\gamma}}  G(a,x),\;\;\;\;\text{for}\;\;\;\;x\in \overline{X},
 \end{equation}
 \noindent
 where 
 \begin{equation}\label{deltahat}
  \hat{\delta}_{a,\lambda}(x)=\bar{\delta}_{0,\lambda}\left(\Psi^{-1}_{a}(x)\right),
 \end{equation}
 and
 \begin{equation}\label{chi1}
  \chi^{a}_\delta=\chi_\delta(\Psi^{-1}_a(x)),
 \end{equation}
 with
 \begin{equation}\label{psi0}
  \Psi_a:B^{+}_{\delta_0}\to \overline{X}.
 \end{equation}
 \vspace{4pt}
 
 \noindent
 Thus, recalling that we are under the assumption \;$G>0$\; (for the definition of \;$G$\; see \eqref{eqgreen}), then \;$\forall$\;\;$a\in M$\; and \;$\forall$\;\;$0<2\delta<\delta_0,$\; we have 
\begin{equation}\label{ual}
u_{a,\lambda}\in W^{1,2}_{y^{1-2\gamma}}(X), \;\;\;\;\text{and}\;\;\;\;u_{a, \l}>0\;\;\;\;\text{in}\;\;\;\;\overline{X}.
\end{equation}
For \;$a_i, a_j\in M,$\; and \;$\l_i, \l_j>0,$ we define
\begin{equation}\label{varepij}
\varepsilon_{ij}=\left[\frac{1}{\frac{\lambda_i}{\lambda_j}+\frac{\lambda_j}{\lambda_i}+\lambda_i\lambda_jG^{\frac{-1}{1-\gamma}}(a_i,a_j)}\right]^{1-\gamma}.
\end{equation}
\vspace{4pt}
 
 \noindent
Moreover,\;\;for \;$0<2\delta<\delta_0,$\;\;$a_i, a_j\in M,$\; and \;$\l_i, \l_j>0,$\; we define
\begin{equation}\label{epij}
\epsilon_{ij}=\oint_{M} u^{\frac{1+\gamma}{1-\gamma}}_{a_i,\lambda_i}u_{a_j,\lambda_j}\;dS_g
\end{equation}
and
\begin{equation}\label{eij}
e_{ij}=d^{*}_{\gamma}\left[\int_{X}y^{1-2\gamma}\nabla_g u_{a_i,\lambda_i} \nabla_gu_{a_j,\lambda_j}\;dV_g+\int_{X}qu_{a_i,\lambda_i} u_{a_j,\lambda_j}\;dV_g\right].
\end{equation}
\vspace{4pt}
 
 \noindent
 To conclude this section, we derive the following $C^{0}$-estimate needed for the energy and the inter-action estimates required for the application of the Barycenter technique of Bahri-Coron\cite{bc} for existence.
 %
 %
 %
 %
 %
 %
 \begin{lem}\label{c0estimate}
 Assuming that \;$\theta>0$\; is small, then there exists \;$C>0$\; such that \;$\forall$\;\;$a \in M,$\;\;$\forall$\;\;$0<2\delta<\delta_0$\; and \;$\forall$\;\;$0<\frac{1}{\l}\le\theta\delta,$\; we have
\begin{equation}\label{C_01}
\begin{split}
    \left|\left(-div_{g}(y^{1-2\gamma}\nabla_{g})+q\right)u_{a,\lambda}(x)\right|\le & C\;y^{1-2\gamma}\left[\left(\frac{1}{\delta^{2}\lambda^{1-\gamma}}+\frac{1}{\delta^{3-2\gamma}\lambda^{1-\gamma}}+\frac{1}{\delta^{6-2\gamma}\lambda^{3-\gamma}}\right) \large{1}_{\{z\in X:\;\delta\le d_g(a,z)\le 2\delta\}}(x)\right.\\
  +&\hat{\delta}_{a,\lambda}(x) \large{1}_{\{z\in X:\; d_g(a,z)\le 2\delta\}}(x)\\ 
  +&\left.\left(\frac{\lambda}{1+\lambda^{2}d_g(a,x)^{2}}\right)^{1-\gamma} \large{1}_{\{z\in X:\;d_g(a,z)\le 2\delta\}}(x)\right],\;\;\;\forall \;\;x\in X,
  \end{split}
\end{equation}
and
\begin{equation}\label{C_02}
\begin{split}
    \left|-d^{*}_{\gamma}\lim_{y\to 0} y^{1-2\gamma}\;\frac{\partial u_{a,\lambda}(z,y) }{\partial y}-u^{\frac{1+\gamma}{1-\gamma}}_{a,\lambda}(z)\right|\le  C&\left[ \left(\frac{\lambda}{1+\lambda^{2}d_g(a,z)^{2}}\right)^{1+\gamma}\large{1}_{\{\rho\in M:\;d_g(a,\rho)\geq \delta\}}(z)\right],\;\;\;\forall \;\;z\in M.
    \end{split}
\end{equation}
 \end{lem}
 \begin{pf}
 To simplifies notation, let us set \;  $G_a(\cdot):=G(a,\cdot)$\; and \;$\bar G_a(\cdot)=c_0G_a(\cdot)$. Then, we have 
\begin{equation}\label{CC_0}
u_{a,\lambda}(x)=\chi^{a}_\delta  (x)\hat{\delta}_{a,\lambda}(x)+(1-\chi^{a}_\delta(x))\frac{\bar{G}_a(x)}{\lambda^{1-\gamma}},\;\;\;\;x\in X.
\end{equation}
In order to deal with \eqref{C_01}, first we need to make the following observation
\begin{equation*}
    \begin{split}
     \left(D_{g}+q\right) u_{a,\lambda}(x)&=\left(D_{g}+q\right)\left[\chi^{a}_\delta(x)\hat{\delta}_{a,\lambda}(x)+\left(1-\chi^{a}_\delta(x)\right)\frac{\bar{G}_a(x)}{\lambda^{1-\gamma}}\right]\\
     &=\left(D_{g}+q\right)\left[\chi^{a}_\delta(x)\left(\hat{\delta}_{a,\lambda}(x)-\frac{\bar{G}_a(x)}{\lambda^{1-\gamma}}\right)\right]\\
     &+\left(D_{g}+q\right)\left[\frac{\bar{G}_a(x)}{\lambda^{1-\gamma}}\right],\;\;\;\;x\in X,  
    \end{split}
\end{equation*}
where \;$D_{g}=-div_{g}(y^{1-2\gamma}\nabla_{g}(\cdot))$.\\
\noindent
Now, since \;$x\in X$\; and \;$a\in M,$\; then \;$\left(-D_{g}+q\right)\left[\frac{\bar{G}_a(x)}{\lambda^{1-\gamma}}\right]=0.$\;
\vspace{4pt}
 
 \noindent
This implies
\begin{equation*}
    \begin{split}
     \left(D_{g}+q\right) u_{a,\lambda}(x)&=\left(D_{g}-D\right)\left[\chi^{a}_\delta(x)\left(\hat{\delta}_{a,\lambda}(x)-\frac{\bar{G}_a(x)}{\lambda^{1-\gamma}}\right)\right]+D\left[\chi^{a}_\delta(x)\left(\hat{\delta}_{a,\lambda}(x)-\frac{\bar{G}_a(x)}{\lambda^{1-\gamma}}\right)\right]\\
     &+q\left[\chi^{a}_\delta(x)\left(\hat{\delta}_{a,\lambda}(x)-\frac{\bar{G}_a(x)}{\lambda^{1-\gamma}}\right)\right]\;\;\;\;x\in X,
     \end{split}
     \end{equation*}
then, we get      
     \begin{equation*}
         \begin{split}
   \left(D_{g}+q\right) u_{a,\lambda}(x)&=\left(D_{g}-D\right)\left[\chi^{a}_\delta(x)\left(\hat{\delta}_{a,\lambda}(x)-\frac{c_0}{\lambda^{1-\gamma}d_g(a,x)^{2-2\gamma}}\right)\right]\\
   &+\left(D_{g}-D\right)\left[\chi^{a}_\delta(x)\left(\frac{c_0}{\lambda^{1-\gamma}d_g(a,x)^{2-2\gamma}}-\frac{\bar{G}_a(x)}{\lambda^{1-\gamma}}\right)\right]\\
   &+D\left[\chi^{a}_\delta(x)\left(\hat{\delta}_{a,\lambda}(x)-\frac{c_0}{\lambda^{1-\gamma}d_g(a,x)^{2-2\gamma}}\right)\right]\\
   &+D\left[\chi^{a}_\delta(x)\left(\frac{c_0}{\lambda^{1-\gamma}d_g(a,x)^{2-2\gamma}}-\frac{\bar{G}_a(x)}{\lambda^{1-\gamma}}\right)\right]\\
   &+q\left[\chi^{a}_\delta(x)\left(\hat{\delta}_{a,\lambda}(x)-\frac{\bar{G}_a(x)}{\lambda^{1-\gamma}}\right)\right],\;\;\;\;x\in X.  
    \end{split}
\end{equation*}
Thus,
$$
\left(D_{g}+q\right) u_{a,\lambda}(x)=\sum_{m=1}^4 I_m
$$
with 
\begin{equation*}
\begin{split}
&I_1=\left(D_{g}-D\right)\left[\chi^{a}_\delta(x)\left(\hat{\delta}_{a,\lambda}(x)-\frac{c_0}{\lambda^{1-\gamma}d_g(a,x)^{2-2\gamma}}\right)\right],\;\;\;\;x\in X,\\
&I_2=D_{g}\left[\chi^{a}_\delta(x)\left(\frac{c_0}{\lambda^{1-\gamma}d_g(a,x)^{2-2\gamma}}-\frac{\bar{G}_a(x)}{\lambda^{1-\gamma}}\right)\right],\;\;\;\;x\in X,\\
&I_3=D\left[\chi^{a}_\delta(x)\left(\hat{\delta}_{a,\lambda}(x)-\frac{c_0}{\lambda^{1-\gamma}d_g(a,x)^{2-2\gamma}}\right)\right],\;\;\;\;x\in X,\\
&I_4=q\left[\chi^{a}_\delta(x)\left(\hat{\delta}_{a,\lambda}(x)-\frac{\bar{G}_a(x)}{\lambda^{1-\gamma}}\right)\right],\;\;\;\;x\in X.
\end{split}
\end{equation*}

\vspace{4pt}

\noindent
We will estimate each of the \;$I_m$'s one at a time. For \;$I_2$,\; We write \;$I_2$\;in the following form 
\begin{equation}\label{I_2}
\begin{split}
I_2&=-div_{g}\left(y^{1-2\gamma}\nabla_{g}\left(\chi^{a}_\delta(x)\left(\frac{c_0}{\lambda^{1-\gamma}d_g(a,x)^{2-2\gamma}}-\frac{\bar{G}_a(x)}{\lambda^{1-\gamma}}\right)\right)\right)\\
&=-y^{1-2\gamma}\Delta_g\chi^{a}_\delta(x)\left[\frac{c_0}{\lambda^{1-\gamma}d_g(a,x)^{2-2\gamma}}-\frac{\bar{G}_a(x)}{\lambda^{1-\gamma}}\right]-2y^{1-2\gamma}\nabla_g\chi^{a}_\delta(x) \nabla_g\left[\frac{c_0}{\lambda^{1-\gamma}d_g(a,x)^{2-2\gamma}}-\frac{\bar{G}_a(x)}{\lambda^{1-\gamma}}\right]\\
&-y^{1-2\gamma}\chi^{a}_{\delta}(x)\Delta_{g}\left[\frac{c_0}{\lambda^{1-\gamma}d_g(a,x)^{2-2\gamma}}-\frac{\bar{G}_a(x)}{\lambda^{1-\gamma}}\right]\\
&-\partial_{y}y^{1-2\gamma}  \partial_{y}\chi^{a}_{\delta}(x) \left[\frac{c_0}{\lambda^{1-\gamma}d_g(a,x)^{2-2\gamma}}-\frac{\bar{G}_a(x)}{\lambda^{1-\gamma}}\right]=J^{2}_{1}+J^{2}_{2}+J^{2}_{3}+J^{2}_{4},\;\;\;\;x\in X.
\end{split}
\end{equation} 
To begin with \;$J^{2}_1,$\; we have 
\begin{equation}\label{J2_1}
J^{2}_1=-y^{1-2\gamma}\Delta_g\chi^{a}_\delta(x)\left[\frac{c_0}{\lambda^{1-\gamma}d_g(a,x)^{2-2\gamma}}-\frac{\bar G_a(x)}{\lambda^{1-\gamma}}\right],\;\;\;\;x\in X.
\end{equation} 
Applying \eqref{bubble-M3} and \eqref{estg}, we derive
\begin{equation}\label{diff2}
\left|\frac{c_0}{\lambda^{1-\gamma}d_g(a,x)^{2-2\gamma}}-\frac{\bar G_a(x)}{\lambda^{1-\gamma}}\right|\le\frac{C}{\lambda^{1-\gamma}},\;\;\;\;x\in X,\;\;\;\;\delta\le d_g(a,z)\le 2\delta.
\end{equation}
Now, we obtain the following by using \eqref{chi} and \eqref{chid}:
\begin{equation}\label{flchiJ1}
|\Delta_g\chi^{a}_\delta(x)|\le \frac{C}{\delta^{2}}\large{1}_{\{z\in X:\;\delta\le d_g(a,z)\le 2\delta\}}(x),\;\;\;\;x\in X.
\end{equation}
Hence, combining  \eqref{J2_1}-\eqref{flchiJ1},  we obtain
\begin{equation}\label{estJ21}
|J^{2}_1|\le C\; y^{1-2\gamma}\left(\frac{1}{\delta^{2}\lambda^{1-\gamma}}\right) \large{1}_{\{z\in X:\;\delta\le d_g(a,z)\le 2\delta\}}(x),\;\;\;\;x\in X.
\end{equation}
\vspace{6pt}
 
 \noindent
In order to estimate \;$J^{2}_2,$\; we first write
\begin{equation}\label{J2_2}
J^{2}_2=-2y^{1-2\gamma}\left<\nabla_g\chi^{a}_\delta(x), \;\nabla_g\left[\frac{c_0}{\lambda^{1-\gamma}d_g(a,x)^{2-2\gamma}}-\frac{\bar G_a(x)}{\lambda^{1-\gamma}}\right]\right>,\;\;\;\;x\in X.
\end{equation}
\vspace{4pt}
 
 \noindent
In the next step, using \eqref{bubble-M3} and \eqref{estgg}, we have the following:
\begin{equation}\label{diff4}
\left|\nabla_g\left[\frac{c_0}{\lambda^{1-\gamma}d_g(a,x)^{2-2\gamma}}-\frac{\bar G_a(x)}{\lambda^{1-\gamma}}\right]\right|\le \frac{C}{\delta^{2-2\gamma}\lambda^{1-\gamma}},\;\;\;\;x\in X,\;\;\;\;\delta\le d_g(a,z)\le 2\delta.
\end{equation}
\vspace{4pt}
 
 \noindent
On the other hand, using \eqref{chi}, we derive 
\begin{equation}\label{fgchiJ2}
\left|\nabla_g \chi^{a}_\delta(x)\right|\le \frac{C}{\delta} \large{1}_{\{z\in X:\;\delta\le d_g(a,z)\le 2\delta\}}(x),\;\;\;\;x\in X.
\end{equation}
Hence, combining \eqref{J2_2}-\eqref{fgchiJ2}, we get
\begin{equation}\label{estJ22}
|J^{2}_2|\le C\; y^{1-2\gamma}\left(\frac{1}{\delta^{3-2\gamma}\lambda^{1-\gamma}}\right) \large{1}_{\{z\in X:\;\delta\le d_g(a,z)\le 2\delta\}}(x),\;\;\;\;x\in X.
\end{equation}
\vspace{6pt}
 
 \noindent
Moreover, for \;$J^{2}_3,$\; we have
\begin{equation}\label{J2_3}
J^{2}_3=-y^{1-2\gamma}\chi^{a}_{\delta}(x)\Delta_{g}\left[\frac{c_0}{\lambda^{1-\gamma}d_g(a,x)^{2-2\gamma}}-\frac{\bar{G}_a(x)}{\lambda^{1-\gamma}}\right],\;\;\;\;x\in X.
\end{equation}
\vspace{4pt}
 
 \noindent
Now, using \eqref{bubble-M3} and \eqref{estggg}, we  obtain
\begin{equation}\label{fgchiJ3}
\left|\nabla^{2}_g\left[\frac{c_0}{\lambda^{1-\gamma}d_g(a,x)^{2-2\gamma}}-\frac{\bar G_a(x)}{\lambda^{1-\gamma}}\right]\right|\le \frac{C}{\delta^{3-2\gamma}\lambda^{1-\gamma}},\;\;\;\;x\in X,\;\;\;\;\delta\le d_g(a,z)\le 2\delta.
\end{equation}
Hence, combining \eqref{J2_3} and \eqref{fgchiJ3}, we have
\begin{equation}\label{estJ23}
|J^{2}_3|\le C\; y^{1-2\gamma}\left(\frac{1}{\delta^{3-2\gamma}\lambda^{1-\gamma}}\right) \large{1}_{\{z\in X:\;\delta\le d_g(a,z)\le 2\delta\}}(x),\;\;\;\;x\in X.
\end{equation}
 
 \vspace{6pt}
 
 \noindent
Now, for \;$J^{2}_4,$\; we have 
\begin{equation}\label{J2_4}
J^{2}_4=-(1-2\gamma)y^{-2\gamma}\partial_{y}\chi^{a}_{\delta}(x)\left[\frac{c_0}{\lambda^{1-\gamma}d_g(a,x)^{2-2\gamma}}-\frac{\bar G_a(x)}{\lambda^{1-\gamma}}\right],\;\;\;\;x\in X.
\end{equation} 
\noindent
 Using \eqref{bubble-M3} and
\eqref{estg}, we derive 
\begin{equation}\label{diff3}
\left|\frac{c_0}{\lambda^{1-\gamma}d_g(a,x)^{2-2\gamma}}-\frac{\bar G_a(x)}{\lambda^{1-\gamma}}\right|\le\frac{C}{\lambda^{1-\gamma}},\;\;\;\;x\in X,\;\;\;\;\delta\le d_g(a,z)\le 2\delta.
\end{equation}
\noindent
Next, using \eqref{chi} and \eqref{chid}, we have the following:
\begin{equation}\label{fgchiJ4}
|\partial_{y}\chi^{a}_\delta(x)|\le C \frac{y}{\delta^{2}}\large{1}_{\{z\in X:\;\delta\le d_g(a,z)\le 2\delta\}}(x),\;\;\;\;x\in X.
\end{equation}
Hence, combining  \eqref{J2_4}-\eqref{fgchiJ4},  we obtain
\begin{equation}\label{estJ24}
|J^{2}_4|\le C\; y^{1-2\gamma}\left(\frac{1}{\delta^{2}\lambda^{1-\gamma}}\right) \large{1}_{\{z\in X:\;\delta\le d_g(a,z)\le 2\delta\}}(x),\;\;\;\;x\in X.
\end{equation}
Thus, using \eqref{estJ21}, \eqref{estJ22}, \eqref{estJ23}, and \eqref{estJ24}, we derive the following for \;$I_2$\;
\begin{equation}\label{I_2}
    |I_2|\le C\; y^{1-2\gamma}\left(\frac{1}{\delta^{2}\lambda^{1-\gamma}}+\frac{1}{\delta^{3-2\gamma}\lambda^{1-\gamma}}\right) \large{1}_{\{z\in X:\;\delta\le d_g(a,z)\le 2\delta\}}(x),\;\;\;\;x\in X.
\end{equation}
\vspace{6pt}
 
 \noindent
In order to estimate \;$I_3$\;, we first write
\begin{equation}
\begin{split}
    I_3&=-div\left(y^{1-2\gamma}\nabla\left(\chi^{a}_\delta(x)\left(\hat{\delta}_{a,\lambda}(x)-\frac{c_0}{\lambda^{1-\gamma}d_g(a,x)^{2-2\gamma}}\right)\right)\right)\\
    &=-y^{1-2\gamma}\Delta\chi^{a}_\delta(x)\left[\hat{\delta}_{a,\lambda}(x)-\frac{c_0}{\lambda^{1-\gamma}d_g(a,x)^{2-2\gamma}}\right]-2y^{1-2\gamma}\nabla\chi^{a}_{\delta}(x)\left[\hat{\delta}_{a,\lambda}(x)-\frac{c_0}{\lambda^{1-\gamma}d_g(a,x)^{2-2\gamma}}\right]\\
    &-\partial_{y}y^{1-2\gamma}\partial_{y}\chi^{a}_{\delta}(x)\left[\hat{\delta}_{a,\lambda}(x)-\frac{c_0}{\lambda^{1-\gamma}d_g(a,x)^{2-2\gamma}}\right]=J^{3}_{1}+J^{3}_{2}+J^{3}_{3},\;\;\;\;x\in X.
\end{split}    
\end{equation}
To start with \;$J^{3}_1,$\; we have 
\begin{equation}\label{J3_1}
J^{3}_1=-y^{1-2\gamma}\Delta\chi^{a}_\delta(x)\left[\hat{\delta}_{a,\lambda}(x)-\frac{c_0}{\lambda^{1-\gamma}d_g(a,x)^{2-2\gamma}}\right],\;\;\;\;x\in X.
\end{equation} 
Now, applying \eqref{bubble-M3} and \eqref{estg}, we derive
\begin{equation}\label{diff2}
\left|\hat{\delta}_{a,\lambda}(x)-\frac{c_0}{\lambda^{1-\gamma}d_g(a,x)^{2-2\gamma}}\right|\le\frac{C}{\delta^{4-2\gamma}\lambda^{3-\gamma}},\;\;\;\;x\in X,\;\;\;\;\delta\le d_g(a,z)\le 2\delta.
\end{equation}
Using \eqref{chi} and \eqref{chid}, we get

\begin{equation}\label{flchi}
|\Delta\chi^{a}_\delta(x)|\le \frac{C}{\delta^{2}}\large{1}_{\{z\in X:\;\delta\le d_g(a,z)\le 2\delta\}}(x),\;\;\;\;x\in X.
\end{equation}
Hence, combining  \eqref{J3_1}-\eqref{flchi},  we obtain
\begin{equation}\label{estJ31}
|J^{3}_1|\le C\; y^{1-2\gamma}\left(\frac{1}{\delta^{6-2\gamma}\lambda^{3-\gamma}}\right) \large{1}_{\{z\in X:\;\delta\le d_g(a,z)\le 2\delta\}}(x),\;\;\;\;x\in X.
\end{equation}
\vspace{6pt}
 
 \noindent
For \;$J^{3}_2,$\; we first write
\begin{equation}\label{J3_2}
J^{3}_2=-2y^{1-2\gamma}\left<\nabla\chi^{a}_\delta(x), \;\nabla\left[\hat{\delta}_{a,\lambda}(x)-\frac{c_0}{\lambda^{1-\gamma}d_g(a,x)^{2-2\gamma}}\right]\right>,\;\;\;\;x\in X.
\end{equation}
\vspace{4pt}
 
 \noindent
In the next step, using \eqref{bubble-M3} and \eqref{estgg}, we have the following:
\begin{equation}\label{diff4}
\left|\nabla\left[\hat{\delta}_{a,\lambda}(x)-\frac{c_0}{\lambda^{1-\gamma}d_g(a,x)^{2-2\gamma}}\right]\right|\le \frac{C}{\delta^{5-2\gamma}\lambda^{3-\gamma}},\;\;\;\;x\in X,\;\;\;\;\delta\le d_g(a,z)\le 2\delta.
\end{equation}
\vspace{4pt}
 
 \noindent
On the other hand, using \eqref{chi}, we derive 
\begin{equation}\label{fgchi}
\left|\nabla \chi^{a}_\delta(x)\right|\le \frac{C}{\delta} \large{1}_{\{z\in X:\;\delta\le d_g(a,z)\le 2\delta\}}(x),\;\;\;\;x\in X.
\end{equation}
Hence, combining \eqref{J3_2}-\eqref{fgchi}, we get
\begin{equation}\label{estJ32}
|J^{3}_2|\le C\; y^{1-2\gamma}\left(\frac{1}{\delta^{6-2\gamma}\lambda^{3-\gamma}}\right) \large{1}_{\{z\in X:\;\delta\le d_g(a,z)\le 2\delta\}}(x),\;\;\;\;x\in X.
\end{equation}

\vspace{6pt}
 
 \noindent
 Now, for \;$J^{3}_3,$\; we have 
\begin{equation}\label{J3_3}
J^{3}_3=-(1-2\gamma)y^{-2\gamma}\partial_{y}\chi^{a}_{\delta}(x)\left[\hat{\delta}_{a,\lambda}(x)-\frac{c_0}{\lambda^{1-\gamma}d_g(a,x)^{2-2\gamma}}\right],\;\;\;\;x\in X.
\end{equation} 
\noindent
Next, using \eqref{bubble-M3} and
\eqref{estg}, we derive 
\begin{equation}\label{diff3}
\left|\hat{\delta}_{a,\lambda}(x)-\frac{c_0}{\lambda^{1-\gamma}d_g(a,x)^{2-2\gamma}}\right|\le\frac{C}{\delta^{4-2\gamma}\lambda^{1-\gamma}},\;\;\;\;x\in X,\;\;\;\;\delta\le d_g(a,z)\le 2\delta.
\end{equation}
\noindent
Using \eqref{chi} and \eqref{chid}, we get
\begin{equation}\label{flchi3}
|\partial_{y}\chi^{a}_\delta(x)|\le C \frac{y}{\delta^{2}}\large{1}_{\{z\in X:\;\delta\le d_g(a,z)\le 2\delta\}}(x),\;\;\;\;x\in X.
\end{equation}
Hence, combining  \eqref{J2_4}-\eqref{flchi3},  we have
\begin{equation}\label{estJ33}
|J^{3}_3|\le C\; y^{1-2\gamma}\left(\frac{1}{\delta^{6-2\gamma}\lambda^{1-\gamma}}\right) \large{1}_{\{z\in X:\;\delta\le d_g(a,z)\le 2\delta\}}(x),\;\;\;\;x\in X.
\end{equation}
Finally, using \eqref{estJ31}, \eqref{estJ32}, and \eqref{estJ33}, we obtain the following for \;$I_3$\;
\begin{equation}\label{I_3}
    |I_3|\le C\; y^{1-2\gamma}\left(\frac{1}{\delta^{6-2\gamma}\lambda^{3-\gamma}}\right) \large{1}_{\{z\in X:\;\delta\le d_g(a,z)\le 2\delta\}}(x),\;\;\;\;x\in X.
\end{equation}
\vspace{6pt}
 
 \noindent
To  estimate \;$I_4,$\; using \eqref{chi} and \eqref{chid}, we get
\begin{equation}\label{I_4}
    |I_4|\le C\;y^{1-2\gamma} \hat{\delta}_{a,\lambda}(x) \large{1}_{\{z\in X:\; d_g(a,z)\le 2\delta\}}(x),\;\;\;\;x\in X.
\end{equation}
\vspace{6pt}
 
 \noindent
To deal with \;$I_1,$\; we first estimate \; $\left(D_g-D\right)\left[\hat{\delta}_{a,\lambda}(x)-\frac{c_0}{\lambda^{1-\gamma}d_g(a,x)^{2-2\gamma}}\right]$\; for \;$\Psi_{a}(B^{+}_{2\delta})$\; (see \eqref{psi}). For this, identifying \;$x \in \Psi_{a}(\zeta)$\; with \;$\zeta\in B^{+}_{2\delta},$\; we have 
\begin{equation}
\begin{split}
D_g\left[\hat{\delta}_{a,\lambda}(x)-\frac{c_0}{\lambda^{1-\gamma}d_g(a,x)^{2-2\gamma}}\right]&=div_{g}\left(y^{1-2\gamma}\nabla_{g}\left(\hat{\delta}_{a,\lambda}(x)-\frac{c_0}{\lambda^{1-\gamma}d_g(a,x)^{2-2\gamma}}\right)\right)\\
&=\frac{1}{\sqrt{|g(x)|}}\partial_{i}\left[g^{i,j}(x)\sqrt{|g(x)|}y^{1-2\gamma}\partial_{j}\left(\bar{\delta}_{0,\lambda}(x)-\frac{c_0}{\lambda^{1-\gamma}d_g(a,x)^{2-2\gamma}}\right)\right]\\
&+\frac{1}{\sqrt{|g(x)|}}\partial_{n}\left[\sqrt{|g(x)|}y^{1-2\gamma}\partial_{n}\left(\bar{\delta}_{0,\lambda}(x)-\frac{c_0}{\lambda^{1-\gamma}d_g(a,x)^{2-2\gamma}}\right)\right]
\end{split}
\end{equation}
 \noindent
for\;$x\in \Psi_{a}(\zeta)$\;, then, we get
\begin{equation}\label{estI1}
    \begin{split}
    D_g\left[\hat{\delta}_{a,\lambda}(x)-\frac{c_0}{\lambda^{1-\gamma}d_g(a,x)^{2-2\gamma}}\right]&=\partial_{i}\log(\sqrt{|g(x)|})g^{i,j}(x)y^{1-2\gamma}\partial_{j}\left(\bar{\delta}_{0,\lambda}(x)-\frac{c_0}{\lambda^{1-\gamma}|x|^{2-2\gamma}}\right)\\
    &+y^{1-2\gamma}\partial_{i}\left[g^{i,j}(x)\partial_{j}\left(\bar{\delta}_{0,\lambda}(x)-\frac{c_0}{\lambda^{1-\gamma}|x|^{2-2\gamma}}\right)\right]\\
    &+\partial_{n}\log(\sqrt{|g(x)|})y^{1-2\gamma}\partial_{n}\left(\bar{\delta}_{0,\lambda}(x)-\frac{c_0}{\lambda^{1-\gamma}|x|^{2-2\gamma}}\right)\\
    &+y^{1-2\gamma}\partial^{2}_{n}\left(\bar{\delta}_{0,\lambda}(x)-\frac{c_0}{\lambda^{1-\gamma}|x|^{2-2\gamma}}\right)\\
    &+(1-2\gamma)y^{-2\gamma}\partial_{n}\left(\bar{\delta}_{0,\lambda}(x)-\frac{c_0}{\lambda^{1-\gamma}|x|^{2-2\gamma}}\right),\text{for}\;\;\;x\in \Psi_{a}(\zeta).
    \end{split}
\end{equation}
Now, we write \eqref{estI1} the following
\begin{equation}\label{estI111}
     \begin{split}
      D_g\left[\hat{\delta}_{a,\lambda}(x)-\frac{c_0}{\lambda^{1-\gamma}d_g(a,x)^{2-2\gamma}}\right]&=g^{i,j}\frac{x_{j}}{r}\partial_{i}\log(\sqrt{|g(x)|})\partial_{r}\left(\bar{\delta}_{0,\lambda}(x)-\frac{c_0}{\lambda^{1-\gamma}|x|^{2-2\gamma}}\right)y^{1-2\gamma}\\
      &+y^{1-2\gamma}\partial_{i}\left(g^{i,j}\frac{x_{j}}{r}\right)\partial_{r}\left(\bar{\delta}_{0,\lambda}(x)-\frac{c_0}{\lambda^{1-\gamma}|x|^{2-2\gamma}}\right)\\&+y^{1-2\gamma}g^{i,j}\frac{x_{i}x_{j}}{r^{2}}\partial^{2}_{r}\left(\bar{\delta}_{0,\lambda}(x)-\frac{c_0}{\lambda^{1-\gamma}|x|^{2-2\gamma}}\right)\\
      &+y^{1-2\gamma}\partial_{n}\log(\sqrt{|g(x)|})\partial_{n}\left(\bar{\delta}_{0,\lambda}(x)-\frac{c_0}{\lambda^{1-\gamma}|x|^{2-2\gamma}}\right)\\
      &+y^{1-2\gamma}\partial^{2}_{n}\left(\bar{\delta}_{0,\lambda}(x)-\frac{c_0}{\lambda^{1-\gamma}|x|^{2-2\gamma}}\right)\\
      &+(1-2\gamma)y^{-2\gamma}\partial_{n}\left(\bar{\delta}_{0,\lambda}(x)-\frac{c_0}{\lambda^{1-\gamma}|x|^{2-2\gamma}}\right),\text{for}\;\;\;x\in \Psi_{a}(\zeta).
     \end{split}
 \end{equation}
Using \eqref{estI111}, \eqref{divD}, and recalling that identifying \;$x\in \Psi_{a}(\zeta)$\; we have for \;$x \in B^{+}_{2\delta}$\; 
 \begin{equation}\label{estI11}
 \begin{split}
     &\left(D_g-D\right)\left[\hat{\delta}_{a,\lambda}(x)-\frac{c_0}{\lambda^{1-\gamma}d_g(a,x)^{2-2\gamma}}\right]\\
     &=y^{1-2\gamma}\left(g^{i,j}(x)\frac{x_i x_j}{r^2}-1\right)\partial^{2}_{r}\left(\bar{\delta}_{0,\lambda}(x)-\frac{c_0}{\lambda^{1-\gamma}|x|^{2-2\gamma}}\right)\\
     &+y^{1-2\gamma}\left(g^{i,j}(x)\frac{x_j}{r}\partial_{i}\log(\sqrt{|g(x)|})+\partial_{i}(g^{i,j}(x)\frac{x_j}{r})-(\frac{n-1}{r})\right)\partial_{r}\left(\bar{\delta}_{0,\lambda}(x)-\frac{c_0}{\lambda^{1-\gamma}|x|^{2-2\gamma}}\right)\\
     &+y^{1-2\gamma}\partial_{n}\log(\sqrt{|g(x)|})\partial_{n}\left(\bar{\delta}_{0,\lambda}(x)-\frac{c_0}{\lambda^{1-\gamma}|x|^{2-2\gamma}}\right),\;\;\;\;x \in B^{+}_{2\delta},
     \end{split}
\end{equation}
where
\begin{equation}\label{divD}
\begin{split}
  D\left[\bar{\delta}_{0,\lambda}(x)-\frac{c_0}{\lambda^{1-\gamma}|x|^{2-2\gamma}}\right]&=y^{1-2\gamma}\partial^{2}_{r}\left(\bar{\delta}_{0,\lambda}(x)-\frac{c_0}{\lambda^{1-\gamma}|x|^{2-2\gamma}}\right)\\
  &+y^{1-2}\left(\frac{n-1}{r}\right)\partial_{r}\left(\bar{\delta}_{0,\lambda}(x)-\frac{c_0}{\lambda^{1-\gamma}d_g(a,x)^{2-2\gamma}}\right)\\
  &+y^{1-2\gamma}\partial^{2}_{n}\left(\bar{\delta}_{0,\lambda}(x)-\frac{c_0}{\lambda^{1-\gamma}|x|^{2-2\gamma}}\right)\\
  &+(1-2\gamma)y^{-2\gamma}\partial_{n}\left(\bar{\delta}_{0,\lambda}(x)-\frac{c_0}{\lambda^{1-\gamma}|x|^{2-2\gamma}}\right).
  \end{split}
\end{equation}
 \noindent
In \eqref{estI11}, \;$n=2,$\; and $\;i=1,..,n,\;$
$$div=\sum_{i=1}^n \frac{\partial}{\partial x_{i}},\;\;\;\;\;\partial_{i}=\frac{\partial}{\partial x_{i}},\;\;\;\;\;\text{and}\;\;\;\;\partial_{r}=\frac{\partial}{\partial r}.$$\;\\
 \noindent
 Since 
 $$D\left[\bar{\delta}_{0,\lambda}(x)-\frac{c_0}{\lambda^{1-\gamma}|x|^{2-2\gamma}}\right]=0,\;\;\;\text{for}\;\;\;\;x \in B^{+}_{2\delta}\cap \R^{3}_{+},$$
 then, \eqref{estI11} implies
 \begin{equation}\label{estI12}
 \begin{split}
     &D_g\left[\hat{\delta}_{a,\lambda}(x)-\frac{c_0}{\lambda^{1-\gamma}d_g(a,x)^{2-2\gamma}}\right]\\
     &=\left(g^{i,j}(x)\frac{x_i x_j}{r^2}-1\right)y^{1-2\gamma}\partial^{2}_{r}\left(\bar{\delta}_{0,\lambda}(x)-\frac{c_0}{\lambda^{1-\gamma}|x|^{2-2\gamma}}\right)\\
     &+\left(g^{i,j}(x)\frac{x_j}{r}\partial_{i}\log(\sqrt{|g(x)|})+\partial_{i}(g^{i,j}(x)\frac{x_j}{r})-(\frac{n-1}{r})\right)y^{1-2\gamma}\partial_{r}\left(\bar{\delta}_{0,\lambda}(x)-\frac{c_0}{\lambda^{1-\gamma}|x|^{2-2\gamma}}\right)\\
     &+y^{1-2\gamma}\partial_{n}\log(\sqrt{|g(x)|})\partial_{i}\left(\bar{\delta}_{0,\lambda}(x)-\frac{c_0}{\lambda^{1-\gamma}|x|^{2-2\gamma}}\right),\;\;\;\;x\in B^{+}_{2\delta}.
     \end{split}
\end{equation}
\vspace{4pt}

 \noindent
Combining \eqref{E_Metric} and \eqref{estI12}, we get
 \begin{equation}\label{estI13}
 \begin{split}
     \left|D_g\left[\hat{\delta}_{a,\lambda}(x)-\frac{c_0}{\lambda^{1-\gamma}d_g(a,x)^{2-2\gamma}}\right]\right|&\le C |x|^{2}y^{1-2\gamma}\partial^{2}_{r}\left(\bar{\delta}_{0,\lambda}(x)-\frac{c_0}{\lambda^{1-\gamma}|x|^{2-2\gamma}}\right)\\
     &+C |x|y^{1-2\gamma}\left(\partial_{r}+\partial_{n}\right)\left(\bar{\delta}_{0,\lambda}(x)-\frac{c_0}{\lambda^{1-\gamma}|x|^{2-2\gamma}}\right),\;\;\;\;x\in B^{+}_{2\delta}.
    \end{split}
\end{equation}
Using Lemma \ref{Appendix} in \eqref{estI13}, we obtain the following estimate for \;$I_1$\; 
\begin{equation}\label{I_1}
    |I_1|\le C\; y^{1-2\gamma}\left(\frac{\lambda}{1+\lambda^{2}d_g(a,x)^{2}}\right)^{1-\gamma} \large{1}_{\{z\in X:\;d_g(a,z)\le 2\delta\}}(x),\;\;\;\;x\in X.
\end{equation}
Hence, the result for \eqref{C_01} follows from \eqref{I_2},\eqref{I_3}, \eqref{I_4},and \eqref{I_1}.\\
\vspace{6pt}
 
\noindent
To prove the formula \eqref{C_02}, we first take \;$z\in M$,\; and choose  \;$x=(z,y)$\; in a
tubular neighborhood of \;$M$\; in \;$\overline{X}$.\; Next, we apply  \eqref{CC_0} to the formula \eqref{C_02} to derive 
\begin{equation*}
\begin{split}
    -d^{*}_{\gamma}\lim_{y\to 0} y^{1-2\gamma}\;\frac{\partial u_{a,\lambda}(x) }{\partial y}-u^{\frac{1+\gamma}{1-\gamma}}_{a,\lambda}(z)&=\left(-d^{*}_{\gamma}\lim_{y\to 0} y^{1-2\gamma}\;\frac{\partial}{\partial y}\right)\left[\chi^{a}_\delta(x)\hat{\delta}_{a,\lambda}(x)+\left(1-\chi^{a}_\delta(x)\right)\frac{\bar{G}_a(x)}{\lambda^{1-\gamma}}\right]\\
    &-u^{\frac{1+\gamma}{1-\gamma}}_{a,\lambda}(z)\\
    &=-d^{*}_{\gamma}\lim_{y\to 0} y^{1-2\gamma}\;\frac{\partial}{\partial y}\left[\chi^{a}_\delta(x)\hat{\delta}_{a,\lambda}(x)\right]\\
    &-d^{*}_{\gamma}\lim_{y\to 0} y^{1-2\gamma}\;\frac{\partial}{\partial y}\left[\left(1-\chi^{a}_\delta(x)\right)\frac{\bar{G}_a(x)}{\lambda^{1-\gamma}}\right]-u^{\frac{1+\gamma}{1-\gamma}}_{a,\lambda}(z).
    \end{split}
\end{equation*}
To continue, we write 
\begin{equation*}
\begin{split}
 -d^{*}_{\gamma}\lim_{y\to 0} y^{1-2\gamma}\;\frac{\partial}{\partial y}\left[\left(1-\chi^{a}_\delta(x)\right)\frac{\bar{G}_a(x)}{\lambda^{1-\gamma}}\right]&=-d^{*}_{\gamma}\lim_{y\to 0} y^{1-2\gamma}\;\frac{\partial}{\partial y}\frac{\bar{G}_a(x)}{\lambda^{1-\gamma}}\left(1-\chi^{a}_\delta(x)\right)\\
 &+d^{*}_{\gamma}\lim_{y\to 0} y^{1-2\gamma}\;\frac{\partial}{\partial y}\chi^{a}_\delta(x)\frac{\bar{G}_a(x)}{\lambda^{1-\gamma}}.
 \end{split}
\end{equation*}
Next, using the definition\;$\chi^{a}_\delta$\; (see \eqref{chi1}), the symmetry of \;$\chi^{a}_\delta$\; in \;$\Psi_{a}(B^{+}_{2\delta})$\; after passing to Euclidean coordinates, and \;$\frac{\partial G(a,z)}{\partial  y}=0$\; for \;$z\in M,$\;\;$z\neq a,$\; we have
\begin{equation*}
   -d^{*}_{\gamma}\lim_{y\to 0} y^{1-2\gamma}\;\frac{\partial}{\partial y}\left[\left(1-\chi^{a}_\delta(x)\right)\frac{\bar{G}_a(z)}{\lambda^{1-\gamma}}\right]=0,\;\;\;\;\text{for}\;\;\;\;z\in M,\;\;\;\;\text{and}\;\;\;\;z\neq a.  
\end{equation*}
Hence for \;$z\in M$\; and \;$z\neq a,$\; we have 
\begin{equation*}
 -d^{*}_{\gamma}\lim_{y\to 0} y^{1-2\gamma}\;\frac{\partial u_{a,\lambda}(x) }{\partial y}-u^{\frac{1+\gamma}{1-\gamma}}_{a,\lambda}(z)=-d^{*}_{\gamma}\lim_{y\to 0} y^{1-2\gamma}\;\frac{\partial}{\partial y}\left[\chi^{a}_\delta(x)\hat{\delta}_{a,\lambda}(x)\right]-u^{\frac{1+\gamma}{1-\gamma}}_{a,\lambda}(z).  
\end{equation*}
Using again the definition of \;$\chi^{a}_\delta$\; and the symmetry of \;$\chi^{a}_\delta$\; in  \;$\Psi_{a}(B^{+}_{2\delta})$\; as before, we obtain 
\begin{equation}\label{C_022}
 -d^{*}_{\gamma}\lim_{y\to 0} y^{1-2\gamma}\;\frac{\partial u_{a,\lambda}(x) }{\partial y}-u^{\frac{1+\gamma}{1-\gamma}}_{a,\lambda}(z)=-d^{*}_{\gamma}\lim_{y\to 0} y^{1-2\gamma}\chi^{a}_\delta(x)\frac{\partial \hat{\delta}_{a,\lambda}(x)}{\partial y}-u^{\frac{1+\gamma}{1-\gamma}}_{a,\lambda}(z),\;\;\;\;z\in M,\;\;\;z\neq a.
\end{equation}
Clearly, \eqref{C_022} is true for \;$z=a.$\; On other hand, since identifying \;$x=(z,y)$\; with \;$\Psi^{-1}_{a}(x),$\; in \;$B^{+}_{\delta}$\; we have  
\begin{equation*}
-d^{*}_{\gamma}\lim_{y\to 0} y^{1-2\gamma}\frac{\partial \hat{\delta}_{a,\lambda}(x)}{\partial y}=-d^{*}_{\gamma}\lim_{y\to 0} y^{1-2\gamma}\frac{\partial \bar{\delta}_{0,\lambda}(x)}{\partial y}=\delta_{0,\lambda}^{\frac{1+\gamma}{1-\gamma}}(z),\;\;\;\;z\in \partial B^{+}_{\delta},
\end{equation*}
then  
\begin{equation*}
     -d^{*}_{\gamma}\lim_{y\to 0} y^{1-2\gamma}\;\frac{\partial u_{a,\lambda}(x) }{\partial y}-u^{\frac{1+\gamma}{1-\gamma}}_{a,\lambda}(z)=0,\;\;\;\;\forall\;\;\;z\in \partial B^{+}_{\delta}.
\end{equation*}
Therefore, \eqref{C_022} implies 
\begin{equation*}
   \left|-d^{*}_{\gamma}\lim_{y\to 0} y^{1-2\gamma}\;\frac{\partial u_{a,\lambda}(z,y) }{\partial y}-u^{\frac{1+\gamma}{1-\gamma}}_{a,\lambda}(z)\right|\le C \left(\frac{\lambda}{1+\lambda^{2}d_g(a,z)^{2}}\right)^{1+\gamma}\large{1}_{\{\rho\in M:\;d_g(a,\rho)\geq \delta\}}(z),\;\;\;\;z\in M.
   \end{equation*}
 \end{pf}
%
%
%
%
\section{PS-sequences and Deformation Lemma}\label{PSS}
In this Section, we discuss the asymptotic behavior of Palais-Smale (PS)-sequences for \;$J_q.$\; We also introduce the neighborhoods of potential critical points at infinity  of $J_q$ and their associated selection maps. As in other applications of the Barycenter Technique of Bahri-Coron\cite{bc} (see \cite{aldawood},\cite{aldawood1}, \cite{moa}, and \cite{jfe2}), we also recall the associated Deformation Lemma.
\vspace{4pt}

\noindent
By some arguments which are classical by now see for example (\cite{sma2}, \cite{sma3}, \cite{fcm2}, and \cite{sma}), we have the following the profile decomposition for (PS)-sequences of \;$J_q.$\;
\begin{lem}\label{psseq}
Suppose that \;$(u_k)\subset W^{1,2}_{y^{1-2\gamma},+}(X)$ is a (PS)-sequences for \;$J_q$, that is \;$\n J_q(u_k) \rightarrow 0$\; and \;$J_q(u_k)\rightarrow c$ \;up to a subsequence, and $\oint_{M} u^{\frac{2}{1-\gamma}}_k\;\;dS_g=c^{\frac{1}{\gamma}}$\; for \;$k\in \N^{*},$\; then up to a subsequence, we have have there exists \;$u_{\infty}\geq 0,$\; an integer \;$p\geq 0,$\; a sequence of points \;${a_{i, k}}\in M,\;\;i=1, \cdots, p,$\; and a sequence of positive numbers \;${\l_{i, k}},$\;\;$i=1, \cdots p,$\; such that\\
1)\\
$$
\left\{
\begin{split}
-div_{g}(y^{1-2\gamma}\nabla_g u_{\infty})+qu_{\infty}&=0& \text{in}\;\;&X,\\
-d^{*}_{\gamma}\lim_{y\to0}y^{1-2\gamma}\frac{\partial u_{\infty}}{\partial y}&=u^{\frac{1+\gamma}{1-\gamma}}_{\infty}&\text{on}\;\;& M.
\end{split}
\right.
$$
2)\\
$$
||u_k-u_{\infty}-\sum_{i=1}^p u_{a_{i, k}, \l_{i, k}}||_q\longrightarrow 0,\;\; \text{as}\;\; k \longrightarrow \infty .
$$
3)\\
$$
J_q(u_k)^{\frac{1}{\gamma}}\longrightarrow J_q(u_{\infty})^{\frac{1}{\gamma}}+pS^{\frac{1}{\gamma}},\;\; \text{as}\;\; k \longrightarrow \infty .
$$
4)\\\
For $i\neq j$, 
$$ 
\frac{\l_{i, k}}{\l_{j, k}}+\frac{\l_{j, k}}{\l_{i, k}}+\l_{i, k}\l_{j, k}G^{\frac{-1}{1-\gamma}}(a_{i, k}, a_{j, k})\longrightarrow +\infty,\;\; \text{as}\;\; k \longrightarrow \infty, 
$$
where \;$G$\; is as in \eqref{eqgreen}, and \;$\|\;\|_{q}$\; is as in \eqref{uqsqrt}.\\
\end{lem}
\vspace{4pt}

\noindent
To introduce the neighborhoods of potential critical points at infinity of \;$J_q$, we first fix
\begin{equation}\label{varepsilon0}
\varepsilon_0>0\;\;\; \text{and}\;\;\;\varepsilon_0 \simeq 0.
\end{equation}
 Furthermore, we choose
 \begin{equation}\label{nu0}
 \nu_0>1\;\;\;\text{and}\;\;\;\nu_0\simeq 1.
 \end{equation}
Then for \;$p\in \N^*$,\; and \;$0<\varepsilon\leq \varepsilon_0$, we define \;$V(p, \varepsilon)$\; the \;$(p, \varepsilon)$-neighborhood of potential critical points at infinity of \;$J_q$\; by
\begin{equation*}
\begin{split}
V(p, \varepsilon):=\{u\in W^{1,2}_{y^{1-2\gamma},+}(X):&\;\;\exists\; a_1, \cdots, a_{p}\in M,\;\;\alpha_1, \cdots, \alpha_{p}>0,\;\;\l_1, \cdots,\l_{p}>0,\;\;\l_i\geq \frac{1}{\varepsilon}\;\;\text{for}\;\;i=1\cdots, p,\\ 
&\Vert u-\sum_{i=1}^{p}\alpha_i u_{a_i, \l_i}\Vert_q\leq \varepsilon,\;\;\frac{\alpha_i}{\alpha_j}\leq \nu_0\;\;\text{and}\;\;\varepsilon_{i, j}\leq \varepsilon\;\;\text{for}\;\;i\neq j=1, \cdots, p\}.
\end{split}
\end{equation*}
\vspace{6pt}

\noindent
Concerning the sets \;$V(p, \varepsilon)$, for every \;$p\in \N^*$ \;there exists \;$0<\varepsilon_p\leq\varepsilon_0$\; such that for every \;$0<\varepsilon\leq \varepsilon_p$, we have
\begin{equation}\label{eq:mini}
\begin{cases}
\forall u\in V(p, \varepsilon)\;\; \text{the minimization problem}\;\;\min_{B_{\varepsilon}^{p}}\Vert u-\sum_{i=1}^{p}\alpha_iu_{a_i, \l_i}\Vert_q \\
\text{has a solution }\;(\bar \alpha, A, \bar \l)\in B_{\varepsilon}^{p} , \text{which is unique up to permutations,}
\end{cases}
\end{equation}
where \;$B^{p}_{\varepsilon}$\; is defined as 
\begin{equation*}
\begin{split}
B_{\varepsilon}^{p}:=\{(&\bar\alpha=(\alpha_1, \cdots, \alpha_p),\;\; A=(a_1, \cdots, a_p),\;\; \bar \l=\l_1, \cdots,\l_p))\in \R^{p}_+\times (M)^p\times (0, +\infty)^{p}\\
&\frac{\alpha_i}{\alpha_j}\leq \nu_0\;\; \text{and}\;\; \varepsilon_{i, j}\leq \varepsilon,\;\;i\neq j=1, \cdots, p\}.
\end{split}
\end{equation*}
We define the selection map \;$s_p$\; via
\begin{equation*}
s_{p}: V(p, \varepsilon)\longrightarrow (M)^p/\sigma_p\\
:
u\longrightarrow s_{p}(u)=A
\;\,\text{and} \,\;A\;\;\text{is given by}\;\,\eqref{eq:mini}.
\end{equation*}
To state the Deformation Lemma needed for the application of the algebraic topological argument of Bahri-Coron\cite{bc} for existence, we first  set
\begin{equation}\label{dfwp}
W_p:=\{u\in W^{1,2}_{y^{1-2\gamma},+}(X)\;:\;J_q(u)\leq (p+1)^{\gamma}\mathcal{S}\},
\end{equation}
 for \;$p\in \N,$\; where \;$\mathcal{S}$\; is as in \eqref{S}. 
\vspace{4pt}

\noindent
As in \cite{aldawood},\cite{aldawood1}, \cite{martndia2},\cite{mmc}, and \cite{nss}, we have Lemma \ref{psseq} implies the following Deformation Lemma.
\begin{lem}\label{deform}
Assuming that \;$J_q$ \;has no critical points, then for every \;$p\in \N^*$, up to taking \;$\varepsilon_p$\; given by \eqref{eq:mini} smaller, we have that for every\; $0<\varepsilon\leq \varepsilon_p$, the topological pair \;$(W_p,\; W_{p-1})$\; retracts by deformation onto \;$(W_{p-1}\cup A_p, \;W_{p-1})$\; with \;$V(p, \;\tilde \varepsilon)\subset A_p\subset V(p, \;\varepsilon)$\; where \;$0<\tilde \varepsilon<\frac{\varepsilon}{4}$\; is a very small positive real number and depends on \;$\varepsilon$.
\end{lem}
%
%
%
%
\section{Self-action estimates}\label{SAE}
In this Section, we establish some sharp self-action estimates that are necessary to apply the Bahri-Coron\cite{bc}'s Barycenter technique for existence.  Specifically, we are going to estimate the numerator and the denominator of \;$J_q.$\; To start, we first estimate the numerator of the functional \;$J_q$\; as follows.
\begin{lem}\label{num}
Assuming that \;$\theta>0$\; is small, then there exists \;$C>0$\; such that \;$\forall$\;\;$a \in M,$\;\;$\forall$\;\;$0<2\delta<\delta_0$\; and \;$\forall$\;\;$0<\frac{1}{\l}\le\theta\delta,$\; we have \\
\begin{equation*}
    \begin{split}
     d^{*}_{\gamma}\int_{X}y^{1-2\gamma}\left(\;\left |\nabla_g u_{a,\lambda}\right|^{2} +q u^{2}_{a,\lambda}\right)\;dV_g &\le \oint_{M}u^{\frac{2}{1-\gamma}}_{a,\lambda}\;dS_g +\frac{C}{\lambda^{2-2\gamma}}\left[1+\delta^{2\gamma-1}+\frac{1}{\delta^{2}\lambda^{2\gamma}}\right],   
    \end{split}
\end{equation*}
where \;$\delta_0$\; is as in \eqref{delta0}.
\end{lem}
\begin{pf}
Setting  
\begin{equation*}
    I=d^{*}_{\gamma}\int_{X}y^{1-2\gamma}\left(\;\left|\nabla_g u_{a,\lambda}\right|^{2} +q u^{2}_{a,\lambda}\right)\;dV_g,
\end{equation*}
we have by Green's first identity
\begin{equation*}
\begin{split}
    I=\oint_{M}u^{\frac{2}{1-\gamma}}_{a,\lambda}\;dS_g
    &+\underbrace{d^{*}_{\gamma}\int_{X}\left(\;-div_{g}(y^{1-2\gamma}\nabla_gu_{a,\lambda})+qu_{a,\lambda}\right)u_{a,\lambda}\;dV_g}_{\mbox{$I_1$}}\\
    &+\underbrace{\oint_{M}\left[-d^{*}_{\gamma}\lim_{y\to 0} y^{1-2\gamma}\;\frac{\partial u_{a,\lambda}}{\partial y}-u^{\frac{1+\gamma}{1-\gamma}}_{a,\lambda}\right]u_{a,\lambda}\;dS_g}_{\mbox{$I_2$}}.
    \end{split}
\end{equation*}
Our plan is to estimate \;$I_1$\; and \;$I_2$\; individually. Using Lemma \ref{c0estimate}, we arrive at the following formula for \;$I_1$\;:
\begin{equation*}
\begin{split}
    |I_1|&\le d^{*}_{\gamma}\int_{X}\left|-div_{g}(y^{1-2\gamma}\nabla_gu_{a,\lambda})+qu_{a,\lambda}\right|u_{a,\lambda}\;dV_g\\
    &\le C \left(\frac{1}{\delta^{2}\lambda^{1-\gamma}}+\frac{1}{\delta^{3-2\gamma}\lambda^{1-\gamma}}+\frac{1}{\delta^{6-2\gamma}\lambda^{3-\gamma}}\right)\int_{X} y^{1-2\gamma} \;u_{a,\lambda}\Large{1}_{\{z\in X:\;\delta\le d_g(a,z)\le 2\delta\}}\;dV_g\\
    &+C\;\int_{X} y^{1-2\gamma} \;\hat{\delta}_{a,\lambda}u_{a,\lambda}\large{1}_{\{z\in X:\;d_g(a,z)\le 2\delta\}}\;dV_g.\\
    &+C\;\int_{X} y^{1-2\gamma}\;\left(\frac{\lambda}{1+\lambda^{2}d_g(a,z)^{2}}\right)^{1-\gamma}u_{a,\lambda}\large{1}_{\{z\in X:\;d_g(a,z)\le 2\delta\}}\;dV_g.
    \end{split}
\end{equation*}
For the first term on the right-hand side of the formula above, we observe
\begin{equation}\label{selfI11}
\begin{split}
    \int_{X}y^{1-2\gamma}\;u_{a,\lambda}&\large{1}_{\{z\in X:\;\delta\le d_g(a,z)\le 2\delta\}}\;dV_g(z)\\
    &\le C \int_{\{z\in X:\;\delta\le d_g(a,z)\le 2\delta\}} y^{1-2\gamma}\;\left(\frac{\lambda}{1+\lambda^{2}d_g(a,z)^{2}}\right)^{1-\gamma}\;dV_g(z)\\
    &\le \frac{C}{\lambda^{1-\gamma}}\int_{B^{+}_{2\delta}\setminus B^{+}_{\delta}} y^{1-2\gamma}\;\left[\frac{1}{1+|y|^{2}+|x|^{2}}\right]^{1-\gamma}\;dy\;dx\\
    &\le \frac{C}{\lambda^{1-\gamma}}\int_{B^{+}_{2\delta}\setminus B^{+}_{\delta}} y^{1-2\gamma}\;\left[\frac{1}{1+|y|^{2}}\right]^{\frac{1-\gamma}{2}}\left[\frac{1}{1+|x|^{2}}\right]^{\frac{1-\gamma}{2}}\;dy\;dx\\
    &\le \frac{C}{\lambda^{1-\gamma}}\int_{B^{+}_{2\delta}\setminus B^{+}_{\delta}}\frac{1}{y^{\gamma}}\frac{1}{|x|^{1-\gamma}}\;dy\;dx\\
    &\le \frac{C}{\lambda^{1-\gamma}}\int_{\delta}^{2\delta}\frac{1}{y^{\gamma}}\;dy\int_{B_{2\delta}\setminus B_{\delta}}\frac{1}{|x|^{1-\gamma}}\;dx\\
    &\le C \frac{\delta^{1-\gamma}}{\lambda^{1-\gamma}}\int_{\delta}^{2\delta}r^{\gamma}\;dr\\
    & \le C \frac{\delta^{2}}{\lambda^{1-\gamma}}.
    \end{split}
\end{equation}
For the second term, we get
\begin{equation}\label{selfI12}
    \begin{split}
    \int_{X} y^{1-2\gamma}\;\hat{\delta}_{a,\lambda}u_{a,\lambda}&\large{1}_{\{z\in X:\;d_g(a,z)\le 2\delta\}}\;dV_g(z)\\
    &\le C\int_{\{z\in X:\;d_g(a,z)\le2\delta\}} y^{1-2\gamma}\left(\frac{\lambda}{1+\lambda^{2}d_g(a,z)^{2}}\right)^{2-2\gamma}\;dV_g(z)\\
    &\le \frac{C}{\lambda^{2-2\gamma}}\int_{B^{+}_{2\delta}} y^{1-2\gamma}\;\left[\frac{1}{1+|y|^{2}+|x|^{2}}\right]^{2-2\gamma}\;dy\;dx\\
    &\le \frac{C}{\lambda^{2-2\gamma}}\int_{B^{+}_{2\delta}} y^{1-2\gamma}\;\left[\frac{1}{1+|y|^{2}}\right]^{\frac{2-2\gamma}{2}}\left[\frac{1}{1+|x|^{2}}\right]^{\frac{2-2\gamma}{2}}\;dy\;dx\\
    &\le \frac{C}{\lambda^{2-2\gamma}}\int_{B^{+}_{2\delta}}\frac{1}{y}\frac{1}{|x|^{2-2\gamma}}\;dy\;dx\\
    &\le \frac{C}{\lambda^{2-2\gamma}}\int_{0}^{2\delta}\frac{1}{y}\;dy\int_{B_{2\delta}}\frac{1}{|x|^{2-2\gamma}}\;dx\\
    &\le \frac{C}{\lambda^{2-2\gamma}}\int_{0}^{2\delta}r^{2\gamma-1}\;dr\\
    & \le C \frac{\delta^{2\gamma}}{\lambda^{2-2\gamma}}.
    \end{split}
    \end{equation}
Now, for the final term, we have
\begin{equation}\label{selfI13}
    \begin{split}
    \int_{X} y^{1-2\gamma}\left(\frac{\lambda}{1+\lambda^{2}d_g(a,z)^{2}}\right)^{1-\gamma}u_{a,\lambda}&\large{1}_{\{z\in X:\;d_g(a,z)\le 2\delta\}}\;dV_g(z)\\
    &\le C\int_{\{z\in X:\;d_g(a,z)\le2\delta\}} y^{1-2\gamma}\left(\frac{\lambda}{1+\lambda^{2}d_g(a,z)^{2}}\right)^{2-2\gamma}\;dV_g(z)\\
    &\le \frac{C}{\lambda^{2-2\gamma}}\int_{B^{+}_{2\delta}} y^{1-2\gamma}\;\left[\frac{1}{1+|y|^{2}+|x|^{2}}\right]^{2-2\gamma}\;dy\;dx\\
    &\le \frac{C}{\lambda^{2-2\gamma}}\int_{B^{+}_{2\delta}} y^{1-2\gamma}\;\left[\frac{1}{1+|y|^{2}}\right]^{\frac{2-2\gamma}{2}}\left[\frac{1}{1+|x|^{2}}\right]^{\frac{2-2\gamma}{2}}\;dy\;dx\\
    &\le \frac{C}{\lambda^{2-2\gamma}}\int_{B^{+}_{2\delta}}\frac{1}{y|x|^{2-2\gamma}}\;dy\;dx\\
    &\le \frac{C}{\lambda^{2-2\gamma}}\int_{0}^{2\delta}\frac{1}{y}\;dy\int_{B_{2\delta}}\frac{1}{|x|^{2-2\gamma}}\;dx\\
    &\le \frac{C}{\lambda^{2-2\gamma}}\int_{0}^{2\delta}r^{2\gamma-1}\;dr\\
    & \le C \frac{\delta^{2\gamma}}{\lambda^{2-2\gamma}}.
    \end{split}
    \end{equation}
Thus, combining \eqref{selfI11}-\eqref{selfI13}, we get
\begin{equation}\label{selfI1}
    |I_1|\le \frac{C}{\lambda^{2-2\gamma}}\left(1+\delta^{2\gamma}+\delta^{2\gamma-1}+\frac{1}{\delta^{4-2\gamma}\lambda^{2}}\right).
\end{equation}
Next, in the case of \;$I_2,$\; we have
\begin{equation*}
    \begin{split}
        |I_2|&\le \oint_{M}\left|-d^{*}_{\gamma}\lim_{y\to 0} y^{1-2\gamma}\;\frac{\partial u_{a,\lambda}}{\partial y}-u^{\frac{1+\gamma}{1-\gamma}}_{a,\lambda}\right|u_{a,\lambda}\;dS_g\\
        &\le C\oint_{ M}\left(\frac{\lambda}{1+\lambda^{2}d_g(a,x)^{2}}\right)^{1+\gamma}u_{a,\lambda}\Large{1}_{\{x\in M:\;d_g(a,x)\geq \delta\}}\;dS_g.
    \end{split}
\end{equation*}
On the right-hand side of the above formula, we observe
\begin{equation}
    \begin{split}
        \oint_{ M}\left(\frac{\lambda}{1+\lambda^{2}d_g(a,x)^{2}}\right)^{1+\gamma} u_{a,\lambda}\Large{1}_{\{x\in M:\;d_g(a,x)\geq \delta\}}\;dS_g&\le C\oint_{\{x\in M:\;d_g(a,x)\geq \delta\}}\left(\frac{\lambda}{1+\lambda^{2}d_g(a,x)^{2}}\right)^{2}\;dS_g\\
        &\le C\oint_{\{x\in M:\;d_g(a,x)\geq \delta\}}\left(\frac{1}{\lambda d_g(a,x)^{2}}\right)^{2}\;dS_g\\
        &\le \frac{C}{\lambda^{2}}\left[\int_{
        B_{\delta_{0}}\setminus B_{\delta}}\frac{1}{|x|^{4}}\;dx+1\right]\\
        &\le \frac{C}{\lambda^{2}}\left[\int_{\R^{2}\setminus B_{\delta}}^{+\infty}r^{-3}\;dr+1\right]\\
        &\le \frac{C}{\lambda^{2}}\left[\int_{\delta}^{+\infty}r^{-3}\;dr+1\right]\\
        &\le C \frac{1}{\delta^{2}\lambda^{2}}.
    \end{split}
\end{equation}
Thus, we have
\begin{equation}\label{selfI2}
    |I_2|\le \frac{C}{\delta^{2}\lambda^{2}}.
\end{equation}
As a result, by combining \eqref{selfI1} and \eqref{selfI2}, we obtain
\begin{equation*}
    \begin{split}
     d^{*}_{\gamma}\int_{X}y^{1-2\gamma}\left(\left |\nabla_g u_{a,\lambda}\right|^{2} +q u^{2}_{a,\lambda}\right)\;dV_g &\le
     \oint_{M}u^{\frac{2}{1-\gamma}}_{a,\lambda}\;dS_g+\frac{C}{\lambda^{2-2\gamma}}\left[1+\delta^{2\gamma-1}+\frac{1}{\delta^{2}\lambda^{2\gamma}}\right].   
    \end{split}
\end{equation*}
\end{pf}\\
\noindent
For the denominator of \;$J_q$, we have.
\begin{lem}\label{denom}
 Assuming that \;$\theta>0$\; is small, then there exists \;$C>0$\; such that \;$\forall$\;\;$a \in M,$\;\;$\forall$\;\;$0<2\delta<\delta_0$\; and \;$\forall$\;\;$0<\frac{1}{\l}\le \theta\delta,$\; we have
\[\oint_{M} u^{\frac{2}{1-\gamma}}_{a,\lambda}\;dS_g=\int_{\R^2} \delta^{\frac{2}{1-\gamma}}_{0,\lambda}\;dx+O\left(\frac{1}{\delta^2\lambda^2}\right),\]
where \;$\delta_0$\; is as in \eqref{delta0}.
\end{lem}

\begin{pf}
We have
\begin{equation}\label{denom1}
\begin{split}
\oint_{M} u^{\frac{2}{1-\gamma}}_{a,\lambda}\;dS_g&=\oint_{\{x\in M:\;d_g(a,x)\le \delta\}} u^{\frac{2}{1-\gamma}}_{a,\lambda}\;dS_g+\oint_{\{x\in M:\;\delta<d_g(a,x)\le 2\delta\}} u^{\frac{2}{1-\gamma}}_{a,\lambda}\;dS_g\\
&+\oint_{\{x\in M:\;d_g(a,x)>2\delta\}} u^{\frac{2}{1-\gamma}}_{a,\lambda}\;dS_g.
\end{split}
\end{equation}
Now, we estimate each of the terms on the right side of the formula for \eqref{denom1}. We will begin by estimating the first term, and we have the following:
\begin{equation}\label{denom2}
\begin{split}
\oint_{\{x\in M:\;d_g(a,x)\le \delta\}} u^{\frac{2}{1-\gamma}}_{a,\lambda}\;dS_g&=\oint_{\{x\in M:\;d_g(a,x)\le \delta\}} \delta^{\frac{2}{1-\gamma}}_{a,\lambda}\;dS_g\\
&=\int_{B_{\delta}}\delta^{\frac{2}{1-\gamma}}_{0,\lambda}\;dx\\
&=\int_{\R^2}\delta^{\frac{2}{1-\gamma}}_{0,\lambda}\;dx-\int_{\{x\in \R^{2}:\;|x|\geq\delta\}} \delta^{\frac{2}{1-\gamma}}_{0,\lambda}\;dx\\
&=\int_{\R^2}\delta^{\frac{2}{1-\gamma}}_{0,\lambda}\;dx-\frac{1}{\lambda^{2}}\int_{\{x\in \R^{2}:\;|x|\geq\delta\}} \frac{1}{|x|^{4}}\;dx\\
&=\int_{\R^2}\delta^{\frac{2}{1-\gamma}}_{0,\lambda}\;dx+O\left(\frac{1}{\delta^2\lambda^2}\right).
\end{split}
\end{equation}
In order to derive the second term from the formula for \eqref{denom1}, we have
\begin{equation}\label{denom3}
\begin{split}
\oint_{\{x\in M:\;\delta<d_g(a,x)\le 2\delta\}} u^{\frac{2}{1-\gamma}}_{a,\lambda}\;dS_g&\le C\oint_{\{x\in M:\;\delta<d_g(a,x)\le 2\delta\}}\left(\frac{\lambda}{1+\lambda^2d_g(a,x)^2}\right)^{2}\;dS_g\\
&\le C \oint_{\{x\in M:\;\delta<d_g(a,x)\le 2\delta\}}\frac{1}{\lambda^{2} d_g(a,x)^{4}}\;dS_g\\ 
&\le \frac{C}{\lambda^{2}}\int_{B_{2\delta}\setminus B_{\delta}}\frac{1}{|x|^{4}}\;dx\\
&\le \frac{C}{\lambda^2}\int_{\delta}^{2\delta} r^{-3}\;dr\\&\le \frac{C}{\delta^2\lambda^2}.
\end{split}
\end{equation}
\vspace{4pt}

\noindent
Next, using \eqref{estg}, we estimate the remaining term as follows 
\begin{equation}\label{denom4}
\begin{split}
\oint_{\{x\in M:\;d_g(a,x)>2\delta\}} u^{\frac{2}{1-\gamma}}_{a,\lambda}\;dS_g&=\oint_{\{x\in M:\;d_g(a,x)>2\delta\}} \left(\frac{c_0}{\lambda^{1-\gamma}}G_a(x)\right)^{\frac{2}{1-\gamma}}dS_g\\&=\frac{C}{\lambda^2}\oint_{\{x\in M:\;d_g(a,x)>2\delta\}} G^{\frac{2}{1-\gamma}}_a(x)\;dS_g\\
&\le C \oint_{\{x\in M:\;d_g(a,x)>2\delta\}} \frac{1}{\lambda^{2}d_g(a,x)^{4}}\;dS_g\\
&\le C\oint_{\{x\in M:\;2\delta< d_g(a,x)\le \delta_{0}\}}\frac{1}{\lambda^{2}d_g(a,x)^{4}}\;dS_g\\
&+C\oint_{\{x\in M:\;d_g(a,x)>2\delta_{0}\}}\frac{1}{\lambda^{2}d_g(a,x)^{4}}\;dS_g\\
&\le \frac{C}{\lambda^{2}}\left[\oint_{\{x\in M:\;2\delta< d_g(a,x)\le \delta_{0}\}}\frac{1}{d_g(a,x)^{4}}\;dS_g+C\right]\\
&\le \frac{C}{\lambda^{2}}\left[\int_{B_{\delta_{0}}\setminus B_{2\delta}}\frac{1}{|x|^{4}}\;dx+C\right]\\
&\le \frac{C}{\lambda^{2}}\left[\int_{2\delta}^{+\infty} r^{-3}\;dr+C\right]\\
&\le\frac{C}{\delta^2\lambda^2}.
\end{split}
\end{equation}
Finally, combining, \eqref{denom1}-\eqref{denom4}, we have \\\
\[\oint_{M} u^{\frac{2}{1-\gamma}}_{a,\lambda}\;dS_g=\int_{\R^2} \delta^{\frac{2}{1-\gamma}}_{0,\lambda}\;dx+O\left(\frac{1}{\delta^2\lambda^2}\right).\]
\end{pf}
\vspace{4pt}

\noindent
In the next Corollary, we establish now the $J_q$-energy estimate of \;$u_{a, \lambda}$\; needed for the application of the Barycenter Technique
of Bahri-Coron\cite{bc}  for existence.
\vspace{4pt}
\begin{cor}\label{sharpenergy}
Assuming that \;$\theta>0$\; is small, then there exists \;$C>0$\; such that \;$\forall$\;\;$a\in M,$\;\;$\forall$\;\;$0<2\delta<\delta_0$\; and \;$\forall$\;\;$0<\frac{1}{\l}\le \theta\delta,$\; we have

\[J_q(u_{a, \l})\leq \mathcal{S}\left( 1+\frac{C}{\lambda^{2-2\gamma}}\left[1+\delta^{2\gamma-1}+\frac{1}{\delta^{2}\lambda^{2\gamma}}\right]\right),\]
where \;$\delta_0$\; is as in \eqref{delta0}.
\end{cor}
\begin{pf}
It follows from the properties of \;$\d_{0, \l}$ (see \eqref{S3} \eqref{S}), Lemma \ref{num}, and Lemma \ref{denom}.
\end{pf}
%
%
%
%
\section{Interaction estimates}\label{IE}
Throughout this Section, we derive some sharp inter-action estimates that are necessary for the algebraic topological argument of Bahri-Coron\cite{bc} for existence . We start with the following technical inter-action estimates. We remark that it will provide the required inter-action estimates between \;$e_{ij}$\; and \;$\epsilon_{ji},$\; see \eqref{epij} and \eqref{eij} for the definitions of \;$e_{ij}$\; and \;$\epsilon_{ji}.$
\begin{lem}\label{interact1}
Assuming that \;$\theta>0$ is small, then there exists \;$C>0$\; such that \;$\forall$\;\;$\ a_i, a_j\in M,$\; \;$\forall$\;\;$0<2\delta<\delta_0,$\; and \;$\forall$\;\;$0<\frac{1}{\l_j},$\;\;$\frac{1}{\l_i}\leq \theta\delta,$\; we have
\begin{equation*}
    \begin{split}
     & d^{*}_{\gamma}\int_{X}\left|\left(-div_{g}(y^{1-2\gamma}\nabla_g)+q\right)u_{a_j,\lambda_j}\right|u_{a_i,\lambda_i}\;dV_g+\oint_{M}\left|-d^{*}_{\gamma}\lim_{y\to 0} y^{1-2\gamma}\;\frac{\partial u_{a_j,\lambda_j}}{\partial y}-u^{\frac{1+\gamma}{1-\gamma}}_{a_j,\lambda_j}\right|u_{a_i,\lambda_i}\;dS_g\\
     &\le C\left[\delta+\delta^{2-2\gamma}+\frac{1}{\delta^{2}\lambda^{2}_j}+\frac{1}{\delta^{4\gamma}\lambda_j^{2\gamma}}\right]\left(\frac{\lambda_{i}}{\lambda_{j}}+\lambda_i\lambda_j d_g(a_{i},a_{j})^{2}\right)^{-(1-\gamma)},
    \end{split}
\end{equation*}
where \;$\delta_0$\; is as in \eqref{delta0}.
\end{lem}
\begin{pf}
By using Lemma \ref{c0estimate}, we arrive at the following results:

\begin{equation}\label{A_j}
\begin{split}
    \underbrace{\left|\left(-div_{g}(y^{1-2\gamma}\nabla_g)+q\right)u_{a_j,\lambda_j}(x)\right|}_{\mbox{$A_j$}}\le & C\;y^{1-2\gamma}\left[\left(\frac{1}{\delta^{2}\lambda^{1-\gamma}_{j}}+\frac{1}{\delta^{3-2\gamma}\lambda^{1-\gamma}_{j}}+\frac{1}{\delta^{6-2\gamma}\lambda^{3-\gamma}_{j}}\right) \large{1}_{\{z\in X:\;\delta\le d_g(a_{j},z)\le 2\delta\}}(x)\right.\\
  +& \hat{\delta}_{a_j,\lambda_j} \large{1}_{\{z\in X:\;d_g(a_{j},z)\le 2\delta\}}(x)\\ 
  +&\left.\left(\frac{\lambda_{j}}{1+\lambda^{2}d_g(a_{j},x)^{2}}\right)^{1-\gamma}\large{1}_{\{z\in X:\;d_g(a_{j},z)\le 2\delta\}}(x)\right],\;\;x\in X,
  \end{split}
\end{equation}
and
\begin{equation}\label{B_j}
    \underbrace{\left|-d^{*}_{\gamma}\lim_{y\to 0} y^{1-2\gamma}\;\frac{\partial u_{a_j,\lambda_j}(z,y)}{\partial y}-u^{\frac{1+\gamma}{1-\gamma}}_{a_j,\lambda_j}(z)\right|}_{\mbox{$B_j$}}\le C\left[\left(\frac{\lambda_j}{1+\lambda^{2}_{j}d_g(a_j,z)^{2}}\right)^{1+\gamma}\Large{1}_{\{\rho\in M:\;d_g(a_j,\rho)\geq \delta\}}(z)\right],\;\;z\in M.
\end{equation}
For \;$x\in\{z\in X:\;\;d_g(a_j,z)\leq 2\d\},$ we have 
$$\left(\frac{\lambda_j}{1+\lambda^{2}_{j}d_g(a_j,x)^{2}}\right)\geq\left(\frac{\lambda_j}{1+4\lambda^2_j \delta^2}\right)\geq\frac{1}{\lambda_j \delta^2}\left[1+O(\frac{1}{\lambda^2_j \delta^2})\right]\geq\frac{1}{2 \lambda_j\delta^2}.$$
This implies $$\frac{1}{\lambda_j^{1-\gamma} \delta^{2-2\gamma}}\leq C \left(\frac{\lambda_j}{1+\lambda^{2}_{j}d_g(a_j,x)^{2}}\right)^{1-\gamma}.$$\\
Thus, using \eqref{A_j} and Lemma \ref{Appendix}, we obtain
\begin{equation}\label{A_j1}
\begin{split}
    A_j&\le C\;y^{1-2\gamma}\;\left(\frac{1}{\delta^{2\gamma}}+\frac{1}{\delta}+\frac{1}{\delta^{4}\lambda^{2}_j}+1\right) \left(\frac{\lambda_j}{1+\lambda^{2}_{j}d_g(a_j,x)^{2}}\right)^{1-\gamma}\large{1}_{\{z\in X:\;d_g(a_j,z)\le 4\delta\}}(x),\;\;x\in X.
    \end{split}
\end{equation}
For \;$B_j$\; given by \eqref{B_j}, we have
\begin{equation}\label{B_j1}
    B_j\le C \left(\frac{\lambda_j}{1+\lambda^{2}_{j}d_g(a_j,z)^{2}}\right)^{1+\gamma}\large{1}_{\{\rho\in M:\;d_g(a_j,\rho)\geq \delta\}}(z),\;\;z\in M.
\end{equation}
Now, using \eqref{A_j1}, we get 
\begin{equation}
\begin{split}
    &\int_{X} A_j\;u_{a_i},_{\lambda_i}\;dV_g\\
    &\le C \left(\frac{1}{\delta^{2\gamma}}+\frac{1}{\delta}+\frac{1}{\delta^{4}\lambda_j^{2}}+1\right) \underbrace{\int_{\{x\in X:\;d_g(a_j,x)\le4\delta\}}y^{1-2\gamma}\left(\frac{\lambda_j}{1+\lambda_j^2d_g(a_j,x)^{2}}\right)^{1-\gamma}\left(\frac{\lambda_i}{1+\lambda_i^2d_g(a_i,x)^{2}}\right)^{1-\gamma}\;dV_g}_{\mbox{$I_1$}}.
\end{split}
\end{equation}
Also, using \eqref{B_j1}, we have 
\begin{equation}
\oint_{M} B_j\;u_{a_i},_{\lambda_i}\;dS_g\le  C\underbrace{\oint_{\{z\in M:\;d_g(a_j,z)\geq\delta\}}\left(\frac{\lambda_j}{1+\lambda_j^2d_g(a_j,z)^{2}}\right)^{1+\gamma}\left(\frac{\lambda_i}{1+\lambda_i^2d_g(a_i,z)^{2}}\right)^{1-\gamma}\;dS_g}_{\mbox{$I_2$}}.
\end{equation}
Denote
\begin{equation*}
    \begin{split}
        &\mathcal{B}={\{x\in X:\;2d_g(a_i,x)\le\frac{1}{\lambda_j}+d_g(a_i,a_j)\}\cap\{x\in X:\;d_g(a_j,x)\le 4\delta\}},\\
        &\mathcal{B}^{c}={\{x\in X:\;2d_g(a_i,x)>\frac{1}{\lambda_j}+d_g(a_i,a_j)\}\cap\{x\in X:\;d_g(a_j,x)\le 4\delta\}}.
    \end{split}
\end{equation*}
In the next step, we are going to estimate \;$I_1$\;as follows
\begin{equation*}
    \begin{split}
        I_1&=\int_{\{x\in X:\;d_g(a_j,x)\le4\delta\}}y^{1-2\gamma}\left(\frac{\lambda_j}{1+\lambda_j^2d_g(a_j,x)^{2}}\right)^{1-\gamma}\left(\frac{\lambda_i}{1+\lambda_i^2d_g(a_i,x)^{2}}\right)^{1-\gamma}\;dV_g\\
        &=\underbrace{\int_{\mathcal{B}}y^{1-2\gamma}\left(\frac{\lambda_j}{1+\lambda_j^2d_g(a_j,x)^{2}}\right)^{1-\gamma}\left(\frac{\lambda_i}{1+\lambda_i^2d_g(a_i,x)^{2}}\right)^{1-\gamma}}_{\mbox{$I^{1}_1$}}\;dV_g\\
        &+\underbrace{\int_{\mathcal{B}^{c}}y^{1-2\gamma}\left(\frac{\lambda_j}{1+\lambda_j^2d_g(a_j,x)^{2}}\right)^{1-\gamma}\left(\frac{\lambda_i}{1+\lambda_i^2d_g(a_i,x)^{2}}\right)^{1-\gamma}}_{\mbox{$I^{2}_1$}}\;dV_g.
    \end{split}
\end{equation*}
To estimate \;$I^{1}_1,$\; we use the triangle inequality, and we have
\begin{equation*}
\begin{split}
    I^{1}_1&\le C\int_{\{z\in X:\;d_g(a_i,z)\le 8\delta\}}y^{1-2\gamma}\left(\frac{\lambda_j}{1+\lambda_j^2d_g(a_j,a_i)^{2}}\right)^{1-\gamma}\left(\frac{\lambda_i}{1+\lambda_i^2d_g(a_i,z)^{2}}\right)^{1-\gamma}\;dV_g(z)\\
    &\le C\frac{\lambda_j^{1-\gamma}}{\left(1+\lambda^2_jd_g(a_j,a_i)^{2}\right)^{1-\gamma}}\int_{\{z\in X:\;d_g(a_i,z)\le 8\delta\}}y^{1-2\gamma}\left(\frac{\lambda_i}{1+\lambda_i^2d_g(a_i,z)^{2}}\right)^{1-\gamma}\;dV_g(z)\\
    &\le C\frac{\lambda_j^{1-\gamma}}{\left(1+\lambda^2_jd_g(a_j,a_i)^{2}\right)^{1-\gamma}}\left(\frac{1}{\lambda_{i}}\right)^{1-\gamma}\int_{B^{+}_{8\delta}} y^{1-2\gamma}\left[\frac{1}{1+|y|^{2}+|x|^{2}}\right]^{1-\gamma}\;dy\;dx\\
    &\le C\frac{\lambda_j^{1-\gamma}}{\left(1+\lambda^2_jd_g(a_j,a_i)^{2}\right)^{1-\gamma}}\left(\frac{1}{\lambda_{i}}\right)^{1-\gamma}\int_{B^{+}_{8\delta}} y^{1-2\gamma}\left[\frac{1}{1+|y|^{2}}\right]^{\frac{1-\gamma}{2}}\left[\frac{1}{1+|x|^{2}}\right]^{\frac{1-\gamma}{2}}\;dy\;dx\\
    &\le C\frac{\lambda_j^{1-\gamma}}{\left(1+\lambda^2_jd_g(a_j,a_i)^{2}\right)^{1-\gamma}}\left(\frac{1}{\lambda_{i}}\right)^{1-\gamma}\int_{B^{+}_{8\delta}}\frac{1}{y^{\gamma}}\frac{1}{|x|^{1-\gamma}}\;dy\;dx\\
     &\le C\frac{\lambda_j^{1-\gamma}}{\left(1+\lambda^2_jd_g(a_j,a_i)^{2}\right)^{1-\gamma}}\left(\frac{1}{\lambda_{i}}\right)^{1-\gamma}\int_{0}^{8\delta}\frac{1}{y^{\gamma}}\;dy\int_{B_{8\delta}}\frac{1}{|x|^{1-\gamma}}\;dx\\
      &\le C\frac{\left(\frac{\lambda_j}{\lambda_i}\right)^{1-\gamma}\delta^{1-\gamma}}{\left(1+\lambda^2_jd_g(a_j,a_i)^{2}\right)^{1-\gamma}}\int_{B_{8\delta}}\frac{1}{|x|^{1-\gamma}}\;dx.
\end{split}
\end{equation*}
So, for \;$I^{1}_1,$\; we have
\begin{equation}\label{I1J1A_j}
    I^{1}_1=O\left(\delta^{2}\left(\frac{\lambda_{i}}{\lambda_{j}}+\lambda_i\lambda_j d_g(a_{i},a_{j})^{2}\right)^{-(1-\gamma)}\right).
\end{equation}
For \;$I^{2}_1,$\; we observe 
\begin{equation*}
\begin{split}
    I^{2}_1&\le C\left(\frac{\lambda_j}{\lambda_i}\right)^{2-2\gamma}\int_{\{z\in X:\;d_g(a_j,z)\le 4\delta\}}y^{1-2\gamma}\left(\frac{\lambda_j}{1+\lambda_j^2d_g(a_j,z)^{2}}\right)^{1-\gamma}\left(\frac{\lambda_i}{1+\lambda_j^2d_g(a_j,a_i)^{2}}\right)^{1-\gamma}\;dV_g(z)\\
    &\le C\frac{\left(\frac{\lambda_j}{\lambda_i}\right)^{2-2\gamma}\lambda_i^{1-\gamma}}{\left(1+\lambda^2_jd_g(a_j,a_i)^{2}\right)^{1-\gamma}}\int_{\{z\in X:\;d_g(a_j,z)\le 8\delta\}}y^{1-2\gamma}\left(\frac{\lambda_j}{1+\lambda_j^2d_g(a_j,z)^{2}}\right)^{1-\gamma}\;dV_g(z)\\
    &\le C\frac{\left(\frac{\lambda_j}{\lambda_i}\right)^{2-2\gamma}\lambda_i^{1-\gamma}}{\left(1+\lambda^2_jd_g(a_j,a_i)^{2}\right)^{1-\gamma}}\left(\frac{1}{\lambda_{j}}\right)^{1-\gamma}\int_{B^{+}_{8\delta}} y^{1-2\gamma}\left[\frac{1}{1+|y|^{2}+|x|^{2}}\right]^{1-\gamma}\;dy\;dx\\
    &\le C\frac{\left(\frac{\lambda_j}{\lambda_i}\right)^{2-2\gamma}\lambda_i^{1-\gamma}}{\left(1+\lambda^2_jd_g(a_j,a_i)^{2}\right)^{1-\gamma}}\left(\frac{1}{\lambda_{j}}\right)^{1-\gamma}\int_{B^{+}_{8\delta}} y^{1-2\gamma}\left[\frac{1}{1+|y|^{2}}\right]^{\frac{1-\gamma}{2}}\left[\frac{1}{1+|x|^{2}}\right]^{\frac{1-\gamma}{2}}\;dy\;dx\\
    &\le C\frac{\left(\frac{\lambda_j}{\lambda_i}\right)^{2-2\gamma}\lambda_i^{1-\gamma}}{\left(1+\lambda^2_jd_g(a_j,a_i)^{2}\right)^{1-\gamma}}\left(\frac{1}{\lambda_{j}}\right)^{1-\gamma}\int_{B^{+}_{8\delta}}\frac{1}{y^{\gamma}}\frac{1}{|x|^{1-\gamma}}\;dy\;dx\\
     &\le C\frac{\left(\frac{\lambda_j}{\lambda_i}\right)^{2-2\gamma}\lambda_i^{1-\gamma}}{\left(1+\lambda^2_jd_g(a_j,a_i)^{2}\right)^{1-\gamma}}\left(\frac{1}{\lambda_{j}}\right)^{1-\gamma}\int_{0}^{8\delta}\frac{1}{y^{\gamma}}\;dy\int_{B_{8\delta}}\frac{1}{|x|^{1-\gamma}}\;dx\\
      &\le C\frac{\left(\frac{\lambda_j}{\lambda_i}\right)^{1-\gamma}\delta^{1-\gamma}}{\left(1+\lambda^2_jd_g(a_j,a_i)^{2}\right)^{1-\gamma}}\int_{B_{8\delta}}\frac{1}{|x|^{1-\gamma}}\;dx.
\end{split}
\end{equation*}
Thus, we get for \;$I^{2}_1$\;
\begin{equation}\label{I1J2A_j}
    I^{2}_1=O\left(\delta^{2}\left(\frac{\lambda_{i}}{\lambda_{j}}+\lambda_i\lambda_j d_g(a_{i},a_{j})^{2}\right)^{-(1-\gamma)}\right).
\end{equation}
Hence, using \eqref{I1J1A_j} and \eqref{I1J2A_j}, we obtain
\begin{equation}\label{IJ1A_j}
     I_1=O\left(\delta^{2}\left(\frac{\lambda_{i}}{\lambda_{j}}+\lambda_i\lambda_j d_g(a_{i},a_{j})^{2}\right)^{-(1-\gamma)}\right).
\end{equation}
\vspace{6pt}

\noindent
For the last step, we are going to estimate \;$I_2.$\; In order to estimate \;$I_2,$\; we first write the following 
\begin{equation*}
    \begin{split}
        I_2&=\oint_{\{z\in M:\;d_g(a_j,z)\geq\delta\}}\left(\frac{\lambda_j}{1+\lambda_j^2d_g(a_j,z)^{2}}\right)^{1+\gamma}\left(\frac{\lambda_i}{1+\lambda_i^2d_g(a_i,z)^{2}}\right)^{1-\gamma}\;dS_g\\
        &=\underbrace{\oint_{\{z\in M:\;2d_g(a_i,z)\le\frac{1}{\lambda_j}+d_g(a_j,a_i\})\cap\{z\in M:\;d_g(a_j,z)\geq\delta)\}}\left(\frac{\lambda_j}{1+\lambda_j^2d_g(a_j,z)^{2}}\right)^{1+\gamma}\left(\frac{\lambda_i}{1+\lambda_i^2d_g(a_i,z)^{2}}\right)^{1-\gamma}}_{\mbox{$I^{1}_2$}}\;dS_g\\
        &+\underbrace{\oint_{\{z\in M:\;2d_g(a_i,z)>\frac{1}{\lambda_j}+d_g(a_j,a_i)\}\cap\{z\in M:\;d_g(a_j,z)\geq\delta\}}\left(\frac{\lambda_j}{1+\lambda_j^2d_g(a_j,z)^{2}}\right)^{1+\gamma}\left(\frac{\lambda_i}{1+\lambda_i^2d_g(a_i,z)^{2}}\right)^{1-\gamma}}_{\mbox{$I^{2}_2$}}\;dS_g.
    \end{split}
\end{equation*}
\noindent
To estimate \;$I^{1}_2,$\; we set
\begin{equation*}
    \begin{split}
    \mathcal{A}&=\{z\in M:2d_g(a_i,z)\le\frac{1}{\lambda_j}+d_g(a_j,a_i)\}\cap\{z\in M:d_g(a_j,z)\geq\delta\},\\
    r_{ij}&=\frac{1}{2}\left(\frac{1}{\lambda_j}+d_g(a_j,a_i)\right),
    \end{split}
\end{equation*}
and observe
\begin{equation*}
\begin{split}
I^{1}_2&\le \frac{C}{ \lambda_j^{-(1+\gamma)}}\oint_{\mathcal{A}}\left(\frac{1}{1+\lambda_j^2d_g(a_j,z)^{2}}\right)^{1-\gamma}\left(\frac{\lambda_i}{1+\lambda_i^2d_g(a_i,z)^{2}}\right)^{1-\gamma}\left(\frac{1}{1+\lambda_j^2d_g(a_j,z)^{2}}\right)^{2\gamma}\;dS_g\\
&\le C\left(\frac{1}{\lambda_i^{1-\gamma}}\right)\left(\frac{1}{\lambda_j^{4\gamma} \delta^{4\gamma}}\right) \frac{\lambda_j^{1+\gamma}}{\left(1+\lambda_j^2d_g(a_j,a_i)^{2}\right)^{1-\gamma}}\oint_{\{z\in M:\;2d_g(a_i,z)\le r_{ij}\}}\frac{1}{d_g(a_i,z)^{2-2\gamma}}\;dS_g\\
&\le C\left(\frac{1}{\lambda_i^{1-\gamma}}\right)\left(\frac{1}{\lambda_j^{4\gamma} \delta^{4\gamma}}\right) \frac{\lambda_j^{1+\gamma}}{\left(1+\lambda_j^2d_g(a_j,a_i)^{2}\right)^{1-\gamma}}\int_{B_{r_{ij}}}\frac{1}{|z|^{2-2\gamma}}\;dz\\
&\le C \left(\frac{1}{\lambda_i}\right)^{1-\gamma}\left(\frac{1}{\lambda_j^{4\gamma}\delta^{4\gamma}}\right) \frac{\lambda^{1+\gamma}_{j}}{\left(1+\lambda_j^2d_g(a_j,a_i)^{2}\right)^{1-\gamma}}\left(\frac{1}{\lambda_j}+d_g(a_j,a_i)\right)^{2\gamma}\\
&\le C \left(\frac{\lambda_j}{\lambda_i}\right)^{1-\gamma}\left(\frac{1}{\lambda_j^{2\gamma}\delta^{4\gamma}}\right) \frac{1}{\left(1+\lambda_j^2d_g(a_j,a_i)^{2}\right)^{1-\gamma}}.
\end{split}
\end{equation*}
This implies 
\begin{equation}\label{I12B_j}
I^{1}_2=O\left(\frac{1}{\delta^{4\gamma}\lambda^{2\gamma}_j}\left(\frac{\lambda_{i}}{\lambda_{j}}+\lambda_i\lambda_j d_g(a_{i},a_{j})^{2}\right)^{-(1-\gamma)}\right).
\end{equation}
For \;$I^{2}_2,$\; we derive
\begin{equation*}
\begin{split}
I^{2}_2&\le C\frac{\lambda_j^{2-2\gamma}}{\lambda_i^{1-\gamma}}\oint_{\{z\in M:\;d_g(a_{j},z)\geq\delta\}}\left(\frac{1}{1+\lambda_j^2d_g(a_{i},a_{j})^{2}}\right)^{1-\gamma}\left(\frac{\lambda_j}{1+\lambda_j^2d_g(a_{j},z)^{2}}\right)^{1+\gamma}\;dS_g\\
&\le C\frac{\lambda_j^{2-2\gamma}}{\lambda_i^{1-\gamma}}\left(\frac{1}{1+\lambda_j^2d_g(a_{i},a_{j})^{2}}\right)^{1-\gamma}\frac{1}{\lambda_j^{1+\gamma}}\oint_{\{z\in M:\;d_g(a_{j},z)\geq\delta\}}\frac{1}{d_g(a_{j},z)^{2+2\gamma}}\;dS_g\\
&\le C\frac{\lambda_j^{2-2\gamma}}{\lambda_i^{1-\gamma}}\left(\frac{1}{1+\lambda_j^2d_g(a_{i},a_{j})^{2}}\right)^{1-\gamma}\frac{1}{\lambda_j^{1+\gamma}}\oint_{\{z\in M:\;\delta \le d_g(a_{j},z)\le \delta_{0}\}}\frac{1}{d_g(a_{j},z)^{2+2\gamma}}\;dS_g\\
&+ C\frac{\lambda_j^{2-2\gamma}}{\lambda_i^{1-\gamma}}\left(\frac{1}{1+\lambda_j^2d_g(a_{i},a_{j})^{2}}\right)^{1-\gamma}\frac{1}{\lambda_j^{1+\gamma}}\oint_{\{z\in M:\;d_g(a_{j},z)\geq \delta_{0}\}}\frac{1}{d_g(a_{j},z)^{2+2\gamma}}\;dS_g\\
&\le C\frac{\lambda_j^{2-2\gamma}}{\lambda_i^{1-\gamma}}\left(\frac{1}{1+\lambda_j^2d_g(a_{i},a_{j})^{2}}\right)^{1-\gamma}\frac{1}{\lambda_j^{1+\gamma}}\left[\int_{B_{\delta_{0}}\setminus B_{\delta}} \frac{1}{|z|^{2+2\gamma}}\;dz+C\right]\\
&\le C\left(\frac{\lambda_j}{\lambda_i}\right)^{1-\gamma}\left(\frac{1}{\lambda_j \delta}\right)^{2\gamma} \frac{1}{\left(1+\lambda_j^2d_g(a_{i},a_{j})^{2}\right)^{1-\gamma}}.
\end{split}
\end{equation*}
Thus, we get for \;$I^{2}_2$\;
\begin{equation}\label{I22B_j}
I^{2}_2=O\left(\left(\frac{1}{\delta\lambda_j}\right)^{2\gamma}\left(\frac{\lambda_{i}}{\lambda_{j}}+\lambda_i\lambda_j d_g(a_{i},a_{j})^{2}\right)^{-(1-\gamma)}\right).
\end{equation}\\
Hence, combining \eqref{I12B_j} and \eqref{I22B_j}, we have
\begin{equation}\label{I2B_j}
I_2=O\left(\left(\frac{1}{\delta^{4\gamma}\lambda_j^{2\gamma}}+\frac{1}{\delta^{2\gamma}\lambda^{2\gamma}_{j}}\right)\left(\frac{\lambda_{i}}{\lambda_{j}}+\lambda_i\lambda_j d_g(a_{i},a_{j})^{2}\right)^{-(1-\gamma)}\right).
\end{equation}\\
Therefore, using \eqref{A_j1}, \eqref{B_j1}, \eqref{IJ1A_j}, and \eqref{I2B_j}, we obtain\\
\begin{equation*}
    \begin{split}
     & d^{*}_{\gamma}\int_{X}\left|\left(-div_{g}(y^{1-2\gamma}\nabla_g)+q\right)u_{a_j,\lambda_j}\right|u_{a_i,\lambda_i}\;dV_g+\oint_{M}\left|-d^{*}_{\gamma}\lim_{y\to 0} y^{1-2\gamma}\;\frac{\partial u_{a_j,\lambda_j}}{\partial y}-u^{\frac{1+\gamma}{1-\gamma}}_{a_j,\lambda_j}\right|u_{a_i,\lambda_i}\;dS_g\\
     &\le C\left[\delta+\delta^{2-2\gamma}+\frac{1}{\delta^{2}\lambda^{2}_j}+\frac{1}{\delta^{4\gamma}\lambda_j^{2\gamma}}\right]\left(\frac{\lambda_{i}}{\lambda_{j}}+\lambda_i\lambda_j d_g(a_{i},a_{j})^{2}\right)^{-(1-\gamma)},
    \end{split}
\end{equation*}
Hence, the proof of Lemma \ref{interact1} is complete.
\end{pf}\\
\vspace{2pt}

\noindent
 Clearly Lemma \ref{interact1} implies the following sharp interaction-estimate relating \;$e_{ij},$\;\;$\epsilon_{ij},$\; and \;$\varepsilon_{ij}$\; (for their
definitions, see \eqref{varepij}-\eqref{eij}).
\begin{cor}\label{interact2}
Assuming that \;$\theta>0$\; is small and \;$\mu_0>0$\; is small, then \;$\forall$\;\;$a_i, a_j\in M$,\;\;$\forall 0<2\delta<\delta_0,$\; and \;$\forall$\;\;$0<\frac{1}{\l_j},$\;\;$\frac{1}{\l_i}\leq \theta\delta$\; such that \;$\varepsilon_{ij}\leq \mu_0,$\; we have 

\[e_{ij}=\epsilon_{ij}+O\left(\delta+\delta^{2-2\gamma}+\frac{1}{\delta^{2}\lambda^{2}_i}+\frac{1}{\delta^{4\gamma}\lambda_i^{2\gamma}}\right)\varepsilon_{ij},\]
where \;$\delta_0$\; is as in \eqref{delta0}.
\end{cor}
\vspace{6pt}

\noindent
In the next lemma, we present a sharp inter-action estimate relating \;$\epsilon_{ji}$\; and \;$\varepsilon_{ij}$.
\begin{lem}\label{interact3}
Assuming that \;$\theta>0$\; is small and \;$\mu_0>0$\; is small, then \;$\forall$\;\;$ a_i, a_j\in M$,\;\;$\forall 0<2\delta<\delta_0,$\; and \;$\forall$\;\;$0<\frac{1}{\l_j}\le\frac{1}{\l_i}\leq \theta\delta$\; such that \;$\varepsilon_{ij}\leq \mu_0,$\; we have 
\[\epsilon_{ji}=c_0^{\frac{2}{1-\gamma}}c_1\varepsilon_{ij}\left[\left(1+O\left(\delta^{2-2\gamma}+\frac{1}{\delta^2\lambda^2_i}\right)\right)\left(1+o_{\varepsilon_{ij}}(1)+O\left(\varepsilon^{\frac{\gamma}{1-\gamma}}_{ij} \delta^{-2\gamma}\right)\right)+O\left(\varepsilon^{\frac{\gamma}{1-\gamma}}_{ij}\delta^{-4}\right)\right],\]
where \;$c_0$\; is as in \eqref{bubble-M3} and \;$c_1$\; is as in \eqref{c3} .
\end{lem}
\noindent
\begin{pf}
By using \eqref{uald}, we have
\[u_{a_i,\lambda_i}(x)=\chi^{a_i}_\delta(x)\hat{\delta}_{a_i,\lambda_i}(x)+(1-\chi^{a_i}_\delta(x))\frac{c_0}{\lambda^{1-\gamma}}G_{a_i}(x),\;\;\;x\in M,\]
with \;$G_{a_i}(x)=G(a_i,x).$\; On the other hand, by using the definition of \;$\hat{\delta}_{a, \l}$\; (see \eqref{deltahat}), we have 
\[\chi^{a_i}_\delta(x)\hat{\delta}_{a_i,\lambda_i}(x)=c_0\chi^{a_i}_\delta(x)\left[\frac{\lambda_i}{1+\lambda^2_i G^{\frac{-1}{1-\gamma}}_{a_i}(x)\frac{d_g(x,a_i)^{2}}{G^{\frac{-1}{1-\gamma}}_{a_i}(x)}}\right]^{1-\gamma},\;\;\;x\in M.\]\\
\vspace{6pt}

\noindent
For \;$x\in\{y\in M:\;d_g(a_i,y)\le 2\delta\},$\; the quantity \;$\left[1+\lambda^2_i G^{\frac{-1}{1-\gamma}}_{a_i}(x)\frac{d_g(x,a_i)^{2}}{G^{\frac{-1}{1-\gamma}}_{a_i}(x)}\right]$\; can be estimated by \\
\begin{equation*}
\begin{split}
1+\lambda^2_i G^{\frac{-1}{1-\gamma}}_{a_i}(x)\frac{d_g(x,a_i)^{2}}{G^{\frac{-1}{1-\gamma}}_{a_i}(x)}&=1+\lambda^2_i G^{\frac{-1}{1-\gamma}}_{a_i}(x)\left(1+O\left(\delta^{2-2\gamma}\right)\right)\\
&=1+\lambda^2_i G^{\frac{-1}{1-\gamma}}_{a_i}(x)+O\left(\lambda^2_i \delta^{2-2\gamma} G^{\frac{-1}{1-\gamma}}_{a_i}(x)\right)\\
&=\left(1+\lambda^2_i G^{\frac{-1}{1-\gamma}}_{a_i}(x)\right)\left[1+O\left(\frac{\lambda^2_i \delta^{2-2\gamma} G^{\frac{-1}{1-\gamma}}_{a_i}(x)}{1+\lambda^2_i G^{\frac{-1}{1-\gamma}}_{a_i}(x)}\right)\right]\\
&=\left(1+\lambda^2_i G^{\frac{-1}{1-\gamma}}_{a_i}(x)\right)\left[1+O\left(\delta^{2-2\gamma}\right)\right].
\end{split}
\end{equation*}\\
Hence, we have
\begin{equation}\label{part1}
\begin{split}
\chi^{a_i}_\delta(x) \delta_{a_i},_{\lambda_i}(x)&=c_0\chi^{a_i}_\delta(x)\left[\frac{\lambda_i}{\left(1+\lambda^2_i G^{\frac{-1}{1-\gamma}}_{a_i}(x)\right)\left[1+O\left(\delta^{2-2\gamma}\right)\right]}\right]^{1-\gamma}\\
&=c_0 \chi^{a_i}_\delta(x)\left[1+O\left(\delta^{2-2\gamma}\right)\right]\left[\frac{\lambda_i}{1+\lambda^2_i G^{\frac{-1}{1-\gamma}}_{a_i}(x)}\right]^{1-\gamma},\;\;\;x\in M.
\end{split}
\end{equation}
Furthermore, we have
\[c_0(1-\chi^{a_i}_\delta(x))\left[\frac{\lambda_i}{1+\lambda^2_i G^{\frac{-1}{1-\gamma}}_{a_i}(x)}\right]^{1-\gamma}=(1-\chi^{a_i}_\delta(x))\frac{c_0}{\lambda_i^{1-\gamma}}G_{a_i}(x)\left[\frac{1}{1+\lambda^{-2}_i G^{\frac{1}{1-\gamma}}_{a_i}(x)}\right]^{1-\gamma},\;\;x\in M.\]
Since on $\{x\in M:\;d_g(x,a_i)\geq \d\}$, we have
\[ \frac{1}{1+\lambda^{-2}_iG^{\frac{1}{1-\gamma}}_{a_i}(x)}=1+O\left(\frac{G^{\frac{-1}{1-\gamma}}_{a_i}(x)}{\lambda^2_i}\right)=1+O\left(\frac{1}{\lambda^2_i \delta^2}\right),\]
then we get
\[c_0(1-\chi^{a_i}_\delta(x))\left[\frac{\lambda_i}{1+\lambda^2_i G^{\frac{-1}{1-\gamma}}_{a_i}(x)}\right]^{1-\gamma}=(1-\chi^{a_i}_\delta(x))\frac{c_0}{\lambda_i^{1-\gamma}}G_{a_i}(x)\left(1+O\left(\frac{1}{\lambda^2_i \delta^2}\right)\right),\;\;\;x\in M.\]
This implies 
\begin{equation}\label{part2}
(1-\chi^{a_i}_\delta(x))\frac{c_0}{\lambda^{1-\gamma}}G_{a_i}(x)=c_0(1-\chi^{a_i}_\delta(x))\left[\frac{\lambda}{1+\lambda^2_i G^{\frac{-1}{1-\gamma}}_{a_i}(x)}\right]^{1-\gamma}\left(1+O\left(\frac{1}{\lambda^2_i \delta^2}\right)\right),\;\;\;x\in M.
\end{equation}
Thus, combining \eqref{part1} and \eqref{part2}, we get
\begin{equation*}
    \begin{split}
       u_{a_i},_{\lambda_i}(x)=c_0\left[\left(1+O\left(\delta^{2-2\gamma}\right)\right)\chi^{a_i}_\delta+(1-\chi^{a_i}_\delta)\left(1+O\left(\frac{1}{\lambda^2_i \delta^2}\right)\right)\right]\left[\frac{\lambda}{1+\lambda^2_i G^{\frac{-1}{1-\gamma}}_{a_i}(x)}\right]^{1-\gamma},\;\;\;x\in M. 
    \end{split}
\end{equation*}
Hence, we obtain the following 
\begin{equation}\label{partf}
u_{a_i, \lambda_i}(x)=c_0\left[1+O\left(\delta^{2-2\gamma}\right)+O\left(\frac{1}{\lambda^2_i\delta^2}\right)\right]\left[\frac{\lambda}{1+\lambda^2_i G^{\frac{-1}{1-\gamma}}_{a_i}(x)}\right]^{1-\gamma},\;\;\;x\in M, 
\end{equation}
\vspace{4pt}
\noindent
Now, we are going to finish our work by using \eqref{partf}. To start, we write our estimate in the following  form:
\begin{equation}\label{J1+J2}
\epsilon_{ji}=\oint_{M} u^{\frac{1+\gamma}{1-\gamma}}_{a_j,\lambda_j}u_{a_i,\lambda_i}\;dS_g=\underbrace{\oint_{B^{\hat{h}}_{\delta}(a_j)} u^{\frac{1+\gamma}{1-\gamma}}_{a_j,\lambda_j}u_{a_i,\lambda_i}\;dS_g}_{\mbox{$J_1$}}+\underbrace{\oint_{M-B^{\hat{h}}_{\delta}(a_j)} u^{\frac{1+\gamma}{1-\gamma}}_{a_j,\lambda_j}u_{a_i,\lambda_i}\;dS_g}_{\mbox{$J_2$}}.
\end{equation}
In order to estimate \;$J_1$\; and \;$J_2,$\; we are going to estimate them individually. To begin with \;$J_2,$\; we have
\begin{equation*}
\begin{split}
\oint_{ M-B^{\hat{h}}_{\delta}(a_j)} u^{\frac{1+\gamma}{1-\gamma}}_{a_j,\lambda_j}u_{a_i,\lambda_i}\;dS_g&\le C\oint_{ M-B^{\hat{h}}_{\delta}(a_j)}\left(\frac{1}{\lambda_j}\right)^{1+\gamma}\left(\frac{1}{\delta}\right)^{2(1+\gamma)}u_{a_i,\lambda_i}\;dS_g\\
&\le C \left(\frac{1}{\lambda_j}\right)^{1+\gamma}\left(\frac{1}{\delta}\right)^{2(1+\gamma)}\oint_{M-B^{\hat{h}}_{\delta}(a_j)}u_{a_i,\lambda_i}\;dS_g\\
&\le C\left(\frac{1}{\lambda_j\d^2}\right)^{1+\gamma}\oint_{M-(B^{\hat{h}}_{\delta}(a_j)\cup B^{\hat{h}}_{\delta}(a_i))}u_{a_i,\lambda_i}\;dS_g\\
&+C\left(\frac{1}{\lambda_j\d^2}\right)^{1+\gamma}\oint_{(M-B^{\hat{h}}_{\delta}(a_j))\cap B^{\hat{h}}_{\delta}(a_i)}u_{a_i,\lambda_i}\;dS_g\\
&\le C\left(\frac{1}{\lambda_j}\right)^{1+\gamma}\left(\frac{1}{\delta}\right)^4\frac{1}{\lambda_i^{1-\gamma}}\\
&+C\left(\frac{1}{\lambda_j}\right)^{1+\gamma}\left(\frac{1}{\delta}\right)^{2(1+\gamma)}\oint_{B^{\hat{h}}_{\delta}(a_i)}\left(\frac{\lambda_i}{1+\lambda^2_id_g(a_i,x)^{2}}\right)^{1-\gamma}\;dS_g\\
&\le C\left(\frac{1}{\lambda_j}\right)^{1+\gamma}\left(\frac{1}{\delta}\right)^4\frac{1}{\lambda_i^{1-\gamma}}+C\left(\frac{1}{\lambda_j}\right)^{1+\gamma}\left(\frac{1}{\delta}\right)^{2(1+\gamma)}\frac{\delta^{2\gamma}}{\lambda_i^{1-\gamma}}\\
&\le C\left(\frac{1}{\lambda_j}\right)^{1+\gamma}\left(\frac{1}{\delta}\right)^4\frac{1}{\lambda_i^{1-\gamma}}\left(1+\delta^2\right)\\
&\le \frac{C}{\lambda^{1+\gamma}_j\lambda_i^{1-\gamma}\delta^4}.
\end{split}
\end{equation*}

\noindent
Thus, we  get for \;$J_2$\;
\begin{equation}\label{J_22}
J_2=O\left(\varepsilon^{\frac{1}{1-\gamma}}_{ij}\frac{1}{\delta^4}\right).
\end{equation}

\noindent
In the next step, using \eqref{partf} for \;$J_1,$\; we have 
\begin{equation*}
\begin{split}
&\oint_{B^{\hat{h}}_{\delta}(a_j)} u^{\frac{1+\gamma}{1-\gamma}}_{a_j,\lambda_j}u_{a_i,\lambda_i}\;dS_g\\
&=c^{\frac{2}{1-\gamma}}_0\oint_{B^{\hat{h}}_{\delta}(a_j)}\left(\frac{\lambda_j}{1+\lambda^2_jd_g(a_j,x)^{2}}\right)^{1-\gamma}\left[1+O\left(\delta^{2-2\gamma}\right)+O\left(\frac{1}{\lambda^2_i\delta^2}\right)\right]\left[\frac{\lambda_i}{1+\lambda^2_i G^{\frac{-1}{1-\gamma}}_{a_i}(x)}\right]^{1-\gamma}\;dS_g\\
&=c^{\frac{2}{1-\gamma}}_0\left[1+O\left(\delta^{2-2\gamma}\right)+O\left(\frac{1}{\lambda^2_i\delta^2}\right)\right]\oint_{B^{\hat{h}}_{\delta}(a_j)}\left(\frac{\lambda_j}{1+\lambda^2_jd_g(a_j,x)^{2}}\right)^{1+\gamma}\left[\frac{\lambda_i}{1+\lambda^2_i G^{\frac{-1}{1-\gamma}}_{a_i}(x)}\right]^{1-\gamma}\;dS_g\\
&=c^{\frac{2}{1-\gamma}}_0 \left[1+O\left(\delta^{2-2\gamma}\right)+O\left(\frac{1}{\lambda^2_i\delta^2}\right)\right]\frac{1}{\lambda_j^{1-\gamma}}\int_{B_{\lambda_j\delta}}\left(\frac{1}{1+|y|^2}\right)^{1+\gamma}\left[\frac{\lambda_i}{1+\lambda^2_i G^{\frac{-1}{1-\gamma}}_{a_i}\left(\Psi_{a_j}\left(\frac{y}{\lambda_j}\right)\right)}\right]^{1-\gamma}\\
&=c^{\frac{2}{1-\gamma}}_0 \left[1+O\left(\delta^{2-2\gamma}\right)+O\left(\frac{1}{\lambda^2_i\delta^2}\right)\right]\int_{B_{\lambda_j\delta}}\left(\frac{1}{1+|y|^2}\right)^{1+\gamma}\left[\frac{1}{\frac{\l_j}{\l_i}+\lambda_i \l_jG^{\frac{-1}{1-\gamma}}_{a_i}\left(\Psi_{a_j}\left(\frac{y}{\lambda_j}\right)\right)}\right]^{1-\gamma}.
\end{split}
\end{equation*}
\vspace{6pt}

\noindent
Recalling that $\lambda_i\le\lambda_j$, then  for \ $\varepsilon_{ij}\sim 0 $, we have\\
1) Either $\varepsilon^{\frac{-1}{1-\gamma}}_{ij}\sim \lambda_i\lambda_jG^{\frac{-1}{1-\gamma}}_{a_i}(a_j).$
\vspace{4pt}

\noindent
2) or $\varepsilon^{\frac{-1}{1-\gamma}}_{ij}\sim\frac{\lambda_j}{\lambda_i}$.
\vspace{6pt}

\noindent
In order to estimate \;$J_1,$\; we first define the following sets
\begin{equation*}
    \begin{split}
        A_1&=\left\{y\in \R^{2}:\;\left(\left|\frac{y}{\lambda_j}\right|\le \epsilon \:G^{\frac{-1}{2-2\gamma}}_{a_i}(a_j)\right)\cap B_{\lambda_j\delta}\right\},\\
        A_2&=\left\{y\in \R^{2}:\;\left(\left|\frac{y}{\lambda_j}\right|\le\epsilon\frac{1}{\lambda_i}\right)\cap B_{\lambda_j\delta}\right\},
    \end{split}
\end{equation*}
and
\[\mathcal{A}=A_1\cup A_2\;,\]
with \;$\epsilon>0$\; very small. Then  by Taylor expansion on \;$\mathcal{A}$, we have
\begin{equation*}
\begin{split}
&\left[\frac{\lambda_j}{\lambda_i}+\lambda_i\lambda_jG^{\frac{-1}{1-\gamma}}_{a_i}\left(\psi_{a_j}\left(\frac{y}{\lambda_j}\right)\right)\right]^{-(1-\gamma)}\\
&=\left[\frac{\lambda_j}{\lambda_i}+\lambda_i\lambda_jG^{\frac{-1}{1-\gamma}}_{a_i}(a_j)\right]^{-(1-\gamma)}\\
&+\left[\left(-\frac{1}{2}\nabla  G^{\frac{-1}{1-\gamma}}_{a_i}\circ \psi_{a_j}(a_j)\lambda_i y\right)\right]\left[\frac{\lambda_j}{\lambda_i}+\lambda_i\lambda_jG^{\frac{-1}{1-\gamma}}_{a_i}\left(\psi_{a_j}\left(\frac{y}{\lambda_j}\right)\right)\right]^{-(2-\gamma)}\\
&+O\left[\left(\frac{\lambda_i}{\lambda_j}\right)|y|^2\right]\left[\frac{\lambda_j}{\lambda_i}+\lambda_i\lambda_jG^{\frac{-1}{1-\gamma}}_{a_i}\left(\Psi_{a_j}\left(\frac{y}{\lambda_j}\right)\right)\right]^{-(2-\gamma)}.
\end{split}
\end{equation*}

\noindent
Thus, we write \;$J_1$\; the following  
\begin{equation}\label{J_11}
J_1=c_0^{\frac{2}{1-\gamma}}\left[1+O\left(\delta^{2-2\gamma}\right)+O\left(\frac{1}{\lambda^2_i\delta^2}\right)\right]\left(\sum_{m=1}^{4} I_m\right),
\end{equation}
with
\begin{equation*}
\begin{split}
I_1&=\left[\frac{\lambda_j}{\lambda_i}+\lambda_i\lambda_jG^{\frac{-1}{1-\gamma}}_{a_i}(a_j)\right]^{-(1-\gamma)}\int_{\mathcal{A}}\left(\frac{1}{1+|y|^2}\right)^{1+\gamma},\\
I_2&=\left[\frac{\lambda_j}{\lambda_i}+\lambda_i\lambda_jG^{\frac{-1}{1-\gamma}}_{a_i}(a_j)\right]^{-(2-\gamma)}\int_{\mathcal{A}}\left(\frac{1}{1+|y|^2}\right)^{1+\gamma}\left[\nabla G^{\frac{-1}{1-\gamma}}_{a_i}\circ \Psi_{a_j}(a_j)\lambda_i y\right],\\
I_3&=\left[\frac{\lambda_j}{\lambda_i}+\lambda_i\lambda_jG^{\frac{-1}{1-\gamma}}_{a_i}(a_j)\right]^{-(2-\gamma)}\int_{\mathcal{A}}\left(\frac{1}{1+|y|^2}\right)^{1+\gamma}O\left[\left(\frac{\lambda_i}{\lambda_j}\right)|y|^2\right],
\end{split}
\end{equation*}
and
\begin{equation*}
\begin{split}
I_4&=\int_{B_{\lambda_j\delta}-\mathcal{A}}\left(\frac{1}{1+|y|^2}\right)^{1+\gamma}\left[\frac{\lambda_j}{\lambda_i}+\lambda_i\lambda_jG^{\frac{-1}{1-\gamma}}_{a_i}\left(\Psi_{a_j}\left(\frac{y}{\lambda_j}\right)\right)\right]^{-(1-\gamma)}.
\end{split}
\end{equation*}

\noindent
Now, let us estimate \;$I_1.$\; We have 
\[I_1=\left[\frac{\lambda_j}{\lambda_i}+\lambda_i\lambda_jG^{\frac{-1}{1-\gamma}}_{a_i}(a_j)\right]^{-(1-\gamma)}\left[c_1+\int_{\R^2-\mathcal{A}}\left(\frac{1}{1+|y|^2}\right)^{1+\gamma}\right],\]
where \;$c_1$\; is as in \eqref{c3}. On the other hand, we set
\begin{equation*}
    \begin{split}
     T_{ij}&=\lambda_j\epsilon G^{\frac{-1}{2-2\gamma}}_{a_i}(a_i,a_j),\\
     L_{ij}&=\epsilon\frac{\lambda_j}{\lambda_i},
    \end{split}
\end{equation*}
and have
\[\int_{\R^2-\mathcal{A}}\left(\frac{1}{1+|y|^2}\right)^{1+\gamma}\le \int_{\R^2-B_{\delta\lambda_j}}\left(\frac{1}{1+|y|^2}\right)^{1+\gamma}+\int_{\R^2-B_{T_{ij}}}\left(\frac{1}{1+|y|^2}\right)^{1+\gamma}\]\\
\noindent
if \;$\varepsilon^{\frac{-1}{1-\gamma}}_{ij}\sim \lambda_i\lambda_jG^{\frac{-1}{1-\gamma}}_{a_i}(a_j),$\; and 
\begin{equation*}
\int_{\R^2-\mathcal{A}}\left(\frac{1}{1+|y|^2}\right)^{1+\gamma}\le \int_{\R^2-B_{\delta\lambda_j}}\left(\frac{1}{1+|y|^2}\right)^{1+\gamma}+\int_{\R^2-B_{L_{ij}}}\left(\frac{1}{1+|y|^2}\right)^{1+\gamma}
\end{equation*}
\noindent
if \;$\varepsilon^{\frac{-1}{1-\gamma}}_{ij}\sim\frac{\lambda_j}{\lambda_i}.$\; We have 
\begin{equation*}
\int_{\R^2-B_{\delta\lambda_j}}\left(\frac{1}{1+|y|^2}\right)^{1+\gamma}=O\left(\frac{1}{\lambda_j^{2\gamma}\delta^{2\gamma}}\right).
\end{equation*}
\noindent
Moreover, if \;$\varepsilon^{\frac{-1}{1-\gamma}}_{ij}\sim \lambda_i\lambda_jG^{\frac{-1}{1-\gamma}}_{a_i}(a_j),$\; then
\begin{equation*}
    \begin{split}
\int_{\R^2-B_{T_{ij}}}\left(\frac{1}{1+|y|^2}\right)^{1+\gamma}&=O\left(\frac{1}{\lambda_j^{2\gamma} G^{\frac{-\gamma}{1-\gamma}}_{a_i}(a_j)}\right)\\
&=O\left(\frac{1}{\lambda_j^{\gamma}\lambda_i^{\gamma} G^{\frac{-\gamma}{1-\gamma}}_{a_i}(a_j)}\right)\\
&=O\left(\varepsilon_{ij}^{\frac{\gamma}{1-\gamma}}\right).
\end{split}
\end{equation*}

\noindent
Furthermore if \;$\varepsilon^{\frac{-1}{1-\gamma}}_{ij}\sim\frac{\lambda_j}{\lambda_i},$ then
\begin{equation*}
\int_{\R^2-B_{L_{ij}}}\left(\frac{1}{1+|y|^2}\right)^{1+\gamma}=O\left(\varepsilon_{ij}^{\frac{\gamma}{1-\gamma}}\right).
\end{equation*}
This implies 
\begin{equation*}
\int_{\R^2-\mathcal{A}}\left(\frac{1}{1+|y|^2}\right)^{1+\gamma}=O\left(\varepsilon_{ij}^{\frac{\gamma}{1-\gamma}}+\frac{1}{\l_j^{2\gamma}\delta^{2\gamma}}\right)=O\left(\varepsilon_{ij}^{\frac{\gamma}{1-\gamma}}+\varepsilon_{ij}^{\frac{\gamma}{1-\gamma}}\delta^{-2\gamma}\right)=O\left(\varepsilon_{ij}^{\frac{\gamma}{1-\gamma}}\delta^{-2\gamma}\right).
\end{equation*}

\noindent
Thus, we get\
\begin{equation*}
    \begin{split}
I_1&=\left[\frac{\lambda_j}{\lambda_i}+\lambda_i\lambda_jG^{\frac{-1}{1-\gamma}}_{a_i}(a_j)\right]^{-(1-\gamma)}\left[c_1+O\left(\varepsilon_{ij}^{\frac{\gamma}{1-\gamma}}\delta^{-2\gamma}\right)\right]\\
&=\varepsilon_{ij}\left(1+o_{\varepsilon_{ij}}(1)\right)\left[c_1+O\left(\varepsilon_{ij}^{\frac{\gamma}{1-\gamma}}\delta^{-2\gamma}\right)\right]
\end{split}
\end{equation*}
\noindent
Hence, we obtain
\begin{equation}\label{I_11}
I_1=c_1\varepsilon_{ij}\left[1+o_{\epsilon_{ij}}(1)+O\left(\varepsilon_{ij}^{\frac{\gamma}{1-\gamma}}\delta^{-2\gamma}\right)\right].
\end{equation}
By symmetry, we have
 \begin{equation}\label{I_22}
 I_2=0.
 \end{equation}
 
 \noindent
 Next, for \;$I_3,$\; we derive
\begin{equation*}
\begin{split}
\int_{\mathcal{A}}\frac{|y|^2}{\left(1+|y^2|\right)^{1+\gamma}}
&\le \int_{B_{T_{ij}}}\frac{|y|^2}{\left(1+|y|^{2}\right)^{1+\gamma}}+\int_{B_{L_{ij}}}\frac{|y|^2}{\left(1+|y|^2\right)^{1+\gamma}}\\
&=O\left(\left(\epsilon\lambda_jG^{\frac{-1}{2-2\gamma}}_{a_i}(a_j)\right)^{2(1-\gamma)}+\left(\epsilon \frac{\lambda_j}{\lambda_i}\right)^{2(1-\gamma)}\right).
\end{split}
\end{equation*}
Thus, we have 
\begin{equation*}
\begin{split}
I_3&=\varepsilon^{\frac{2-\gamma}{1-\gamma}}_{ij}\left(\frac{\lambda_i}{\lambda_j}\right)\left(1+o_{\varepsilon_{ij}}(1)\right)\left[O\left(\left(\epsilon\lambda_jG^{\frac{-1}{2-2\gamma}}_{a_i}(a_j)\right)^{2(1-\gamma)}+\left(\epsilon \frac{\lambda_j}{\lambda_i}\right)^{2(1-\gamma)}\right)\right]\\
&=\varepsilon^{\frac{2-\gamma}{1-\gamma}}_{ij}\left(1+o_{\varepsilon_{ij}}(1)\right)\left[O\left(\left(\sqrt{\lambda_i\lambda}_jG^{\frac{-1}{2-2\gamma}}_{a_i}(a_j)\right)^{2(1-\gamma)}+\left(\sqrt{\frac{\lambda_j}{\lambda_i}}\right)^{2(1-\gamma)}\right)\right].
\end{split}
\end{equation*}

\noindent
Hence, we obtain
\begin{equation}\label{I_33}
I_3=O\left(\varepsilon^{\frac{1}{1-\gamma}}_{ij}\right).
\end{equation}

\noindent
Finally, we estimate \;$I_4$\; as follows.\\\\
If \;$\varepsilon^{\frac{-1}{1-\gamma}}_{ij}\sim\frac{\lambda_j}{\lambda_i},$\; then\\
\begin{equation}\label{part1i4}
I_4\le C\varepsilon_{ij}\int_{B_{\lambda_j\delta}-\mathcal{A}}\left(\frac{1}{1+|y|^2}\right)^{1+\gamma}\le C\varepsilon_{ij}\left(\frac{\lambda_j}{\lambda_i}\right)^{-2\gamma}\le C \varepsilon^{\frac{1+\gamma}{1-\gamma}}_{ij}.
\end{equation}\\
If \;$\varepsilon^{\frac{-1}{1-\gamma}}_{ij}\sim \lambda_i\lambda_jG^{\frac{-1}{1-\gamma}}_{a_i}(a_j),$\; then we argue as follows.
In case \;$d_g(a_i,a_j)\geq2\delta,$\; since 
\[G_{a_i}\left(\Psi_{a_j}\left(\frac{y}{\lambda_j}\right)\right)\le C\delta^{-2(1-\gamma)}\]
 for \;$y\in B_{\l_j\d},$\; then we have 
\begin{equation*}
\begin{split}
I_4&\le C\int_{B_{\lambda_j\delta}-\mathcal{A}}\left(\frac{1}{1+|y|^2}\right)^{1+\gamma} \frac{1}{\left(\lambda_i\lambda_j\delta^{2}\right)^{1-\gamma}}\\
&\le\frac{C}{\left(\lambda_i\lambda_j\delta^{2}\right)^{1-\gamma}}\left(\frac{1}{\lambda_j^{2\gamma}G_{a_i}^{\frac{\gamma}{1-\gamma}}(a_j)}\right)\\&\le \frac{C}{\left(\lambda_i\lambda_j\right)^{\gamma}G^{\frac{\gamma}{1-\gamma}}_{a_i}(a_j)}\frac{G^{-1}_{a_i}(a_j)}{\left(\lambda_i\lambda_j\delta^{2}\right)^{1-\gamma}G_{a_i}^{-1}(a_j)}\\&\le C\varepsilon^{\frac{\gamma}{1-\gamma}}_{ij} \varepsilon_{ij} \delta^{-2(1-\gamma)}.
\end{split}
\end{equation*}

\noindent
Thus, when \; $d_g(a_i,a_j)\geq 2\d,$\; we have 
\begin{equation}\label{part2i4}
I_4=O\left(\varepsilon^{\frac{1}{1-\gamma}}_{ij}\delta^{-2(1-\gamma)}\right).
\end{equation}
In case \;$d_g(a_i,a_j)<2\delta,$\; we first observe that 
$$
B_{\l_j\d}\setminus \mathcal{A}\subset A_1\cup A_2
$$
with 
$$
A_1=\left\{y\in \R^{2}:\;\epsilon\lambda_jG^{\frac{-\gamma}{1-\gamma}}_{a_i}(a_j)\le|y|\le E\lambda_jG^{\frac{-\gamma}{1-\gamma}}_{a_i}(a_j)\right\}
$$
and
$$
A_2=\left\{y\in \R^{2}:\;E\lambda_jG^{\frac{-\gamma}{1-\gamma}}_{a_i}(a_j)\le|y|\le \lambda_j\delta\right\},
$$
where \;\;$0<\epsilon<E.$\;\\\\
\noindent
Thus, we have 
\begin{equation}\label{i4ine}
I_4\le I^1_4+I^2_4,
\end{equation}
with
\begin{equation*}
I_4^1=\int_{A_1} \left(\frac{1}{1+|y|^2}\right)^{1+\gamma}\left[\frac{\lambda_j}{\lambda_i}+\lambda_i\lambda_jG^{\frac{-1}{1-\gamma}}_{a_i}\left(\Psi_{a_j}\left(\frac{y}{\lambda_j}\right)\right)\right]^{-(1-\gamma)},
\end{equation*}
and
\begin{equation*}
I_4^2=\int_{A_2} \left(\frac{1}{1+|y|^2}\right)^{1+\gamma}\left[\frac{\lambda_j}{\lambda_i}+\lambda_i\lambda_jG^{\frac{-1}{1-\gamma}}_{a_i}\left(\Psi_{a_j}\left(\frac{y}{\lambda_j}\right)\right)\right]^{-(1-\gamma)}.
\end{equation*}
For now, we first estimate \;$I_4^1$\; as follows:
\begin{equation*}
\begin{split}
I_4^1&\le C\left[1+\lambda_j^2G^{\frac{-1}{1-\gamma}}_{a_i}(a_j)\right]^{-(1+\gamma)}\int_{\{y\in \R^{2}:\;|y|\le E\lambda_jG^{\frac{-1}{2-2\gamma}}_{a_i}(a_j)\}}\left[\frac{\lambda_j}{\lambda_i}+\lambda_i\lambda_jG^{\frac{-1}{1-\gamma}}_{a_i}\left(\Psi_{a_j}\left(\frac{y}{\lambda_j}\right)\right)\right]^{-(1-\gamma)}\\
&\le C\left[1+\lambda_j^2G^{\frac{-1}{1-\gamma}}_{a_i}(a_j)\right]^{-(1+\gamma)}\left(\frac{\lambda_i}{\lambda_j}\right)^{1-\gamma}\int_{\{y\in \R^{2}:\;|y|\le E\lambda_jG^{\frac{-1}{2-2\gamma}}_{a_i}(a_j)\}}\left[1+\lambda^2_iG^{\frac{-1}{1-\gamma}}_{a_i}\left(\Psi_{a_j}\left(\frac{y}{\lambda_j}\right)\right)\right]^{-(1-\gamma)}\\
&\le C\left[1+\lambda_j^2G^{\frac{-1}{1-\gamma}}_{a_i}(a_j)\right]^{-(1+\gamma)}\left(\frac{\lambda_i}{\lambda_j}\right)^{1-\gamma}\int_{\{y\in \R^{2}:\;|y|\le E\lambda_jG^{\frac{-1}{2-2\gamma}}_{a_i}(a_j)\}}\left[1+\lambda^2_i\left|\Psi^{-1}_{a_i}\circ\Psi_{a_j}\left(\frac{y}{\lambda_j}\right)\right|^2\right]^{-(1-\gamma)}\\
&\le C\left[1+\lambda_j^2G^{\frac{-1}{1-\gamma}}_{a_i}(a_j)\right]^{-(1+\gamma)}
\left(\frac{\lambda_i}{\lambda_j}\right)^{1-\gamma}\left(\frac{\lambda_j}{\lambda_i}\right)^{2}\int_{\{z\in \R^{2}:\;|z|\le \bar E\lambda_iG^{\frac{-1}{2-2\gamma}}_{a_i}(a_j)\}}\left[\frac{1}{1+|z|^2}\right]^{1-\gamma}\\
&\le C \left[\frac{\lambda_i}{\lambda_j}+\lambda_i\lambda_jG^{\frac{-1}{1-\gamma}}_{a_i}(a_j)\right]^{-(1+\gamma)}\left(\lambda_iG^{\frac{-1}{2-2\gamma}}_{a_i}(a_j)\right)^{2\gamma}
\\&
\le C \varepsilon^{\frac{1+\gamma}{1-\gamma}}_{ij}\left(\sqrt{\lambda_i\lambda_j} G^{\frac{-1}{2-2\gamma}}_{a_i}(a_j)\right)^{2\gamma},
\end{split}
\end{equation*}

\noindent
where \;$\bar E$\; is a positive constant. So we obtain
\begin{equation}\label{i14fin}
I^1_4=O\left(\varepsilon^{\frac{1}{1-\gamma}}_{ij}\right).
\end{equation}
\noindent
For \;$I_4^2,$\; we have
\begin{equation*}
\begin{split}
I^2_4&=\int_{A_2}\left(\frac{1}{1+|y|^2}\right)^{1+\gamma}\left[\frac{\lambda_j}{\lambda_i}+\lambda_i\lambda_jG^{\frac{-1}{1-\gamma}}_{a_i}\left(\Psi_{a_j}\left(\frac{y}{\lambda_j}\right)\right)\right]^{-(1-\gamma)}\\
&\le C\int_{\{y\in \R^{2}:\;|y|\geq E\lambda_jG^{\frac{-1}{2-2\gamma}}_{a_i}(a_j)\}}\left(\frac{1}{1+|y|^2}\right)^{1+\gamma}\left[\frac{\lambda_j}{\lambda_i}+\lambda_i\lambda_jG^{\frac{-1}{1-\gamma}}_{a_i} (a_j )\right]^{-(1-\gamma)}\\
&\le C\left[\frac{\lambda_j}{\lambda_i}+\lambda_i\lambda_jG^{\frac{-1}{1-\gamma}}_{a_i} (a_j )\right]^{-(1-\gamma)}\left(\frac{1}{\lambda_jG^{\frac{-1}{2-2\gamma}}_{a_i} (a_j )}\right)^{2\gamma}.
\end{split}
\end{equation*}

\noindent
This implies
 \begin{equation}\label{i24fin}
 I^2_4=O\left(\varepsilon^{\frac{1}{1-\gamma}}_{ij}\right).
 \end{equation}
Thus, combining \eqref{i4ine} and \eqref{i24fin}, we have that if \;$d_g(a_i,a_j)<2\d,$\; then 
\begin{equation}\label{i4pcclose}
I_4=O\left(\varepsilon^{\frac{1}{1-\gamma}}_{ij}\right).
\end{equation}

\noindent
Now, using \eqref{part2i4} and \eqref{i4pcclose}, we infer that  in case \;$\varepsilon_{i, j}^{\frac{-1}{1-\gamma}}\simeq \l_i\l_jG_{a_i}^{\frac{-1}{1-\gamma}}(a_j),$\; 
\begin{equation}\label{i4pc}
I_4=O\left(\varepsilon^{\frac{1}{1-\gamma}}_{i, j}\delta^{-2(1-\gamma)}\right).
\end{equation}
\noindent
Finally combining \eqref{part1i4}-\eqref{i4pc}, we get
\begin{equation}\label{I_44}
    I_4=O\left(\varepsilon^{\frac{1}{1-\gamma}}_{i, j}\delta^{-2(1-\gamma)}\right).
\end{equation} 

\noindent
Using \eqref{J_11}-\eqref{I_33}, and \eqref{I_44}, we obtain the following for \;$J_1$\; (see \eqref{J1+J2})  
\begin{equation}\label{J_111}
J_1=c_0^{\frac{2}{1-\gamma}}\left[1+O\left(\delta^{2-2\gamma}+\frac{1}{\lambda^2_i\delta^2}\right)\right]\left[c_1\varepsilon_{ij}\left(1+o_{\varepsilon_{ij}}(1)+O\left(\varepsilon^{\frac{\gamma}{1-\gamma}}_{ij} \delta^{-2\gamma}\right)\right)\right].
\end{equation}
\noindent
Thus, using \eqref{J1+J2}, \eqref{J_22}, and \eqref{J_111}, we arrive to 
\begin{equation*}
\begin{split}
\oint_{M} u^{\frac{1+\gamma}{1-\gamma}}_{a_j,\lambda_j}u_{a_i,\lambda_i}\;dS_g=&c_0^{\frac{2}{1-\gamma}}\left[1+O\left(\delta^{2-2\gamma}+\frac{1}{\lambda^2_i\delta^2}\right)\right]\left[c_1\varepsilon_{ij}\left(1+o_{\varepsilon_{ij}}(1)+O\left(\varepsilon^{\frac{\gamma}{1-\gamma}}_{ij} \delta^{-2\gamma}\right)\right)\right]\\+&O\left(\varepsilon^{\frac{1}{1-\gamma}}_{ij}\frac{1}{\delta^4}\right).
\end{split}
\end{equation*}

\noindent
Therefore, we obtain
\begin{equation}\label{eqfinal}
\begin{split}
\oint_{M} u^{\frac{1+\gamma}{1-\gamma}}_{a_j,\lambda_j}u_{a_i,\lambda_i}\;dS_g=&c_0^{\frac{2}{1-\gamma}}c_1\varepsilon_{ij}\left[\left(1+O\left(\delta^{2-2\gamma}+\frac{1}{\lambda^2_i\delta^2}\right)\right)\left(1+o_{\varepsilon_{ij}}(1)+O\left(\varepsilon^{\frac{\gamma}{1-\gamma}}_{ij} \delta^{-2\gamma}\right)\right)\right]\\+&O\left(\varepsilon_{ij}^{\frac{1}{1-\gamma}}\frac{1}{\delta^4}\right).
\end{split}
\end{equation}
Hence, recalling \;$\epsilon_{ji}=\oint_{M} u^{\frac{1+\gamma}{1-\gamma}}_{a_j,\lambda_j}u_{a_i,\lambda_i}\;dS_g,$\;  then the result follows from \eqref{eqfinal}
\end{pf}
\vspace{8pt}

\noindent
The following corollary is equivalent to Lemma \ref{interact3} when the indexes \;$i$\; and \;$j$\; are switched in Lemma 5.3. We decide to present the following corollary due to the fact that its form is more consistent with our presentation of the Barycenter Technique of Bahri-Coron\cite{bc}, which is based on the work \cite{martndia2}  as it was done in \cite{nss}.
\begin{cor}\label{interact4}
Assuming that \;$\theta>0$\; is small and \;$\mu_0>0$\; is small then \;$\forall$\;\;$ a_i, a_j\in M,$\;\;$\forall\; 0<2\delta<\delta_0,$\; and \;$\forall$\;\;$0<\frac{1}{\l_i}\le\frac{1}{\l_j}\leq \theta\delta$\; such that \;$\varepsilon_{ij}\leq \mu_0,$\; we have
\[\epsilon_{ij}=c_0^{\frac{2}{1-\gamma}}c_1\varepsilon_{ij}\left[\left(1+O\left(\delta^{2-2\gamma}+\frac{1}{\lambda^2_j\delta^2}\right)\right)\left(1+o_{\varepsilon_{ij}}(1)+O\left(\varepsilon^{\frac{\gamma}{1-\gamma}}_{ij} \delta^{-2\gamma}\right)\right)\right]\\+O\left(\varepsilon_{ij}^{\frac{1}{1-\gamma}}\delta^{-4}\right),\]
where \;$\delta_0$\; is as in \eqref{delta0}.
\end{cor}
\vspace{8pt}

\noindent
In the following lemma, we present some sharp high-order inter-action estimates that are needed  for the application of the algebraic topological argument of Bahri-Coron\cite{bc} for existence. To begin, we are going to derive the balanced high-order inter-action estimate shown below.
\begin{lem}\label{interact5}
Assuming that \;$\theta>0$\; is small and \;$\mu_0>0$\; is small then \;$\forall$\;\;$ a_i, a_j\in M$,\;\;$\forall 0<2\delta<\delta_0,$\; and \;$\forall$\;\;$0<\frac{1}{\l_i},$\;\;$\frac{1}{\l_j}\leq \theta\delta$\; such that \;$\varepsilon_{ij}\leq \mu_0,$\; we have
\[\oint_{M} u^{\frac{1}{1-\gamma}}_{a_i,\lambda_i} u^{\frac{1}{1-\gamma}}_{a_j,\lambda_j}\;dS_g=O\left(\varepsilon^{{\frac{1}{1-\gamma}}}_{ij} \delta^{-4}\log\left(\varepsilon^{\frac{-1}{2-2\gamma}}_{ij}\delta^{-1}\right)\right).\]
\end{lem}
\begin{pf}
By symmetry, we can assume without loss of generality (w.l.o.g) that\;$\lambda_j\le\lambda_i.$\; Thus we have \\
1) Either \;$\varepsilon^{{\frac{-1}{1-\gamma}}}_{ij}\sim \lambda_i\lambda_jG^{{\frac{-1}{1-\gamma}}}_{a_i}(a_j).$\;
\vspace{4pt}

\noindent
2) Or \;$\varepsilon^{{\frac{-1}{1-\gamma}}}_{ij}\sim\frac{\lambda_i}{\lambda_j}.$\;
\vspace{6pt}

\noindent
Now, if \;$d_g(a_i,a_j)\geq2\delta,$\; then we have 
\begin{equation}\label{i1}
\begin{split}
I&:=\oint_{M} u^{{\frac{1}{1-\gamma}}}_{a_i,\lambda_i} u^{{\frac{1}{1-\gamma}}}_{a_j,\lambda_j}\;dS_g\\
&\le C\oint_{B^{\hat{h}}_{\delta}(a_i)}\left(\frac{\lambda_i}{1+\lambda^2_id_g(a_i,x)^{2}}\right) \left(\frac{\lambda_j}{1+\lambda^2_jG^{{\frac{-1}{1-\gamma}}}_{a_j}(x)}\right)\;dS_g\\
&+\frac{C}{\lambda_i\delta^{2}}\oint_{B^{\hat{h}}_{\delta}(a_j)}\left(\frac{\lambda_j}{(1+\lambda^2_jd_g(a_j,x)^{2}}\right)\;dS_g+\frac{C}{\lambda_i\lambda_j\delta^{4}}\\
&\le\underbrace{C\int_{B_{\delta\lambda_i}}\frac{1}{1+|y|^2}\left(\frac{1}{\frac{\lambda_i}{\lambda_j}+\lambda_i\lambda_jG^{{\frac{-1}{1-\gamma}}}_{a_j}\left(\Psi_{a_i}\left(\frac{y}{\lambda_i}\right)\right)}\right)}_{\mbox{$I_1$}}\\
&+\frac{C}{\lambda_i\lambda_j\delta^{2}}\int_{B_{\delta\lambda_j}}\left(\frac{1}{1+|y|^2}\right)+\frac{C}{\lambda_i\lambda_j\delta^{4}}\\&
\le C I_1+\frac{C}{ \lambda_i\lambda_j \delta^{2}}\left[\log\left(\lambda_j\delta\right)+C\right]+\frac{C}{\lambda_i\lambda_j\delta^{4}}\\&
\le C I_1+\frac{C}{\lambda_i\lambda_j\delta^{4}}\log\left(\lambda_j\right).
\end{split}
\end{equation}
\noindent
Now, we estimate \;$I_1$\; as follows.
\vspace{4pt}

\noindent
If \;$\varepsilon^{{\frac{-1}{1-\gamma}}}_{ij}\sim\frac{\lambda_i}{\lambda_j},$\; then we get
\[I_1\le C \varepsilon^{\frac{1}{1-\gamma}}_{ij}\left[\log\left(\lambda_i\delta\right)+C\right].\]
So, for \;$I$\; we have
\begin{equation*}
\begin{split}
I&\le C \varepsilon^{\frac{1}{1-\gamma}}_{ij}\left[\log\left(\lambda_i\delta\right)+C\right]+\frac{C}{\lambda_i\lambda_j\delta^{4}}\log\left(\lambda_i\lambda_j\right)
\\&\le \frac{C}{\delta^{4}}\varepsilon^{{\frac{1}{1-\gamma}}}_{ij}\log\left(\varepsilon^{{\frac{-1}{1-\gamma}}}_{ij}G^{{\frac{1}{1-\gamma}}}_{a_i}(a_j)\right)\\
&=O\left(\varepsilon^{\frac{1}{1-\gamma}}_{ij} \delta^{-4}\log\left(\varepsilon^{\frac{-1}{2-2\gamma}}_{ij}\delta^{-1}\right)\right).
\end{split}
\end{equation*}

\noindent
If \;$\varepsilon^{{\frac{-1}{1-\gamma}}}_{ij}\sim \lambda_i\lambda_jG^{{\frac{-1}{1-\gamma}}}_{a_i}(a_j),$\; then we get
\begin{equation*}
I_1\le \frac{C}{\lambda_i\lambda_j\delta^{2}}\left[\log\left(\lambda_i\delta\right)+C\right].
\end{equation*}
\noindent
So, for \;$I$\; we have
\begin{equation*}
I\le\frac{C}{\lambda_i\lambda_j\delta^{4}}\left[\log\left(\lambda_i\lambda_j\right)\right].
\end{equation*}
This implies 
\begin{equation*}
I\le \frac{C}{\delta^{4}}\varepsilon^{{\frac{1}{1-\gamma}}}_{ij}\log\left(\varepsilon^{{\frac{-1}{1-\gamma}}}_{ij}G^{{\frac{1}{1-\gamma}}}_{a_i}(a_j)\right).
\end{equation*} 

\noindent
Hence, for \;$d_g(a_i,a_j)\geq 2\d,$\; we obtain
\begin{equation}\label{estfar}
I=O\left(\varepsilon^{{\frac{1}{1-\gamma}}}_{ij} \delta^{-4}\log\left(\varepsilon^{\frac{-1}{2-2\gamma}}_{ij}\delta^{-1}\right)\right).
\end{equation}

\noindent
On the other hand, arguing as above, if \;$d_g(a_i,a_j)<2\delta,$\; then we have also 
\begin{equation*}
\begin{split}
I&\le I_1+\frac{C}{\lambda_i\lambda_j\delta^{2}}\left[\log\left(\lambda_j\right)\right]+\frac{C}{\lambda_i\lambda_j\delta^{4}}\\&
\le I_1+\frac{C}{\lambda_i\lambda_j\delta^{4}}\log\left(\lambda_i\lambda_j\right),
\end{split}
\end{equation*}
where \;$I_1$\; is as in \eqref{i1}. Thus, if \;$\varepsilon_{i,j}^{{\frac{-1}{1-\gamma}}}\simeq \frac{\l_i}{\l_j},$\; then
\begin{equation*}
  \begin{split}
      I\le I_1+\frac{C}{\delta^{4}}\left(\frac{\lambda_j}{\lambda_i}\right)\frac{1}{\lambda_j}\left[\log(\frac{\lambda_i}{\lambda_j})+\log(\lambda^{2}_j)\right].
  \end{split}  
\end{equation*}
This implies
\begin{equation*}
    I\leq  I_1+\frac{C}{\d^{4}}\varepsilon_{ij}^{\frac{1}{1-\gamma}}\log \left(\varepsilon_{ij}^{\frac{-1}{2-2\gamma}}\right).
\end{equation*}
\vspace{5pt}
\noindent 
Next, if \;$\varepsilon_{i,j}^{{\frac{-1}{1-\gamma}}}\simeq \l_i\l_j G^{{\frac{-1}{1-\gamma}}}_{a_i}(a_j),$\; then we get
\begin{equation*}
\begin{split}
 I&\le I_1+\frac{1}{\lambda_i\lambda_j\delta^{4}G^{{\frac{-1}{1-\gamma}}}_{a_i}(a_j)}\left[\log\left(\lambda_i\lambda_jG^{{\frac{-1}{1-\gamma}}}_{a_i}(a_j)\right)+\log\left(G^{{\frac{1}{1-\gamma}}}_{a_i}(a_j)\right)\right]G^{{\frac{-1}{1-\gamma}}}_{a_i}(a_j)\\
 &\le I_1+\frac{C}{\delta^{4}}\varepsilon_{ij}^{\frac{1}{1-\gamma}}\log\left(\varepsilon_{ij}^{\frac{-1}{2-2\gamma}}\right).
\end{split}
\end{equation*}

\noindent
In order to arrive at an estimate for \;$I_1,$\;  we first begin by defining the sets as follows:
\begin{equation*}
\begin{split}
&A_1=\left\{y\in \R^{2}:\;|y|\le \epsilon\lambda_i \sqrt{G^{{\frac{-1}{1-\gamma}}}_{a_j}(a_i)+\frac{1}{\lambda^2_j}}\right\},\\&
A_2=\left\{y\in \R^{2}:\;\epsilon\lambda_i \sqrt{G^{{\frac{-1}{1-\gamma}}}_{a_j}(a_i)+\frac{1}{\lambda^2_j}}\le |y|\le E \lambda_i\sqrt{G^{{\frac{-1}{1-\gamma}}}_{a_j}(a_i)+\frac{1}{\lambda^2_j}}\right\},\\&
A_3=\left\{y\in \R^{2}:\;E\lambda_i \sqrt{G^{{\frac{-1}{1-\gamma}}}_{a_j}(a_i)+\frac{1}{\lambda^2_j}}\le |y|\le 4\lambda_i\delta\right\},
\end{split}
\end{equation*}
with \;$0<\epsilon<E<\infty$. Clearly by  the definition of \;$I_1$ (see \eqref{i1}), we have
\[I_1\le \int_{A_1}L_{ij}+\int_{A_2}L_{ij}+\int_{A_3}L_{ij},\]
where 
\[L_{ij}=\left(\frac{1}{1+|y|^2}\right)\left(\frac{1}{\frac{\lambda_i}{\lambda_j}+\lambda_i\lambda_jG^{{\frac{-1}{1-\gamma}}}_{a_j}\left(\Psi_{a_i}\left(\frac{y}{\lambda_i})\right)\right)}\right).\]

\noindent
Now, let us estimate for \;$\int_{A_1}L_{ij}.$\; We have 
\begin{equation*}\begin{split}
\int_{A_1}L_{ij} &\le C \varepsilon^{{\frac{1}{1-\gamma}}}_{ij}\int_{A_1}\left(\frac{1}{1+|y|^2}\right)\\
&\le C\varepsilon^{{\frac{1}{1-\gamma}}}_{ij} \log\left(\sqrt{\frac{\lambda_i}{\lambda_j}}\sqrt{\lambda_i\lambda_jG^{{\frac{-1}{1-\gamma}}}_{a_j}(a_i)+\frac{\lambda_i}{\lambda_j}}\right)\\&
\le C\varepsilon^{{\frac{1}{1-\gamma}}}_{ij} \log\left(\varepsilon^{\frac{-1}{2-2\gamma}}_{ij}\right).
\end{split}
\end{equation*}

\noindent
Next, for \;$\int_{A_2}L_{i,j},$\; we derive
\begin{equation*}
\begin{split}
\int_{A_2}L_{ij}&\le C \left(\frac{1}{\left(\frac{\lambda_i}{\lambda_j}\right)^2+\lambda^2_iG^{{\frac{-1}{1-\gamma}}}_{a_j}(a_i)}\right) \int_{A_2}\left(\frac{1}{\frac{\lambda_i}{\lambda_j}+\lambda_i\lambda_jG^{{\frac{-1}{1-\gamma}}}_{a_j}\left(\Psi_{a_i}\left(\frac{y}{\lambda_i}\right)\right)}\right)\\&
\le C \left(\frac{\lambda_j}{\lambda_i}\right) \varepsilon^{{\frac{1}{1-\gamma}}}_{ij}\int_{\{y\in\R^{2}:\;|y|\leq E\l_i\sqrt{G_{a_j}^{{\frac{-1}{1-\gamma}}}(a_i)+\frac{1}{\l_j^2}}\}}\left(\frac{1}{\frac{\lambda_i}{\lambda_j}+\frac{\lambda_j}{\lambda_i}\left|\lambda_i\Psi^{-1}_{a_j}\circ\Psi_{a_i}\left(\frac{y}{\lambda_i}\right)\right|^2}\right)\\
&\le C \left(\frac{\lambda_j}{\lambda_i}\right)\varepsilon^{{\frac{1}{1-\gamma}}}_{ij}\int_{\{y\in\R^{2}:\;|y|\leq \bar E\l_i\sqrt{G_{a_j}^{{\frac{-1}{1-\gamma}}}(a_i)+\frac{1}{\l_j^2}}\}}\left(\frac{1}{\frac{\lambda_i}{\lambda_j}+\frac{\lambda_j}{\lambda_i}\left|y\right|^2}\right)\\
&\le C \left(\frac{\lambda_j}{\lambda_i}\right)^{2}\left(\frac{\lambda_j}{\lambda_i}\right)^{-2}\varepsilon^{{\frac{1}{1-\gamma}}}_{ij}\int_{\{y\in\R^{2}:\;|y|\leq \bar E\l_j\sqrt{G_{a_j}^{{\frac{-1}{1-\gamma}}}(a_i)+\frac{1}{\l_j^2}}\}}\left(\frac{1}{1+|y|^2}\right)\\&
\le C\varepsilon^{{\frac{1}{1-\gamma}}}_{ij}\log\left(\varepsilon^{\frac{-1}{2-2\gamma}}_{ij}\right).
\end{split}
\end{equation*}

\noindent
Finally, we estimate for \;$\int_{A_3}L_{i,j}.$\; We get 
\begin{equation*}
\begin{split}
\int_{A_3}L_{ij}&\le\int_{A_3}\left(\frac{1}{1+|y|^2}\right)\left(\frac{1}{\frac{\lambda_i}{\lambda_j}+\frac{\lambda_j}{\lambda_i}|y|^2}\right)\\&
\le C \left(\frac{\lambda_i}{\lambda_j}\right)\int_{A_3}\frac{1}{|y|^4}\\&
\le C \left(\frac{\lambda_i}{\lambda_j}\right)\frac{1}{\left(\lambda^{2}_iG^{{\frac{-1}{1-\gamma}}}_{a_j}(a_i)+\left(\frac{\lambda_i}{\lambda_j}\right)^2\right)}\\&
\le C \varepsilon^{{\frac{1}{1-\gamma}}}_{ij}.
\end{split}
\end{equation*}
\noindent
Therefore, we have  
\begin{equation*}
    \begin{split}
      I_1\leq C\varepsilon^{{\frac{1}{1-\gamma}}}_{ij}\log{\left(\varepsilon^{\frac{-1}{2-2\gamma}}_{ij}\right)}.  
    \end{split}
\end{equation*}

\noindent
This implies for \;$d_g(a_i,a_j)<2\d,$\;  we have 
\begin{equation*}
    \begin{split}
       I=O\left(\frac{\varepsilon^{{\frac{1}{1-\gamma}}}_{ij}}{\delta^{4}} \log\left(\varepsilon^{\frac{-1}{2-2\gamma}}_{ij}\right)\right).
    \end{split}
\end{equation*}
 
\noindent
Hence, combining with the estimate for \;$d_g(a_i,a_j)\geq 2\d$\; (see \eqref{estfar}), we have 
\[\oint_{M} u^{\frac{1}{1-\gamma}}_{a_i,\lambda_i} u^{\frac{1}{1-\gamma}}_{a_j,\lambda_j}\;dS_g=O\left(\frac{\varepsilon^{{\frac{1}{1-\gamma}}}_{ij}}{\delta^{4}}\log\left(\varepsilon^{\frac{-1}{2-2\gamma}}_{ij}\delta^{-1}\right)\right).\]
\end{pf}\\

\noindent
\vspace{4pt}
Finally, we establish a sharp unbalanced high-order inter-action estimate that is required for the application of the Barycenter Technique of Bahri-Coron\cite{bc} for existence. 

\noindent
\begin{lem}\label{interact6}
Assuming that \;$\theta>0$\; is small and \;$\mu_0>0$\; is small, then \;$\forall$\;\;$ a_i, a_j\in M$,\;\;$\forall 0<2\delta<\delta_0,$\; and \;$\forall$\;\;$0<\frac{1}{\l_i}\le\frac{1}{\l_j}\leq \theta\delta$\; such that \;$\varepsilon_{ij}\leq \mu_0,$\; we have
\[\oint_{M} u^{\alpha}_{a_i,\lambda_i} u^{\beta}_{a_j,\lambda_j}\;dS_g=O\left(\frac{\varepsilon_{ij}^{\beta}}{\d^{4}}\right).\]
where \;$\delta_0$\; is as in \eqref{delta0}, \;\(\alpha+\beta=\frac{2}{1-\gamma},\)\; and \;\(\alpha>\frac{1}{1-\gamma}>\beta>0.\)
\end{lem}
\begin{pf}
Let \;$\hat{\alpha}=(1-\gamma)\alpha$\; and \;$\hat{\beta}=(1-\gamma)\beta.$\; Then we have \;$\hat{\alpha}+\hat{\beta}=2.$\; Now, since \;$\lambda_j\le \lambda_i,$\; then for \;$\varepsilon_{ij}\sim 0$\; we have\\
1) Either \;$\varepsilon^{\frac{-1}{1-\gamma}}_{ij}\sim \lambda_i\lambda_jG^{\frac{-1}{1-\gamma}}_{a_i}(a_j).$\;
\vspace{4pt}

\noindent
2) Or \;$\varepsilon^{\frac{-1}{1-\gamma}}_{ij}\sim\frac{\lambda_i}{\lambda_j}.$\;
\vspace{4pt}

\noindent
To continue, we write
\[\oint_{M} u^{\alpha}_{a_i,\lambda_i} u^{\beta}_{a_j,\lambda_j}\;dS_g=\underbrace{\oint_{B^{\hat{h}}_{\delta}(a_i)} u^{\alpha}_{a_i, \lambda_i} u^{\beta}_{a_j,\lambda_j}\;dS_g}_{\mbox{$I_1$}}+\underbrace{\oint_{{M}-B^{\hat{h}}_{\delta}(a_i)} u^{\alpha}_{a_i, \lambda_i} u^{\beta}_{a_j, \lambda_j}\;dS_g}_{\mbox{$I_2$}}\]
\noindent
and estimate \;$I_1$\; and \;$I_2.$\; For \;$I_2,$\; we have 
\begin{equation*}\begin{split}
I_2&= \oint_{\left({M}-B^{\hat{h}}_{\delta}(a_i)\right)\cap B^{\hat{h}}_{\delta}(a_j)} u^{\alpha}_{a_i,\lambda_i} u^{\beta}_{a_j,\lambda_j}\;dS_g+\oint_{{M}-\left(B^{\hat{h}}_{\delta}(a_i)\cup B^{\hat{h}}_{\delta}(a_j)\right)} u^{\alpha}_{a_i,\lambda_i} u^{\beta}_{a_j,\lambda_j}\;dS_g\\
&\le C\oint_{\left({M}-B^{\hat{h}}_{\delta}(a_i)\right)\cap B^{\hat{h}}_{\delta}(a_j)} \left(\frac{\lambda_i}{1+\lambda^2_iG^{\frac{-1}{1-\gamma}}_{a_i}(x)}\right)^{\hat{\alpha}}\left(\frac{\lambda_j}{1+\lambda^2_jd_g(a_j,x)^{2}}\right)^{\hat{\beta}}\;dS_g\\
&+C\oint_{{\partial M}-\left(B^{\hat{h}}_{\delta}(a_i)\cup B^{\hat{h}}_{\delta}(a_j)\right)}\left(\frac{\lambda_i}{1+\lambda^2_iG^{\frac{-1}{1-\gamma}}_{a_i}(x)}\right)^{\hat{\alpha}}\left(\frac{\lambda_j}{1+\lambda^2_jG^{\frac{-1}{1-\gamma}}_{a_j}(x)}\right)^{\hat{\beta}}\;dS_g\\
&\le \frac{C}{\lambda^{\hat{\alpha}}_i\lambda^{2-\hat{\beta}}_j\delta^{\alpha}}\int_{B_{\delta\lambda_j}}\left(\frac{1}{1+|y|^2}\right)^{\hat{\beta}}+\frac{C}{\lambda^{\hat{\alpha}}_i\lambda^{\hat{\beta}}_j\delta^{4}}\\&
\le \frac{C}{\lambda^{\hat{\alpha}}_i\lambda^{2-\hat{\beta}}_j\delta^{\alpha}}\left(\frac{1}{\lambda_j \delta}\right)^{2\hat{\beta}-2}+\frac{C}{\lambda^{\hat{\alpha}}_i\lambda^{\hat{\beta}}_j\delta^{4}}.
\end{split}
\end{equation*}

\noindent
Thus, we have for \;$I_2$\;
\begin{equation}\label{esti2int}
I_2\le \frac{C}{\lambda^{\hat{\alpha}}_i\lambda^{\hat{\beta}}_j\delta^{4}}.
\end{equation}
\noindent
Next, for \;$I_1$\; we have 
\begin{equation*}
\begin{split}
I_1&=\oint_{B^{\hat{h}}_{\delta}(a_i)}\left(\frac{\lambda_i}{(1+\lambda^2_id_g(a_i,x)^{2}}\right)^{\hat{\alpha}}\left(\frac{\lambda_j}{1+\lambda^2_jG^{\frac{-1}{1-\gamma}}_{a_j}(x)}\right)^{\hat{\beta}}\;dS_g\\
&= \int_{B_{\delta\lambda_i}} \left(\frac{1}{1+|y|^2}\right)^{\hat{\alpha}}\left[\frac{1}{\frac{\lambda_i}{\lambda_j}+\lambda_i\lambda_jG^{\frac{-1}{1-\gamma}}_{a_j}\left(\Psi_{a_i}\left(\frac{y}{\lambda_i}\right)\right)}\right]^{\hat{\beta}}.
\end{split}\end{equation*}
\noindent
Thus, if \;$\varepsilon^{\frac{-1}{1-\gamma}}_{ij}\sim \frac{\lambda_i}{\lambda_j},$\; then 
\begin{equation*}
\begin{split}
I_1&\le C \varepsilon^{\frac{\hat{\beta}}{1-\gamma}}_{ij}\left[\left(\frac{1}{\lambda_i\delta}\right)^{2\hat{\alpha}-2}+C\right]\\&\le C \varepsilon^{\beta}_{ij}.
\end{split}
\end{equation*}
\noindent
If \;$\varepsilon^{\frac{-1}{1-\gamma}}_{ij}\sim \lambda_i \lambda_j G^{\frac{-1}{1-\gamma}}_{a_i}(a_j)$\; and \;$d_g(a_i,a_j)\geq2\delta,$\; then we have 
\begin{equation*}
\begin{split}
I_1&\le C \left(\frac{1}{\lambda_i\lambda_j\delta^{2}}\right)^{\hat{\beta}}\left[\left(\frac{1}{\lambda_i\delta}\right)^{2\hat{\alpha}-2}+C\right]\\&
\le \frac{C}{\delta^{2}}\left(\frac{1}{\lambda_i\lambda_j}\right)^{\hat{\beta}}\le \frac{C}{\delta^{2}}\left[\left(\frac{1}{\lambda_i\lambda_jG^{\frac{-1}{1-\gamma}}_{a_i}(a_j)}\right)^{1-\gamma}\right]^{\beta}
\\&\le \frac{C}{\delta^{2}}\varepsilon^{\beta}_{ij}.
\end{split}
\end{equation*}
\noindent
Now, if \;$\varepsilon^{\frac{-1}{1-\gamma}}_{ij}\sim \lambda_i \lambda_j G^{\frac{-1}{1-\gamma}}_{a_i}(a_j)$\; and \;$d_g(a_i,a_j)<2\delta,$\; then we get
\[I_1\le C\int_{B_{\delta\lambda_i}} \left(\frac{1}{1+|y|^2}\right)^{\hat{\alpha}}\left[\frac{1}{\frac{\lambda_i}{\lambda_j}+\lambda_i\lambda_j\left|\Psi^{-1}_{a_j}\circ\Psi_{a_i}\left(\frac{y}{\lambda_i}\right)\right|^{2}}\right]^{\hat{\beta}}.\]
Next, we define 
\[B=\left\{y\in\R^{2}:\;\frac{1}{2}d_g(a_i,a_j)\le \frac{|y|}{\lambda_i}\le2d_g(a_i,a_j)\right\}\]
and have
\begin{equation*}
\begin{split}
I_1&\le C\int_{B}\left(\frac{1}{1+|y|^2}\right)^{\hat{\alpha}}\left[\frac{1}{\frac{\lambda_i}{\lambda_j}+\lambda_i\lambda_j\left|\Psi^{-1}_{a_j}\circ\Psi_{a_i}\left(\frac{y}{\lambda_i}\right)\right|^{2}}\right]^{\hat{\beta}}\\
&+C\int_{B_{\delta\lambda_i}-B} \left(\frac{1}{1+|y|^2}\right)^{\hat{\alpha}}\left[\frac{1}{\frac{\lambda_i}{\lambda_j}+\lambda_i\lambda_j\left|\Psi^{-1}_{a_j}\circ\Psi_{a_i}\left(\frac{y}{\lambda_i}\right)\right|^{2}}\right]^{\hat{\beta}}.
\end{split}
\end{equation*}
\noindent
For the second term, we have
\begin{equation*}
\begin{split}
\int_{B_{\delta\lambda_i}-B} \left(\frac{1}{1+|y|^2}\right)^{\hat{\alpha}}\left[\frac{1}{\frac{\lambda_i}{\lambda_j}+\lambda_i\lambda_j\left|\Psi^{-1}_{a_j}\circ\Psi_{a_i}\left(\frac{y}{\lambda_i}\right)\right|^{2}}\right]^{\hat{\beta}}
&\le C \varepsilon^{\beta}_{ij}\left[\left(\frac{1}{\lambda_i\delta}\right)^{\alpha(1-\gamma)-2}+C\right]\\&\le C \varepsilon^{\beta}_{ij}.
\end{split}
\end{equation*}
\noindent
For the first term, we have
\begin{equation*}
\begin{split}
&\int_{B}\left(\frac{1}{1+|y|^2}\right)^{\hat{\alpha}}\left[\frac{1}{\frac{\lambda_i}{\lambda_j}+\lambda_i\lambda_j\left|\Psi^{-1}_{a_j}\circ\Psi_{a_i}\left(\frac{y}{\lambda_i}\right)\right|^{2}}\right]^{\hat{\beta}}\\
&\le C\left(\frac{1}{1+\lambda^2_id^{2}_g(a_i,a_j)}\right)^{\hat{\alpha}}\int_{\{y\in\R^{2}:\;|y|\le 2\lambda_id_g(a_i,a_j)\}}\left[\frac{1}{\frac{\lambda_i}{\lambda_j}+\frac{\lambda_j}{\lambda_i}\left|\lambda_i\Psi^{-1}_{a_j}\circ\Psi_{a_i}\left(\frac{y}{\lambda_i}\right)\right|^2}\right]^{\hat{\beta}}\\
&\le C\left(\frac{1}{1+\lambda^2_id^{2}_g(a_i,a_j)}\right)^{\hat{\alpha}}\int_{\{z\in\R^{2}:\;|z|\le 4\lambda_id_g(a_i,a_j)\}}\left[\frac{1}{\frac{\lambda_i}{\lambda_j}+\frac{\lambda_j}{\lambda_i}|z|^2}\right]^{\hat{\beta}}\\
&\le C\left(\frac{1}{\frac{\lambda_j}{\lambda_i}+\lambda_i\lambda_jd^{2}_g(a_i,a_j)}\right)^{\alpha(1-\gamma)}\int_{\{z\in\R^{2}:\;|z|\le 4\lambda_jd_g(a_i,a_j)\}}\left[\frac{1}{1+|z|^2}\right]^{\hat{\beta}}.
\end{split}
\end{equation*}
\noindent
If \;$\lambda_jd_g(a_i,a_j)$\; is bounded, then we get
\begin{equation*}
\begin{split}
I_1&\le C\left(\frac{1}{\frac{\lambda_j}{\lambda_i}+\lambda_i\lambda_jd_g(a_i,a_j)^{2}}\right)^{\alpha(1-\gamma)}\\&\le C\varepsilon^{\beta}_{ij}. 
\end{split}
\end{equation*}
\noindent
If \;$\lambda_jd_g(a_i,a_j)$\; is unbounded, then we get
\begin{equation*}
\begin{split}
I_1&\le C\left(\frac{1}{\frac{\lambda_j}{\lambda_i}+\lambda_i\lambda_jd_g(a_i,a_j)^{2}}\right)^{\alpha(1-\gamma)}\left(\lambda_jd_g(a_i,a_j)^{2}\right)^{2-2\hat{\beta}}\\
&\le C \left(\frac{1}{1+\lambda^2_id_g(a_i,a_j)^{2}}\right)^{\hat{\alpha}+\hat{\beta}-1}\left(\frac{\lambda_i}{\lambda_j}\right)^{\hat{\beta}}\\
&\le C \left(\frac{1}{\frac{\lambda_j}{\lambda_i}+\lambda_i\lambda_jd_g(a_i,a_j)^{2}}\right)^{\hat{\beta}}\left(\frac{1}{1+\lambda^2_id_g(a_i,a_j)^{2}}\right)^{\hat{\alpha}-1}\\&\le C \varepsilon^{\beta}_{ij}.
\end{split}
\end{equation*}
\noindent
Thus, we have for \;$I_1$\;
\begin{equation}\label{esti1f}
I_1\le \frac{C}{\delta^{2}} \varepsilon^{\beta}_{ij}.
\end{equation}
\noindent
On the other hand, using the estimate for \;$I_2$\; (see \eqref{esti2int}), we have
\begin{equation}\label{esti2f}
I_2=O\left(\frac{\varepsilon^{\beta}_{ij}}{\delta^{4}}\right).
\end{equation}
Hence, combining \eqref{esti1f} and \eqref{esti2f}, we have 
\[\oint_{M} u^{\alpha}_{a_i,\lambda_i} u^{\beta}_{a_j,\lambda_j}\;dS_g=O\left(\frac{\varepsilon_{ij}^{\beta}}{\d^{4}}\right).\]
\end{pf}
\vspace{4pt}
%
%
%
%
%
\section{Algebraic topological argument}\label{ATA}
In this Section, we present the algebraic topological argument for existence. We start by fixing some notations from algebraic topology.  
For  a topological space \;$Z$\; and \;$Y$\; a subspace of \;$Z,$\;\; $H_*(Z, Y)$\; stands for the relative homology with \;$\Z_2$\; coefficients of the topological pair \;$(Z, Y)$. For\;$f: (Z, Y)\longrightarrow (W, X)$ \; a map with \;$(Z, Y)$\; and \;$(W, X)$\; topological pairs, \;$ f_*$ denotes the induced map in relative homology.
\vspace{4pt}

\noindent
Furthermore, we discuss some algebraic topological tools needed for our application of the Barycenter Technique of Bahri-Coron\cite{bc} for existence. We start with recalling the space of formal the  barycenter of \;$ M.$\; For \;$p\in \N^*$,\; the set of formal barycenters of \;$M$\; of order \;$p$\; is defined as 
 \begin{equation}\label{eq:barytop}
B_{p}(M)=\{\sum_{i=1}^{p}\alpha_i\d_{a_i}\;:\;a_i\in M, \;\alpha_i\geq 0,\;\; i=1,\cdots, p,\;\,\sum_{i=1}^{p}\alpha_i=1\},\;\;\text{and}\;\;B_0(M)=\emptyset,
\end{equation}
where \;$\delta_{a}$\; for \;$a\in M$\; is the Dirac measure at \;$a.$\; Since \;$dim(M)=2,$\; then we have the existence of \;$\Z_2$\; orientation classes  (see \cite{bc} and \cite{kk})
\begin{equation}\label{orientation_classes}
w_p\in H_{3p-1}(B_{p}(M), B_{p-1}(M)), \;\;\;\;\;p\in \N^*.
\end{equation}\\
Now to continue, we fix \;$\delta$\; small such that \;$0<2\d< \delta_0$ \;where \;$\delta_0$ \; is as in \eqref{delta0}. Moreover, we choose \;$\theta_0>0$\; and smalll. After this, we let \;$\l$\; varies such that \;$0<\frac{1}{\l}\leq \theta_0\delta$\; and associate for every \;$p\in \N^*$\; the map
$$
f_p(\l): B_p(M)\longrightarrow W^{1,2}_{y^{1-2\gamma},+}(X) 
$$
defined by the formula
$$
f_p(\l)(\sigma)=\sum_{i=1}^p\alpha_i u_{a_i,\l}, \;\;\;\;\sigma=\sum_{i
=1}^p\alpha_i\d_{a_i}\in B_p(M),
$$
where \;$u_{a_i, \l}$\; is as in \eqref{ual} with \;$a$\; replaced by \;$a_i$\;.\\
\vspace{6pt}

\noindent
As in Proposition 3.1 in  \cite{martndia2} and Proposition 6.3 in \cite{nss}, using Corollary \ref{sharpenergy}, Corollary \ref{interact2}, Corollary \ref{interact4}, Lemma \ref{interact5}, and  Lemma \ref{interact6}, we have the following multiple-bubble estimate.
\begin{pro}\label{eq:baryest}
There exist \;$\bar C_0>0$\; and \;$\bar c_0>0$\; such that for every \;$p\in \N^*$, $p\geq 2$ and every \;$0<\varepsilon\leq \varepsilon_0$, there exists \;$\l_p:=\l_p(\varepsilon)$ such that for every \;$\l\geq \l_p$ and for every $\sigma=\sum_{i=1}^p\alpha_i\delta_{a_i}\in B_p(M)$, we have
\begin{enumerate}
 \item
If \;$\sum_{i\neq j}\varepsilon_{i, j}> \varepsilon$\; or there exist \;$i_0\neq j_0$\; such that \;$\frac{\alpha_{i_0}}{\alpha_{j_0}}>\nu_0$, then
 $$
J_q(f_p(\l)(\sigma))\leq p^{\gamma}\mathcal{S}.
$$ 
 \item
If \;$\sum_{i\neq j}\varepsilon_{i, j}\leq \varepsilon$\; and for every \;$i\neq j$\; we have \;$\frac{\alpha_{i}}{\alpha_j}\leq\nu_0$, then
$$
J_q(f_p(\l)(\sigma))\leq p^{\gamma}\mathcal{S}\left(1+ \frac{\bar C_0}{\lambda^{2-2\gamma}}-\bar c_{0}\frac{(p-1)}{\lambda^{2-2\gamma}}\right),
$$
where \;$\mathcal{S}$\; is as in \eqref{S}, \;$\varepsilon_{ij}$\; is as in \eqref{varepij}, \;$\l_i=\l_j=\l,$\;\;$\varepsilon_0$\; is as in \eqref{varepsilon0}, and \;$\nu_0$\; is as in \eqref{nu0}.
\end{enumerate}

\end{pro}
\vspace{6pt}

\noindent
As in Lemma 4.2  in \cite{martndia2} and Lemma 6.4 in \cite{nss}, we have the selection map \;$s_1$ (see \eqref{eq:mini}), Lemma \ref{deform} and Corollary \ref{sharpenergy} imply the following topological result.
\begin{lem}\label{eq:nontrivialf1}
Assuming that \;$J_q$\; has no critical points, then there exists \;$\bar \l_1>0$\;  such that for every \;$\l\geq \bar \l_1$, we have
$$
f_1(\l)\; : \;(B_1(M),\; B_0( M))\longrightarrow (W_1, \;W_0)
$$
is well defined and satisfies
$$
(f_1(\l))_*(w_1)\neq 0\;\;\text{in}\;\;H_{3}(W_1, \;W_0).
$$
\end{lem}
\vspace{6pt}

\noindent
As in Lemma 4.3  in \cite{martndia2} and Lemma 6.5 in \cite{nss}, we have the selection map \;$s_p$ (see \eqref{eq:mini}), Lemma \ref{deform} and Proposition \ref{eq:baryest} imply the following recursive topological result.
\begin{lem}\label{eq:nontrivialrecursive}
Assuming that \;$J_q$\; has no critical points, then  there exists \;$\bar \l_p>0$\; such that for every \;$\l\geq \bar\l_p$, we have 
$$
f_{p+1}(\l): (B_{p+1}(M),\; B_{p}(M))\longrightarrow (W_{p+1}, \;W_{p})
$$
and 
$$
f_p(\l): (B_p(M), \;B_{p-1}(M))\longrightarrow (W_p, \; W_{p-1})
$$
are well defined and satisfy
$$(f_p(\l))_*(w_p)\neq 0\;\; \text{in}\;\; \;\;H_{3p-1}(W_p, \;W_{p-1})$$ implies
$$(f_{p+1}(\l))_*(w_{p+1})\neq 0\;\; \text{in} \;\;H_{3(p+1)-1}(W_{p+1}, \;W_{p}).$$
\end{lem}
\vspace{6pt}

\noindent
Finally, as in Corollary 3.3 in \cite{martndia2} and Lemma 6.6 in \cite{nss}, we clearly have that Proposition \ref{eq:baryest}  implies the following result. 
\begin{lem}\label{eq:largep}
Setting \;$$\bar p_{0}:=[1+\frac{\bar C_{0}}{\bar c_{0}} ]+2$$ with \;$\bar C_0$\; and \;$\bar c_0$\; as in Proposition \ref{eq:baryest} and recalling \eqref{dfwp}, we have  there exists \;$\hat \l_{\bar p_0}>0$\; such that \;$\forall\l\geq\hat \l_{\bar p_0}$, 
$$
f_{\bar p_{0}}(\l)(B_{\bar p_{0}}(M))\subset W_{{\bar p}_0-1}.
$$
\end{lem}

\noindent
\begin{pfn} {of Theorem \ref{thm1}} \\
As in \cite{martndia2} and \cite{nss}, the theorem follows by a contradiction argument from Lemma \ref{eq:nontrivialf1} - Lemma \ref{eq:largep}.
\end{pfn}

\section{Appendix}\label{APP}
In this Section, using the explicit expression of \;$\bar{\delta}_{a,\lambda}$\; (see \eqref{deltabar}) or Corollary 3.2 in \cite{mmc}, we have the following technical estimates.
\begin{lem}\label{Appendix}
Recalling the definition of \;$\bar{\delta}_{a,\lambda}$\; see \eqref{deltabar}, and setting \;$\zeta=(\bar{\zeta},\zeta_3)$\; with \;$\bar{\zeta}=(\zeta_1,\zeta_2),$\; we have on \;$\bar{\R}^{3}_{+}$\;
\end{lem}
 \begin{equation*}
     \begin{split}
         &\;\;\bar{\delta}_{0,\lambda}(\zeta)=O\left(\left(\frac{\lambda}{1+\lambda^{2}|\zeta|^{2}}\right)^{1-\gamma}\right),\\
         &\;\;\partial_{\zeta_{3}}\bar{\delta}_{0,\lambda}(\zeta)=O\left(\lambda^{\gamma} \zeta_{3}^{2\gamma-1}\left(\frac{\lambda}{1+\lambda^{2}|\zeta|^{2}}\right)\right),\\
         &\;\;\nabla_{\bar{\zeta}}\bar{\delta}_{0,\lambda}(\zeta)=O\left(\sqrt{\lambda}\left(\frac{\lambda}{1+\lambda^{2}|\zeta|^{2}}\right)^{3-2\gamma}\right),\\
         &\;\;\nabla^{2}_{\bar{\zeta}}\bar{\delta}_{0,\lambda}(\zeta)=O\left(\lambda \left(\frac{\lambda}{1+\lambda^{2}|\zeta|^{2}}\right)^{2-\gamma}\right),\;\;\text{for}\;\; \gamma \in (0,1).
     \end{split}
 \end{equation*}
 

\vspace{6pt}
 
 \noindent
\begin{thebibliography}{99}
\bibitem{aldawood} Aldawood. Mohammed and Ndiaye. Cheikh.B,{\em The Brezis-Nirenberg problem on non-contractible bounded domains of $\R^{3}$}, arXiv preprint arXiv:2202.05147 (2022).

\bibitem{aldawood1}Aldawood. Mohammed and Ndiaye. Cheikh.B,{\em Cherrier-Escobar problem for the elliptic Schroedinger-to-Neumann map,} arXiv preprint arXiv:2206.08838 (2022).


\bibitem{moa} Ahmedou. Mohameden Ould., {\em Conformal deformations of Riemannian metrics via “critical
point theory at infinity”: the conformally flat case with umbilic boundary},Contemporary Mathematics. 350 (2004),
1–18.


\bibitem{sma1} Almaraz. S{\'e}rgio., {\em An existence theorem of conformal scalar flat metrics on manifolds with boundary}, Pacific Journal of Mathematics 248 (2010) , no. 1, 1–22.

\bibitem{sma2} Almaraz. S{\'e}rgio., {\em Convergence of scalar-flat metrics on manifolds with boundary under a Yamabe-type flow }, J. Differential. Equations, 259 (2015), no 7, 2626–2694.

\bibitem{sma3} Almaraz. S{\'e}rgio., {\em The asymptotic behavior of Palais-Smale sequences on manifolds with boundary},  J. Differential Equations 251 (2011), no. 7, 1813–1840. 

\bibitem{bc} Bahri. Abbas and  Coron.  Jean-Michel., {\em On a nonlinear elliptic equation involving the critical Sobolev exponent: the effect of the topology of the domain}, Comm. Pure Appl. Math. (1988), 41-3 , 253-294.

\bibitem{bb} Bahri. Abbas and Brezis. H., {\em Non-Linear Elliptic Equations on closed Riemannian Manifolds with the Sobolev Critical Exponent},Topics in Geometry: In Memory of Joseph D’Atri,(1996), 1-100.

\bibitem{bre1} Brendle. Simon., {\em Convergence of the Yamabe flow for arbitrary initial energy}, J. Diff. Geom. 69 (2005), 217-278.

\bibitem{brenir} Br{\'e}zis. H and Nirenberg. Louis., {\em Positive solutions of nonlinear elliptic equations involving critical Sobolev exponents},Communications on pure and applied mathematics, 36 1983), no. 4, 437–477.

\bibitem{brezis2}  Br{\'e}zis. H., {\em Elliptic equations with limiting Sobolev exponent--the impact of topology},Communications on pure
and applied mathematics 39 (1986), no. S1, S17–S39.


\bibitem{cgg}Chang. Sun-Yung Alice and Gonz\'alez.M. d. M., {\em Fractional Laplacian in conformal geometry}, Adv.Math 226
(2011), no. 2, 1410–1432.

\bibitem{pc} Cherrier. Pascal., {\em Problemes de Neumann non lin{\'e}aires sur les vari{\'e}t{\'e}s riemanniennes}, 
Journal of Functional
Analysis 57 (1984), no. 2, 154–206.

\bibitem{jfe1} Escobar. Jos{\'e} F., {\em Conformal deformation of a Riemannian metric to a scalar flat metric with constant mean curvature on the boundary}, Annals of Mathematics (1992), 1–50.

\bibitem{jfe2} Escobar. Jos{\'e} F., {\em Conformal deformation of a Riemannian metric to a constant scalar curvature metric with constant mean curvature on the boundary}, Indiana University Mathematics Journal (1996), 917–943.

\bibitem{yfmd}Fang, Yi and Gonz\'alez.M. d. M., {\em Asymptotic behavior of Palais-Smale sequences associated
with fractional Yamabe-type equations},Pacific Journal of Mathematics 278 (2015), no. 2, 369–405.

\bibitem{sk}Kim. Seunghyeok, Musso. Monica, and Wei. Juncheng., {\em Existence theorems of the fractional
yamabe problem}, Analysis \& PDE 11(2017), no. 1, 75–113.


\bibitem{fcm1} Marques. Fernando C., {\em Existence results for the Yamabe problem on manifolds with
boundary}, Indiana University
mathematics journal (2005), 1599–1620.

\bibitem{fcm2} Marques. Fernando C., {\em Conformal deformations to scalar-flat metrics with constant mean curvature on the boundary}, Communications in Analysis and Geometry 15 (2007), no. 2, 381–405.


\bibitem{gm}M. d. M. Gonz\'alez, Mazzeo. Rafe, and Sire. Yannick., {\em  Singular solutions of fractional order Laplacians Journal of Geometric Analysis}, Journal of Geometric Analysis 22 (2012), 845–863.
 
 \bibitem{gg}M. d. M. Gonz\'alez and Qing. Jie., {\em Fractional conformal Laplacians and fractional Yamabe problems}, Analysis
\& PDE 6 (2013), no. 7, 1535–1576.

\bibitem{martndia2} Mayer. Martin and Ndiaye. Cheikh.B, {\em Barycenter technique and the Riemann mapping problem of Cherrier-Escobar}, Journal of Differential Geometry 107 (2017), no. 3, 519–560.

\bibitem{mmc}Mayer. Martin and Ndiaye. Cheikh.B, {\em Fractional Yamabe problem on locally flat conformal infinities of Poincare-Einstein manifolds},arXiv preprint arXiv:1701.05919 (2017).


\bibitem{mmn} Mayer. Martin and Ndiaye. Cheikh.B,{\em symptotics of the Poisson kernel and Green's functions of the fractional conformal Laplacian}, Discrete and Continuous Dynamical Systems 42 (2022), no. 10, 5037-5062.


\bibitem{mr}Mazzeo. Rafe., {\em The Hodge cohomology of a conformally compact metric}, Journal of differential geometry 28 (1988), no. 2, 309–339.


\bibitem{nss} Ndiaye. Cheikh.B, Sire. Yannick., and Sun. Liming., {\em Uniformizations Theorems: Between Yamabe and Paneitz}, Pacific Journal of Mathematics 314 (2021), no. 1, 115–159.

\bibitem{gpap} Palatucci. Giampiero. and Pisante. Adriano., {\em A Global Compactness type result for Palais-Smale
sequences in fractional Sobolev spaces}, Nonlinear Analysis: Theory, Methods \& Applications 117 (2015), 1–7.



\bibitem{kk} Sadok. Kallel. and R. Karoui., {\em Symmetric joins and weighted barycenters}, Advanced Nonlinear Studies 11 (2006), 117 - 143.



\bibitem{sc} Schoen. Richard., {\em  Conformal deformation of a Riemannian metric to constant scalar curvature}, Journal of Differential Geometry 20 (1984), no. 2, 479–495.


\bibitem{sma} Struwe. Michael., {\em A global compactness result for elliptic boundary value problems involving limiting nonlinearities}, Mathematische zeitschrift 187 (1984), 511–517.












\end{thebibliography}
\end{document}